\documentclass[10 pt,letter,reqno]{amsart}
\usepackage{amssymb}
\usepackage[usenames,dvipsnames]{color}
\usepackage{hyperref}
\hypersetup{
    colorlinks=true,       			
    linkcolor=blue,          			
    citecolor=blue,		 		
    filecolor=blue,      				
    urlcolor=blue           			
}

\newtheorem*{thm}{Theorem}
\newtheorem*{prop}{Proposition}
\newtheorem*{lemma}{Lemma}
\newtheorem*{corollary}{Corollary}

\newenvironment{pf}{\paragraph{{\sc Proof}}}{\qed\par\medskip}
\theoremstyle{definition}
\newtheorem*{defn}{Definition}
\newtheorem*{rem}{Remark}

\numberwithin{equation}{section}

\newcommand {\IC}{\mathbb{C}}
\newcommand {\IN}{\mathbb{N}}                          
\newcommand {\IP}{\mathbb{P}}
\newcommand {\IR}{\mathbb{R}}

\newcommand {\IZ}{\mathbb{Z}}

\newcommand {\bfa}{\mathbf a}

\newcommand {\bfI}{\mathbf I}
\newcommand {\bfJ}{\mathbf J}

\newcommand {\A}{\mathcal A}
\newcommand {\B}{\mathcal B}
\newcommand {\C}{\mathcal C}
\newcommand {\D}{\mathcal D}
\newcommand {\E}{\mathcal E}
\newcommand {\F}{\mathcal F}
\newcommand {\G}{\mathcal G}
\renewcommand {\H}{\mathcal H}
\newcommand {\I}{\mathcal I}
\newcommand {\J}{\mathcal J}
\newcommand {\K}{\mathcal K}

\newcommand {\N}{\mathcal N}
\renewcommand {\O}{\mathcal O}

\newcommand {\Q}{\mathcal Q}
\newcommand {\R}{\mathcal R}

\newcommand {\U}{\mathcal U}
\newcommand {\V}{\mathcal V}

\newcommand {\Z}{\mathcal Z}

\renewcommand {\a}{\mathfrak a}

\newcommand {\g}{\mathfrak{g}}
\newcommand {\h}{\mathfrak h}
\renewcommand {\ll}{\mathfrak l}

\newcommand {\gX}{\mathfrak X}

\newcommand {\sfh}{\mathsf{h}}
\newcommand {\nablah}{\mathsf{h}}
\newcommand {\sfk}{\mathsf{k}}
\newcommand {\sfA}{\mathsf{A}}

\newcommand {\sfPhi}{\mathsf{\Phi}}

\newcommand {\ad}{\operatorname{ad}}
\newcommand {\Ad}{\operatorname{Ad}}

\newcommand {\Aut}{\operatorname{Aut}}
\newcommand {\End}{\operatorname{End}}
\newcommand {\Hom}{\operatorname{Hom}}
\renewcommand {\Im}{\operatorname{Im}}
\newcommand {\Ker}{\operatorname{Ker}}
\newcommand {\Rep}{\operatorname{Rep}}

\newcommand {\id}{\operatorname{id}}

\renewcommand {\sl}[1]{\mathfrak{sl}_{#1}}

\newcommand {\Ug}{U\g}

\newcommand {\Uhg}{U_{\hbar}\g}

\newcommand {\Uhsl}[1]{U_{\hbar}\sl{#1}}

\newcommand {\calpha}{C_{\alpha}}
\newcommand {\nablac}{\nabla_{C}}

\newcommand {\nablak}{\nabla_{\kappa}}
\newcommand {\wtnablak}{\wt{\nabla}_{\kappa}}

\newcommand {\reg}{_{\operatorname{reg}}}
\newcommand {\hreg}{\h\reg}

\newcommand {\ie}{{\it i.e., }}
\newcommand {\eg}{{\it e.g.}, }
\newcommand {\fd}{finite--dimensional }

\newcommand {\lhs}{left--hand side }
\newcommand {\rhs}{right--hand side }

\newcommand {\wrt}{with respect to }

\newcommand {\ol}{\overline}
\newcommand {\wt}{\widetilde}
\newcommand {\wh}{\widehat}

\newcommand {\mns}{maximal nested set }
\newcommand {\mnss}{maximal nested sets }

\newcommand {\supp}{\operatorname{supp}}
\newcommand {\zsupp}{\operatorname{\mathfrak{z}supp}}

\newcommand {\DCP}{De Concini--Procesi }
\renewcommand {\DJ}{Drinfeld--Jimbo }

\newcommand {\EK}{Etingof--Kazhdan }

\newcommand {\Ho}{Hochschild }
\newcommand {\KD}{Drinfeld--Kohno }
\newcommand {\KM}{Kac--Moody }

\newcommand {\KZ}{Knizhnik--Zamolodchikov  }

\newcommand {\DKKZ}{_
{\scriptscriptstyle{\operatorname{DKZ}}}}
\newcommand {\KKZ}{_
{\scriptscriptstyle{\operatorname{KZ}}}}

\newcommand {\qW}{quantum Weyl group }

\newcommand {\fml}{[\negthinspace[\hbar]\negthinspace]}

\newcommand {\bin}[3]{\begin{bmatrix}#1\\#2\end{bmatrix}_{#3}}
\newcommand {\half}[1]{\frac{#1}{2}}

\newcommand {\cor}[1]{\alpha_{#1}^{\vee}}	
\renewcommand {\root}[1]{\alpha_{#1}}
\newcommand {\cow}[1]{\lambda^{\vee}_{#1}}

\newcommand {\aand}{\qquad\text{and}\qquad}

\newcommand {\qc}{quasi--Coxeter }
\newcommand {\qca}{quasi--Coxeter algebra }
\newcommand {\qcas}{quasi--Coxeter algebras }
\newcommand {\ba}{bialgebra }

\newcommand {\qba}{quasibialgebra }
\newcommand {\qbas}{quasibialgebras }
\newcommand {\qtqba}{quasitriangular quasibialgebra }
\newcommand {\qtqbas}{quasitriangular quasibialgebras }

\newcommand {\qt}{quasitriangular }

\newcommand {\qcqtqba}{quasi--Coxeter quasitriangular quasibialgebra }
\newcommand {\qcqtqbas}{quasi--Coxeter quasitriangular quasibialgebras }

\newcommand {\Dg}{D_{\g}}

\newcommand {\gD}{\g_{D}}
\newcommand {\gB}{\g_B}

\newcommand {\BD}{B_D}
\newcommand {\ND}{\N_D}

\newcommand {\dH}{d_{H}}

\newcommand {\rD}{r_D}
\newcommand {\rB}{r_B}

\newcommand {\Bh}{_{B,\hbar}}

\newcommand {\DKZ}{_{D}^{\KKZ}}
\newcommand {\BKZ}{_B^{\KKZ}}
\newcommand {\ih}{_i^\hbar}

\newcommand {\ic}{_{i,\scriptscriptstyle{C}}}

\newcommand {\ul}[1]{\underline{#1}}

\newcommand {\Alt}{\operatorname{Alt}}

\newcommand {\veps}{\varepsilon}

\newcommand {\op}{^{\scriptscriptstyle{\operatorname{op}}}}

\newcommand {\VV}{\mathbb{V}}

\newcommand {\dprime}{^{\prime\prime}}

\newcommand {\Veck}{\mathbf{Vec_\sfk}}
\newcommand {\Mod}[1]{\mathbf{Mod_{fd}}(#1)}
\newcommand {\bfF}{\mathbf{F}}

\newcommand {\comment}[1]{\footnote{\textcolor{blue}{#1}}}
\newcommand {\Omit}[1]{}
\newcommand {\rlim}[1]{{r\negthinspace\lim_{#1}}}
\newcommand {\isom}{\stackrel{\sim}{\longrightarrow}}

\newcommand {\FMTV}{Felder--Markov--Tarasov--Varchenko }
\newcommand {\ds}[1]{\displaystyle{#1}}

\newcommand {\Kalpha}{\mathcal K_\alpha}
\newcommand {\Kalphap}{\mathcal K_{\alpha,+}}
\newcommand {\Kalpham}{\mathcal K_{\alpha,-}}
\newcommand {\Kalphapm}{\mathcal K_{\alpha,\pm}}

\newcommand {\ii}{^{(i)}}
\newcommand {\jj}{^{(j)}}
\newcommand {\kk}{^{(k)}}
\newcommand {\nn}{^{(n)}}

\newcommand {\oo}{^{(1)}}
\newcommand {\twtw}{^{(2)}}
\newcommand {\thth}{^{(3)}}
\newcommand {\IH}{\mathbb H}
\newcommand {\uzeta}{{\ul{\zeta}}}
\newcommand {\ugamma}{{\ul{\gamma}}}

\newcommand {\iCn}{\imath\C_n}
\newcommand {\iC}{\imath\C}
\renewcommand {\Re}{\operatorname{Re}}

\newcommand {\Left}[1]{#1^\ell}

\newcommand {\ulz}{\ul{z}}
\newcommand {\dlog}[1]{\frac{d#1}{#1}}

\newcommand {\Cno}{\IC^n_0}
\newcommand {\bigp}{(x_1\cdots x_n)}

\newcommand {\EE}{\End_{\IC\fml}(U\g^{\otimes n}\fml_o)}

\newcommand {\aBC}{_{\cdot\infty(\cdot\cdot)}}
\newcommand {\ABc}{_{(\cdot\cdot)\infty\cdot}}

\newcommand {\err}{r}
\newcommand {\JJ}{\ol{J}}

\newcommand {\Ugo}[1]{U\g^{\otimes #1}\fml^o}

\newcommand {\cowo}[1]{{\cow{#1}}^{(1)}}

\newcommand {\onetwo}[1]{{#1}^{(1)}+{#1}^{(2)}}
\newcommand {\oll}{\ol{\mathfrak{l}}}
\newcommand {\olg}{{\ol{\mathfrak{g}}}}
\newcommand {\olh}{{\ol{\mathfrak{h}}}}
\newcommand {\olC}{\ol{\C}}

\newcommand {\olmu}{{\ol{\mu}}}
\newcommand {\olK}{\ol{\K}}
\newcommand {\olF}{\ol{F}}
\newcommand {\olcalF}{\ol{\F}}
\newcommand {\sfQ}{{\mathsf Q}}

\newcommand {\muone}{\mu^{(1)}}
\newcommand {\cowone}[1]{{\cow{#1}}^{(1)}}
\newcommand {\olPhi}{\ol{\Phi}}
\newcommand {\dloga}{\frac{d\alpha}{\alpha}}

\newcommand {\Qnablac}{\Q^\nabla_C}
\newcommand {\Qnablak}{\Q^\nabla_\kappa}
\newcommand {\Qhc}{\Q^\hbar_C}
\newcommand {\Qhk}{\Q^\hbar_\kappa}

\newcommand {\Sikh}[1]{S_{#1,\kappa}^\hbar}
\newcommand {\Sich}[1]{S_{#1,C}^\hbar}
\newcommand {\Siknabla}[1]{S_{#1,\kappa}^\nabla}
\newcommand {\Sicnabla}[1]{S_{#1,C}^\nabla}

\newcommand {\QQch}{\Q\Q_C^\hbar}
\newcommand {\QQdh}{\Q\Q_C^\nabla}

\newcommand {\ot}[1]{^{\otimes #1}}
\newcommand {\Amod}[1]{\operatorname{Mod}(#1)}
\newcommand {\mon}{x_1\cdots x_n}
\newcommand {\Oint}{{\mathcal O}^{\scriptstyle{\operatorname{int}}}}

\begin{document}

\renewcommand {\comment}[1]{}

\title[Quasi--Coxeter algebras and the Casimir connection]
{Quasi--Coxeter quasitriangular quasibialgebras and the Casimir connection}
\author[V. Toledano Laredo]{Valerio Toledano Laredo}
\address{Department of Mathematics,
Northeastern University,
360 Huntington Avenue,
Boston MA 02115.}
\email{V.ToledanoLaredo@neu.edu}
\thanks{Work supported in part by NSF grants DMS--0707212, DMS-0854792 and DMS--1206305}
\date{January 2016}
\begin{abstract}
Let $\g$ be a complex, semisimple Lie algebra. We prove the existence of a \qcqtqba
structure on the enveloping algebra of $\g$, which binds the \qca structure underlying
the Casimir connection of $\g$ and the \qtqba one underlying its KZ equations. This
implies in particular that the monodromy of the rational Casimir connection of $\g$ is
described by the \qW operators of the quantum group $\Uhg$. \Omit{Along the way, we give
a canonical transcendental construction of twists killing Drinfeld's KZ associator using
an ODEs with regular singularities at $0,-1$ and an irregular singularity at $\infty$.}
\end{abstract}
\Omit{
arXiv abstract:
Let g be a complex, semisimple Lie algebra. We prove the existence of a quasi-Coxeter,
quasitriangular quasibialgebra structure on the enveloping algebra of g, which binds the
quasi-Coxeter structure underlying the Casimir connection of g and the quasitriangular
quasibialgebra one underlying its KZ equations. This implies in particular that the monodromy
of the rational Casimir connection of g is described by the quantum Weyl group operators
of the quantum group U_h(g).
}
\maketitle

\setcounter{tocdepth}{1}
\tableofcontents

\section{Introduction}


\subsection{} 

Let $\g$ be a complex, semisimple Lie algebra. De Concini \cite{DC}, and
independently the author \cite{TL1}, conjectured that the monodromy of the
Casimir connection of $\g$ is described by the \qW group operators of the
quantum group $\Uhg$, in a way analogous to the \KD theorem \cite{Dr4}.

This conjecture was proved in \cite{TL1} for the Lie algebra $\sl{n}$. For an
arbitrary $\g$, it was reduced in \cite{TL2} to a structural statement about
the enveloping algebra $U\g$, namely the existence of a \qcqtqba structure
on $U\g\fml$ binding the \qc structure underlying the 
Casimir connection of $\g$, to the \qtqba structures underlying the 
\KZ connections of its standard Levi subalgebras. 

The goal of the present paper is to establish the existence of such a structure,
and therefore prove the monodromy conjecture for any semisimple Lie algebra.

\subsection{} 

Let $\h\subset\g$ be a Cartan subalgebra, $\sfPhi\subset\h^*$
the corresponding  root system, and $\hreg=\h\setminus\bigcup
_{\alpha\in\sfPhi}\Ker(\alpha)$ the set of regular elements in $\h$.
Fix an non--degenerate, invariant bilinear form $(\cdot,\cdot)$ on
$\g$. The Casimir connection of $\g$ is a connection on the
holomorphically trivial bundle $\VV$ on $\hreg$ with fibre a
given \fd representation $V$ of $\g$. It is given by
\[\nablac=
d-\frac{\nablah}{2}
\sum_{\alpha\in\sfPhi_+}\frac{d\alpha}{\alpha}\cdot\calpha\]
where $\nablah\in\IC$ is a deformation parameter, $\alpha$
ranges over a chosen system of positive roots $\sfPhi_+$,\footnote
{the connection is independent of the choice of $\sfPhi_+$}
and $\calpha$ is the Casimir operator induced by the restriction
of $(\cdot,\cdot)$ to the three--dimensional subalgebra $\sl{2}
^\alpha\subseteq\g$ corresponding to the root $\alpha$.
The connection $\nablac$ is flat for any $\nablah$ \cite{
MTL,TL1,DC,FMTV}, and can be made equivariant \wrt
the Weyl group $W$ of $\g$ \cite{MTL,TL1}. Its monodromy
defines a one--parameter family of actions $\mu_\sfh$ of
the braid group $B_W=\pi(\hreg/W)$ on $V$ depending on
$\sfh$, which deforms the action of (the Tits extension of)
$W$. We denote by $\mu:B_W\to G(V\fml)$ the formal
Taylor series of $\mu_\sfh$ \wrt the parameter $\hbar=
\pi\iota\sfh$.

\subsection{} 

Let $\Uhg$ be the Drinfeld--Jimbo quantum group corresponding
to $\g$, thought of as a topological Hopf algebra over $\IC\fml$.
Let $\V$ be a quantum deformation of $V$, that is a topologically
free $\IC\fml$--module such that $\V/\hbar\V$ is isomorphic to
$V$ as a $\g$--module. Since $\V$ is integrable, the braid
group $B_W$ acts on $\V$ though the \qW operators of $\Uhg$
\cite{Lu}. The main result of the present paper is the following.

\begin{thm}\label{th:intro monodromy thm}
The monodromy $\mu:B_W\to GL(V\fml)$ of the Casimir connection
on $V$ is equivalent to the \qW action of the braid group $B_W$
on $\V$.
\end{thm}

\subsection{} 

Recall that a \qtqba is an algebra $A$ over a commutative ring $\sfk$
endowed with two morphisms, the coproduct $\Delta:A\to A\ot{2}$ and
counit $\veps:A\to \sfk$, and two distinguished invertible elements, the
$R$--matrix $R\in A\ot{2}$ and associator $\Phi\in A\ot{3}$ \cite{Dr3}.
The relations satisfied by $\Delta,\veps,R$ and $\Phi$ are designed so as
to endow the category of $A$--modules with the structure of a braided
tensor category. In particular, for any $V\in\Amod{A}$, there is a family
of actions $\rho_b:B_n\to GL(V\ot{n})$ of the $n$--strand braid group
on the $n$--fold tensor product of $V$, which are labelled by the choice
of a complete bracketing $b$ of the non--associative monomial
$x_1\cdots x_n$. The actions corresponding to different choices of 
$b$ are canonically isomorphic, via intertwiners built out of the
associator $\Phi$. For example, for $n=3$, the action of $\Phi$ on
$V\ot{3}$ intertwines $\rho_{((x_1x_2)x_3)}$ and $\rho_{(x_1(x_2x_3))}$.

\subsection{} 

In a similar spirit, a \qca $A$ is designed so that a module over
it carries a family of canonically equivalent representations of the
braid group $B_W$ of a given irreducible Coxeter group $W$ \cite
{TL2}. Central to this notion are the {\it \mnss}on the Coxeter graph
$D$ of $W$, which generalise complete bracketings to an arbitrary
Coxeter type. These were introduced by \DCP \cite{DCP1,DCP2},
and consist of maximal collections of connected subgraphs of $D$
which are pairwise compatible, that is such that either one is contained
in the other, or they are orthogonal, in the sense that they have
disjoint vertex sets and are not linked by an edge of $D$. When $W$ is of type
$\sfA_{n-1}$, with the vertices of $D$ labelled $1,\ldots,n-1$ as
in \cite[Planche 1]{Bo}, mapping a connected subdiagram $B$
with vertices $\{i,\ldots,j-1\}$ to the pair of parentheses $x_1\cdots
x_{i-1}(x_i\cdots x_j) x_{j+1}\cdots x_n$, and noting that $B,B'
\subseteq D$ are compatible precisely when the corresponding
pairs of parentheses are consistent, yields a bijection between
\mnss on $D$ and complete bracketings of the monomial $x_1
\cdots x_n$.

\subsection{} 

A \qca is endowed with three sets of data.
First, a collection of subalgebras $A_B$ labelled by the connected
subgraphs $B\subseteq D$, which are such that $A_{B_1}\subseteq
A_{B_2}$ if $B_1\subseteq B_2$, and $[A_{B_1},A_{B_2}]=0$
if $B_1$ and $B_2$ are orthogonal. Next, invertible elements
$\Phi_{\G\F}\in A$ called associators. These are labelled by
pairs of \mnss on $D$, and satisfy in particular the transitivity
relations $\Phi_{\H\G}\Phi_{\G\F}=\Phi_{\H\F}$. Finally, invertible
elements $S_i$ labelled by the vertices of $D$, and called local
monodromies. This data satisfies various compatibility relations,
in particular a version of the braid relations defining $B_W$.
They give rise to a family of actions $\lambda_\F:B_W\to GL
(V)$ on any $V\in\Amod{A}$, which are labelled by the \mnss
on $D$, with $\Phi_{\G\F}\in A$ intertwining $\lambda_\F$ and
$\lambda_\G$.


\subsection{} 

A \qcqtqba of type $W$ is a \qc algebra $A$ additionally
endowed with a coproduct $\Delta$ and counit $\veps$,
such that each diagrammatic subalgebra $A_B$ has a
\qtqba structure of the form $(\Delta,\veps,R_B,\Phi_B)$,
for a given $R$--matrix $R_B\in A_B\ot{2}$ and associator
$\Phi_B\in A_B\ot{3}$. This gives rise in particular to
a family of commuting representations of the groups
$B_n,B_W$ on the tensor power $V^{\otimes n}$ of
any $V\in\Amod{A}$, specifically
\[\rho_{B,b}:B_n\to GL(V)\aand \lambda_{\F,b}:B_W\to GL(V)\]
The first family is determined by the choice of a subdiagram
$B\subseteq D$ and a bracketing $b$ of the monomial
$\mon$, and arises by restricting $V$ to the \qtqba $A_B$, with the
representations $\pi_{B,b}$ and $\pi_{B',b'}$ equivalent
provided $B=B'$. The second arises from the action of
$A$ on $V\ot{n}$ determined by the choice of $b$, and
depends on the choice of a \mns $\F$ on $D$, with
$\lambda_{b,\F}$ equivalent to $\lambda_{b',\G}$ for
any $b,b'$ and $\F,\G$.

\subsection{} 

A \qcqtqba $A$ possesses an additional piece of data, which
binds the associators $\Phi_B$ coming from the \qtqba structure
on each diagrammatic $A_B$, to the associators $\Phi_{\G\F}$
coming from the \qc structure on $A$. This welding of the \qc
and \qtqba structures is what gives the examples of interest
the rigidity required to determine the monodromy of the Casimir
connection.

The additional data consists of {\it relative twists}, which are
elements $F_{(B,\alpha)}\in A_B\ot{2}$ labelled by a connected
subdiagram $B$ and a vertex $\alpha\in B$. These twists are
required to satisfy two identities. The first is the twist equation
\begin{equation}\label{eq:Phi F_B}
(\Phi_B)_{F_{(B;\alpha)}}=\Phi_{B\setminus\alpha}
\end{equation}
together with the requirement that $F_{(B;\alpha)}$ commute
with $\Delta(A_{B\setminus\alpha})$.\footnote
{The twist $\Phi_F$ of an associator $\Phi$ by a twist $F$ is
equal to $1\otimes F\cdot \id\otimes\Delta(F)\cdot\Phi\cdot
\Delta\otimes\id(F)^{-1}\cdot F^{-1}\otimes 1$. If $B\setminus
\alpha$ is not connected, the associator $\Phi_{B\setminus\alpha}$
is taken to be the (commuting) product $\prod_i \Phi_{B_i}$, where
$B_i$ runs over the connected components of $B\setminus
\alpha$.} This amounts to asking that $F_{(B;\alpha)}$ defines
a tensor structure on the restriction functor from the monoidal
category of $A_B$--modules with associativity constraints given
by the associator $\Phi_B$, to that of $A_{B\setminus\alpha}
$--modules with associativity constraints given by $\Phi_{B
\setminus\alpha}$.

By restricting in stages from $A_D$ to $A_\emptyset=\sfk$,
the relative twists allow to define a tensor structure on the\
forgetful functor from the monoidal category of $A$--modules
with associativity constraints given by $\Phi_D$ to $\sfk$--modules,
which depends on the choice of a \mns $\F$ on $D$, as follows.
For any element $B$ of $\F$, the collection of maximal elements
of $\F$ properly contained in $B$ covers all but one of the vertices
$\alpha^B_\F$ of $B$ (and consists in fact of the connected
components of $B\setminus\alpha$
\footnote{In type $\sfA_{n-1}$, if $B$ is the diagram with vertices
$i,\ldots,j-1$, the elements of $\F$ properly contained in $B$ give
a bracketing of the monomial $x_i\cdots x_j$, which necessarily
contains two maximal pairs of parentheses of the form
$(x_i\cdots x_{k-1})(x_k\cdots x_j)$, and $\alpha^B_\F$ is
the vertex $k$.}). Define the twist $F_F\in A\ot{2}$ by
\begin{equation}\label{eq:factorised F}
\ds{F_\F=\stackrel{\longrightarrow}{\prod}_{B\in\F} F_{(B;\alpha^B_\F)}}
\end{equation}
where the product is taken with $F_{(B,\alpha)}$ to the
right of $F_{(B',\alpha')}$ if $B'\subset B$.\footnote{this
does not specify the order of the factors uniquely, but any
two orders satisfying the above requirements are easily
seen to give rise to the same product.} Then, it follows
from \eqref{eq:Phi F_B} that $(\Phi_D)_{F_\F}=1^{\otimes 3}$. 

The second identity satisfied by the relative twists requires
that the tensor structures $F_\F$ corresponding to the choices
of different \mnss be isomorphic, with the isomorphism given
by the associators $\Phi_{\G\F}$ of the \qc structure, that is
\begin{equation}\label{eq:F DCP}
F_\G=\Phi_{\G\F}^{\otimes 2}\cdot F_\F\cdot\Delta(\Phi_{\G\F})^{-1}
\end{equation}

\Omit{
Option: Discuss/mention in passing the coproduct identity: $\Delta(S_i)
=R_i^{21}\cdot S_i\otimes S_i$ ? This identity is not the main
point. Rigidity does not come from the coproduct identity, but
rather from the fact that the equation Phi_F=Phi' has a unique
solution in F if step is corank 1. The coproduct identity is only
used to show that the gauge transforming one F to another fixes
the local monodromies.}

\subsection{} 

The quantum group $\Uhg$ associated to a complex semisimple Lie
algebra $\g$ is a \qcqtqba in a very simple way. For any subdiagram
$B$ of the Dynkin diagram $D$ of $\g$, the subalgebra $\Uhg_B$ is
the quantum group corresponding to the generators of $\Uhg$ labelled
by the vertices of $B$, endowed with its universal $R$--matrix $R_B$
and trivial associator $\Phi_B=1^{\otimes 3}$. The associators $\Phi_
{\G\F}$ are all trivial, and the local monodromies $S_i$ are given by
Lusztig's \qW group operators. It was proved in \cite[Thm 8.3]{TL2}
that this structure can be transferred to an isomorphic one on $\Ug\fml$
(which, however, does not have trivial associators). Moreover, it was
also proved in \cite[Thm. 9.1]{TL2} that \qcqtqba structures on $\Ug\fml$
are rigid, that is determined by their $R$--matrices $R_B$ and local
monodromies $S_i$.\footnote{This is not the case of \qca structures
on $\Ug\fml$.} Thus, Theorem \ref{th:intro monodromy thm} can be
proved by showing the existence of a \qcqtqba structure on $\Ug\fml$
which binds the \qc structure underlying the monodromy of the Casimir
connection of $\g$, to the \qtqba structure underlying that of the \KZ
connections of its standard Levi subalgebras.

\subsection{} 

It is a well--known, and beautiful observation of Drinfeld's that the
monodromy of the KZ equations of $\g$ gives rise to a \qtqba structure
on $\Ug\fml$ \cite{Dr4}. The corresponding $R$--matrix is the monodromy
$R\KKZ=e^{\hbar\Omega}$ of the KZ equations on $n=2$ points, and
the associator $\Phi\KKZ$ the ratio $\Psi_{(x_1x_2)x_3}^{-1}\cdot\Psi
_{x_1(x_2x_3)}$ of the solutions of the KZ equations on $n=3$ points
corresponding to the asymptotic zones $z_1-z_3>>z_1-z_2$ and
$z_1-z_3>>z_2-z_3$ respectively. The associativity constraints
relating the copies of the $n$--fold tensor power $V^{\otimes n}$
of a representation $V$ corresponding to two bracketings $b,b'$
can be expressed in terms of the associator $\Phi\KKZ$, as in any
monoidal category, or more directly obtained as the ratio $\Psi_{b'}
^{-1}\cdot\Psi_b$ of the solutions of the KZ equations on $n$ points
corresponding to the asymptotic zones determined by $b$ and $b'$.

\subsection{} 

Similarly, the monodromy of the Casimir connection of $\g$ gives
rise to a \qc structure on $\Ug\fml$. This relies on the \DCP construction
of a compactification of $\hreg$, where the root hyperplanes are
replaced by a normal crossings divisor \cite{DCP1,DCP2}. The
irreducible components of the divisor which intersect the closure
of the Weyl chamber are labelled by the connected subdiagrams
of $D$. The \mnss on the Dynkin diagram $D$ of $\g$ label the
points at infinity, that is the non--empty intersection of a maximal
collection of these components. Near each of those, one can
construct a canonical fundamental solution $\Psi_\F$ having good
asymptotics. The associators $\Phi_{\G\F}$ then arise as the ratios
$\Psi_\G^{-1}\cdot \Psi_\F$.

\subsection{} 

The previous paragraphs suggest that the relative twists of
a \qcqtqba structure on  $\Ug\fml$ which binds the structures
coming from the KZ and Casimir equations might also arise
by comparing appropriate solutions of a flat connection. This
is indeed the case, as we explain below.

Since the associator $\Phi_{\G\F}$ of the \qc structure underlying
the Casimir connection $\nablac$ is equal to $\Psi_\G^{-1}\cdot
\Psi_\F$, where the latter are the \DCP fundamental solutions
of $\nablac$, the compatibility relation \eqref{eq:F DCP} may
be rewritten as
\begin{equation}\label{eq:the observation}
\Psi_\G\ot{2}\cdot F_\G\cdot\Delta(\Psi_\G)^{-1}=
\Psi_\F\ot{2}\cdot F_\F\cdot\Delta(\Psi_\F)^{-1}
\end{equation}
Either side defines a holomorphic function $F:\hreg\to\Ug
\ot{2}\fml$ which is independent of the choice of a \mns
on $D$ by \eqref{eq:the observation}, satisfies the differential
equation
\begin{equation}\label{eq:Casimir for F}
dF=
\half{\sfh}\sum_{\alpha>0}\dloga\left((\onetwo{\calpha})F-
F\Delta(\calpha)\right)
\end{equation}
and the twist relation
$(\Phi\KKZ)_F=1\ot{3}$. 
We shall call such an $F$ a {\it differential twist}.
Given $F$, the twists $F_\F$ may be recovered as
\[F_\F=(\Psi_\F\ot{2})^{-1}\cdot F\cdot\Delta(\Psi_\F)\]

\subsection{}\label{ss:centraliser} 

We show in Section \ref{se:diffl twist} that the requirement that
the twists $F_\F$ possess the factorised form \eqref{eq:factorised F},
where the factors $F_{(B;\alpha)}$ satisfy the twist relation \eqref
{eq:Phi F_B}, is equivalent to the following {\it centraliser property}
of the differential twist $F$. This assumes that a differential twist
$F_{\g_B}$ is given for any diagrammatic subalgebra $\g_B\subseteq
\g$, and expresses a compatibility  between $F_{\g_B}$ and $F_
{\g_{B\setminus\alpha}}$, for any vertex $\alpha$ of $B$. Specifically,
consider the asymptotics $\rlim{\alpha\to\infty}F_{\g_B}(\mu)$
of $F_{\g_B}(\mu)$, as the coordinate $\alpha$ of $\mu
\in\h_B$ goes to $\infty$.\footnote{the existence $\rlim{\alpha\to
\infty}F_{\g_B}$ relies on the fact that the Casimir connection has
regular singularities.} These asymptotics are a function of the image
$\ol{\mu}$ of $\mu$ in $\h_{B\setminus\alpha}$, thought
of as a quotient of $\h_B$. They solve the Casimir equation
\eqref{eq:Casimir for F} for the Lie algebra $\g_{B\setminus\alpha}$,
which are the limit of the Casimir equations of $\g_B$ as $\alpha\to
\infty$. The centraliser property is the requirement that the element
$F_{(B;\alpha)}\in\Ug\fml\ot{2}$ defined by the equation
\begin{equation}\label{eq:intro centraliser}
\rlim{\alpha\to\infty}F_{\g_B}(\mu)=
F_{(B;\alpha)}\cdot F_{\g_{B\setminus\alpha}}(\ol{\mu})
\end{equation}
be invariant under $\g_B$. This implies in particular that $F_{(B;
\alpha)}$ does not depend on $\ol{\mu}$. 

\subsection{} 

The construction of a differential twist satisfying the centraliser
property can, in turn, be obtained from that of an appropriate
{\it fusion operator}.\footnote{The terminology is borrowed
from the work of Etingof (see \eg \cite{E}, and references
therein). The relation of our construction to \cite{EV} is
discussed in \ref{ss:EV}.} The latter is a solution of a
dynamical version of the KZ equations in $n=2$ points,
namely
\begin{equation}\label{eq:intro DKZ}
\frac{dJ}{dz}=\left(\sfh\frac{\Omega}{z}+\ad\muone\right)J
\end{equation}
where $z=z_1-z_2$, and $\mu\in\hreg$. The dynamical
KZ equations arise naturally in Conformal Field Theory
(see, \eg \cite{EFK}). They were studied in more detail
by \FMTV in \cite{FMTV}, where it was shown that they
are bispectral to, that is commute with, a dynamical version
of the Casimir connection \wrt the variable $\mu$.

The presence of the dynamical term $\ad\muone$ creates an
irregular singularity at $z=\infty$ of Poincar\'e rank 1. Assuming
$\mu$ to be real, so that the Stokes rays of \eqref{eq:intro DKZ}
all lie on the real axis, we construct in Section \ref{se:fusion} two
canonical solutions $J_\pm(z,\mu)$ with values in $\Ug\ot{2}\fml$,
which have the form
\[J_\pm(z,\mu)=H_\pm(z,\mu)\cdot z^{\sfh\Omega_\h}\]
where $H_\pm(z,\mu)$ is a holomorphic function on the upper
(resp. lower) half--plane which possesses an asymptotic expansion
of the form $1^{\otimes 2}+O(z^{-1})$ as $z\to\infty$ with $0<<|\arg z|<<\pi$,
and $\Omega=t_a\otimes t^a$, with $\{t_a\},\{t^a\}$ dual bases
of $\h$ \wrt $(\cdot,\cdot)$, is the projection of $\Omega$ on
the kernel of $\ad(\muone)$. The construction of $J_\pm$
and the study of its analytic properties require some care,
since the equation \eqref{eq:intro DKZ} takes values in the
infinite--dimensional algebra $\Ug\ot{2}\fml$, and is the main
technical contribution of the paper.

Given the fusion operator $J_\pm(z,\mu)$, a differential twist
can be obtained as either of the ratios
\[F_\pm(\mu)=J_\pm(z,\mu)^{-1}\cdot J_0(z,\mu) \]
where $J_0(z,\mu)$ is the unique solution of \eqref{eq:intro DKZ}
which is asymptotic to $(1^{\otimes 2}+O(z))\cdot z^{\hbar\Omega}$ near
$z=0$.\footnote{The existence of $J_0$ is straightforward since
$z=0$ is a regular singularity of \eqref {eq:intro DKZ}.}

\subsection{} 

The proof that $F_\pm(\mu)$ kills the associator is fairly standard
(see \eg \cite{E,EE}), provided an analogue of the fusion operator
can be constructed for the dynamical KZ equation in $n=3$ points.
Even though the corresponding KZ equations have irregular singularities
at $z_i=\infty$, a solution can still be constructed by simultaneously scaling
all variables $z_1,z_2,z_3$, and sending the scale to infinity. Here,
we crucially exploit a beautiful fact, 
which is that the dynamical KZ equations in
$n$ variables {\it abelianise at infinity}. Roughly speaking, this
means that the connection satisfied by the asymptotics of
solutions on the divisor where all $z_i-z_j$ are infinite is the
abelian KZ equations
\[d-\sfh\sum_{i<j}d\log(z_i-z_j)\Omega_{ij}^\h\]
Contrary to its non-abelian analogue, this equation possesses a
canonical solution, namely $\prod_{i<j}(z_i-z_j)^{\sfh\Omega_{ij}^\h}$,
which leads to the construction of a multicomponent fusion operator
of the form
\[J_\pm=H_\pm(z_1,\ldots,z_n)\cdot\prod_{i<j}(z_i-z_j)^{\sfh\Omega_{ij}^\h}\]

\subsection{} 

The fact that the twists $F_\pm(\mu)$ satisfy the Casimir equations
\ref{eq:Casimir for F} follows from the fact that $J_\pm(z,\mu)$
satisfies a dynamical version of the Casimir equations, namely
\[d_\h J=
\half{\sfh}
\sum_{\alpha>0}\dloga\left(\Delta(\calpha)J-J(\onetwo{\calpha})\right)+
z\ad(d\muone)J\]
which expresses the fact that, when $z_1-z_2=\infty$, the dynamical
Casimir connection on the tensor product $V_1\otimes V_2$ of two
representations becomes the tensor product of the (non--dynamical)
Casimir connections on each factor. The above equation is 
a consequence of the bispectrality of the dynamical KZ and Casimir
equations, 
together with the fact that $J$ varies smoothly in $\mu\in\hreg$, which
follows from our analysis.

\subsection{} 

In Section \ref{se:centraliser}, we show that the differential twists 
$F_\pm(\mu)$ satisfy the centraliser property. This follows by
relating the (irregular) asymptotics of the fusion operator of $\g$
when one of the simple roots coordinates $\alpha_i$ tends to $\infty$,
to the fusion operator of the corank 1 subalgebra $\g_{D\setminus
\alpha_i}$. \Omit{The latter, in turn, rests on the fact that the limit, as
$\alpha_i\to\infty$ of the dynamical KZ equation \eqref{eq:intro DKZ}
are essentially the corresponding equations for $\g_{D\setminus
\alpha_i}$, namely
\[\frac{dJ}{dz}=
\left(\sfh\frac{\Omega_{\g_{D\setminus\alpha_i}}+\Lambda}{z}+\ad\ol{\mu}^{(1)}\right)J\]
where $\Lambda=\cow{i}\otimes\cow{i}/\|\cow{i}\|^2$, and $\cow{i}$
is the fundamental coweight corresponding to $\alpha_i$.}

We revisit these calculations in Appendix \ref{app:revisited}, and
give a direct construction of the relative twists arising from the
fusion operator which is similar in spirit to Drinfeld's construction
of the KZ associator. This gives, to the best of our knowledge,
the first canonical transcendental construction of a twist killing
the KZ associator $\Phi\KKZ$.

Specifically, we show that the twist $F_{(D;\alpha_i)}$ defined
by the factorisation relation \eqref{eq:intro centraliser} can be
realised as the constant relating the canonical fundamental
solutions at 0 and $\infty$ of an ODE with regular singularities
at $0,-1$ and an irregular singularity at $\infty$. The ODE is
defined \wrt the blowup coordinate $v=z\alpha_i$, 
and is given by
\[\frac{dG}{dv}=
\left(\frac{\sfh\Omega}{v}
+\frac{\sfh}{2}\frac{\Delta(\K_D-\K_{D\setminus\alpha_i})-2\Omega}{v+1}
+\ad\cowone{i}
\right)G\]
where, for any subdiagram $B\subseteq D$, $\K_B$ is the
truncated (Cartan--less) Casimir operator of $\g_B$,
and $\cow{i}$ is the coroot corresponding to $\alpha_i$.

\subsection{}\label{ss:EV} 

Our construction of the differential twist is very close, in spirit at least,
to Etingof and Varchenko's study of the fusion operator $\wh{J}(w,
\lambda)$ of the affine Lie algebra $\wh{\g}$ \cite{EV}. The latter
satisfies the trigonometric KZ equations \wrt $w$, together with the
dynamical difference equations \wrt $\lambda$, where the latter are
a system of difference equations which degenerates to the Casimir
connection. The regularised limit $\wh{J}\lambda)$ of $\wh{J}(z,\lambda)$
as $z\to 1$ kills the KZ associator in the shifted sense, that is satisfies
\[\wh{J}_{12,3}(\lambda)\wh{J}_{1,2}(\lambda-h^{(3)})=
\wh{J}_{1,23}(\lambda)\wh{J}_{2,3}(\lambda)\Phi\KKZ\]
(see, \eg \cite{EE,E}).  A construction of a differential twist might
therefore arise by taking an appropriate scaling limit of $\wh{J}
(\lambda)$ as $\lambda$ goes to infinity (a process which would
kill the shift in the above equation). Controlling the asymptotics of
$\wh{J}$ at $\lambda=\infty$ seems difficult, however, since $\wh
{J}(\lambda)$ only satisfies a system of difference equations \wrt
$\lambda$.

Rather than purse this path, we give in this paper a direct construction
of a solution of the rational dynamical KZ equations, which can be
thought of as (and probably is) a degeneration of the fusion operator
of $\hat{\g}$. Our $J_\pm(z,\mu)$ should in fact be a fusion operator
for the current algebra $\g[t]$. One further difference with \cite{EV} is
that, unlike $\wh{J}(w,\lambda)$, $J_\pm(z,\mu)$ is not constructed via
representation theory, specifically the fusion construction for loop
modules of $\wh{\g}$, but purely using differential equations. It seems
an interesting question to construct our fusion operator from representation
theory. Such a construction should be obtained by replacing Verma
modules by the irregular Wakimoto modules for $\wh{\g}$ considered
in \cite{FFT,Fe}, since these give rise to the Casimir connection \cite{Fe}.

\subsection{} 

The monodromy theorem proved in this paper is extended to the
case of an arbitrary symmetrisable \KM algebra in \cite{ATL1,ATL2,
ATL3}. The approach is close in spirit to that \cite{TL2}, but differs
very significantly in the details of two out of the three steps of the
construction, namely the transfer of braided \qc structure from
the category $\Oint_\hbar$ of integrable, highest weight $\Uhg
$--modules to an isomorphic structure on the corresponding
category $\Oint\fml$ for $\Ug\fml$, and the proof that such 
structures are rigid. The last step, namely the construction of
a braided \qc structure on $\Oint\fml$ which accounts for the
monodromy of the Casimir equations of $\g$ and that of the
KZ equations of its Levi subalgebras is carried out in \cite{ATL3}
by using the construction of the fusion operator and differential
twist given in this paper.

\Omit{
. Specifically, since $\Ug\fml$ and $\Uhg$ are not isomorphic
as algebras if $\g$ is infinite--dimensional, the transfer of \qcqtqba
structure on $\Uhg$ to an isomorphic one on $\Ug\fml$ obtained
in \cite{TL2} is replaced by one of the induced structures on the
categories $\Oint_\hbar$, $\Oint\fml$ of integrable, highest weight
modules over $\Uhg$ and $\Ug\fml$ respectively. This involves developing the notion of {\it
braided, \qc tensor category}, and modifying the \EK quantisation
functor so that it gives an equivalence of $\Oint_\hbar$ and $\Oint
\fml$ compatible with the inclusion of standard Levi subalgebras
\cite{ATL1}. The rigidity of such structures is proved in \cite{ATL2},
and involves a complex which replaces the Hochschild complex of
the completion of $\Ug$ \wrt category $\O$, since the latter does not
have a well--behaved cohomology. The new complex arises from the
PROP consisting of a bialgebra with a root space decomposition, and
leads to a far stronger rigidity result. The final step in the proof namely
the construction of a braided \qc tensor structure on $\Oint\fml$ which
accounts for the monodromy of the Casimir equations of $\g$ and that
of the KZ equations of its Levi subalgebras is obtained in \cite{ATL3}
relies on a 
}

\subsection{Outline of the paper}

Section \ref{se:top prel} contains some preliminary material required
to study differential equations with values in infinite--dimensional
filtered vector spaces.

In Section \ref{se:qca}, we review the definition of \qtqbas and of
\qcas together with the fact that the monodromy of the KZ and Casimir
connections respectively define such structures on the enveloping
algebra $U\g$ of a complex, semisimple Lie algebra $\g$.

In Section \ref{se:diffl twist}, we introduce the notion of differential twist
for $\g$, and show that it gives rise to a \qcqtqba structure on $U\g$,
which interpolates between the \qtqba structure underlying the monodromy
of the KZ connection and the \qc structure underlying that of the Casimir
one.

In Section \ref{se:fusion}, we construct a fusion operator for $\g$
as a joint solution of the coupled KZ--Casimir equations on $2$
points, with prescribed asymptotics when $z_1-z_2\to\infty$. 

In Section \ref{se:diffl from fusion}, we obtain a differential twist for
$\g$ as the regularised limit of the fusion operator when $z\to 0$,
and prove that it kills the KZ associator.

In Section \ref{se:centraliser}, we relate the differential twists for
$\g$ and for a corank 1 Levi subalgebra, and use this to prove
that the differential twist of $\g$ satisfies the centraliser property
described in \ref{ss:centraliser}.

This property is used in Section \ref{se:diffl qcqtqba}, to show that
the differential twist arising from the fusion operator it gives rise to
a \qcqtqba structure on $U\g$ interpolating between the \qtqba
structure underlying the KZ equations and the \qca one underlying
the Casimir connection.

Section \ref{se:Uhg} collects some facts about the quantum group
$\Uhg$, and in particular the fact that it possesses a \qcqtqba structure
which accounts for both its $R$--matrix and quantum Weyl group
representations.

Section \ref{se:monodromy} contains the main result of this paper,
namely the equivalence of $\Uhg$ and $U\g\fml$ as \qcqtqbas and
the immediate corollary that the monodromy of the Casimir connection
is described by \qW group operators.

Appendix \ref{se:basic ODE} contains a detailed discussion of
the solutions of a linear, scalar ODE with an irregular singularity
at $\infty$, which plays a similar role in this paper than Drinfeld's
ODE $\frac{df}{dz}=(\frac{A}{z}+\frac{B}{z-1})f$ does for
the construction of the KZ associators.

Appendix \ref{app:revisited} gives an alternative proof that the
differential twist obtained in Section \ref{se:diffl from fusion}
possesses the centraliser property. As a corollary, we obtain
a canonical, transcendental construction of a twist killing the
KZ associator $\Phi\KKZ$.

The final Appendix \ref{app:joint} gives an alternative proof
that the coupled KZ--Casimir equations are integrable.

\subsection{Acknowledgments}

This paper was a very long time in the making, as its main result was
first announced in 2005. I am very grateful to my late colleague Andrei
Zelevinsky, and to Edward Frenkel for amicably, but firmly, encouraging
me to complete it.
I am also grateful to Claude Sabbah for correspondence and discussions
on irregular singularities in 2 dimensions.

\section{Filtered algebras} 
\label{se:top prel}

\subsection{} 

Let $A$ be an algebra over $\IC$ endowed with an ascending filtration
\[\IC=A_0\subset A_1\subset\cdots\]
such that $A_m\cdot A_n\subset A_{m+n}$. Let $o=\{o_k\}_{k\in\IN}$
be a sequence of non--negative integers, $\hbar$ a formal variable,
and consider the subspace $A\fml^o\subset A\fml$ defined by
\[A\fml^o=\{\sum_{k\geq 0}a_k\hbar^k|\,a_k\in A_{o_k}\}\]
Note that:
\begin{enumerate}
\item $A\fml^o$ is a (closed) $\IC\fml$--submodule of $A\fml$ if $o$ is
increasing.\footnote{$A\fml^o$ is then the $\hbar$--adic completion of the Rees algebra
of $A$ corresponding to the filtration $A_{o_0}\subset A_{o_1}\subset\cdots$}
\item $A\fml^o$ is a subalgebra of $A\fml$ if $o$ is subadditive, that is such
that $o_k+o_l\leq o_{k+l}$ for any $k,l\in\IN$. This implies in particular that
$o$ is increasing, and that $o_0=0$.
\end{enumerate}

\Omit{
\begin{rem}
The smallest subadditive sequence $\{o_k\}$ with a given $o_1$ is
$o_k=k\cdot o_1$.
\end{rem}
}

If the subspaces $A_k\subset A$ are finite--dimensional, so are the quotients
\[A\fml^o/(\hbar^{p+1}A\fml\cap A\fml^o)\cong
A_{o_0}\oplus\hbar A_{o_1}\oplus \cdots\oplus\hbar^p A_{o_p}\]
Assuming this, we shall say that a map $F:X\to A\fml^o$, where $X$ is a
topological space (resp. a smooth or complex manifold), is continuous
(resp. smooth or holomorphic) if each of its truncations $F_p:X\to A\fml^
o/(\hbar^{p+1}A\fml\cap A\fml^o)$ are.

\subsection{}

We shall mainly be interested in the following situation: $A=\Ug^{\otimes n}$
endowed with the standard order filtration given by $\deg(x\ii)=1$ for $x\in\g$,
where
\[x\ii=1^{\otimes(i-1)}\otimes x\otimes 1^{\otimes (n-i)}\]
The sequence $o$ will be chosen subadditive, and such that $o_1\geq 2$ in
order for $\hbar\Omega_{ij},\hbar\Delta^{(n)}(\Kalpha)\in A\fml^o$. Note that $\g\cap U\g\fml^o=\{0\}$ since $o_0=0$, but that the
adjoint action of $\g$ on $\Ug^{\otimes n}$ induces one on by derivations on
$\Ug^{\otimes n}\fml^o$. Note also that $\Ug\fml^o$ is {\bf not} a Hopf algebra,
since $\Delta:\Ug\fml^o\to\Ug^{\otimes 2}\fml^o\supsetneq\left(\Ug\fml^o\right)
^{\otimes 2}$.

\Omit{
\begin{lemma}
$A\fml^o$ is a topologically free $\IC\fml$--module.
\end{lemma}
\begin{pf}
$A\fml^o$ is clearly torsion--free and separated since $\bigcap_n\hbar^n A\fml^o
\subset \bigcap_n\hbar^n A\fml=0$. To show completeness, note that
$\hbar^n A\fml^o=\{\sum_{k\geq 0}a_k\hbar^k|\,a_k\in A_{k-n}\}$ where 
$A_{-m}=0$ for $m>0$. Thus, $A\fml^o/\hbar^n A\fml^o$ is isomorphic
to
\[\{\sum_{k\geq 0}a_k\hbar^k|\,a_k\in A_k, k<n\,\text{and}\,a_k\in A_k/A_{k-n}, k\geq n\}\]
from which the completeness of $A\fml^o$ follows.
\end{pf}
}

\subsection{}

Let $\A$ be a $\IC\fml$--module and consider the natural map
$\imath:\A\to\ds{\lim_{\longleftarrow}}\,\A/\hbar^n\A$. Recall that
$\A$ is {\it separated} if $\imath$ is injective, and {\it complete}
if it is surjective. By definition, $\A$ is {\it topologically free} if it is
separated, complete and torsion--free.

Consider now the map
\[\imath:\End_{\IC\fml}(\A)\to
\lim_{\longleftarrow}\End_{\IC\fml}(\A/\hbar^n\A)\]

\begin{lemma} Assume that $\A$ is separated. Then,
\begin{enumerate}
\item 
$\imath$ is injective.
\item If $\A$ is complete, $\imath$ is surjective.
\end{enumerate}
\end{lemma}
\begin{pf}
(1) If $T\in\End_{\IC\fml}$ is such that $\imath T=0$, then $T(\A)\subset
\bigcap_n\hbar^n\A=0$.
(2) Let $\{T_n\}\in\lim_n\End_{\IC\fml}(\A/\hbar^n\A)$. For any $a\in\A$,
the sequence $\{T_n a\}$ lies in $\lim_n\A/\hbar^n\A$ and is therefore
the image of a unique element $a'\in\A$. The assignment $a\to Ta=a'$
is easily seen to define an element of $\End_{\IC\fml}(\A)$ which projects
to each of the $T_n$.
\end{pf}

\begin{corollary}
If $\A$ is topologically free, the map $\imath:\End_{\IC\fml}(\A)\to
\ds{\lim_{\longleftarrow}}\End_{\IC\fml}(\A/\hbar^n\A)$ is an isomorphism.
\end{corollary}

\Omit{
\begin{rem}
If $\A\subseteq\B$ are $\IC\fml$--modules, then $\A$ is separated
(resp. torsion--free) is $\B$ is. However, the completeness of $\B$
does not imply that of $\A$. Indeed, if $V$ is an infinite--dimensional,
complex vector space, then $\A=V\otimes_{\IC}\IC\fml\subsetneq
V\fml=\B$, and $\B$ is complete, while $\A$ is not.
\end{rem}
}

\section{Quasi--Coxeter algebras} 
\label{se:qca}

We review in this section the definition of \qtqbas following
\cite{Dr3}, and of quasi--Coxeter algebras following \cite{TL2},
to which we refer for more details. We also explain how the
monodromy of the \KZ (KZ) and Casimir connections of a
complex, semisimple Lie algebra $\g$ respectively give rise
to a \qtqba and \qca structure on the enveloping algebra
$U\g\fml$ of $\g$.

\subsection{Quasitriangular quasibialgebras \cite{Dr3}}

\subsubsection{}\label{sss:qba recap}

Recall that a \qba $(A,\Delta,\veps,\Phi)$ is an algebra $A$ endowed with algebra
homomorphisms $\Delta:A\rightarrow A^{\otimes 2}$ and $\veps:A\rightarrow
k$ called the coproduct and counit, and an invertible element $\Phi\in
A^{\otimes 3}$ called the associator which satisfy, for any $a\in A$
\begin{align*}
\id\otimes\Delta(\Delta(a))
&=
\Phi\cdot\Delta\otimes\id(\Delta(a))\cdot\Phi^{-1}
\\
\id^{\otimes 2}\otimes\Delta(\Phi)\cdot
\Delta\otimes\id^{\otimes 2}(\Phi)
&=
1\otimes\Phi\cdot
\id\otimes\Delta\otimes\id(\Phi)\cdot
\Phi\otimes 1\\
\veps\otimes\id\circ\Delta
&=\id\\
\id\otimes\veps\circ\Delta
&=\id\\
\id\otimes\veps\otimes\id(\Phi)
&=1
\end{align*}

A twist of a \qba $A$ is an invertible element $F\in A^
{\otimes 2}$ satisfying $$\veps\otimes\id(F)=1=\id\otimes\veps(F)$$
Given such an $F$, the twisting of $A$ by $F$ is the \qba $(A,\Delta_F,
\veps,\Phi_F)$, where the coproduct $\Delta_F$ and associator $\Phi
_F$ are given by
\begin{align*}
\Delta_F(a)&=F\cdot\Delta(a)\cdot F^{-1}\\
\Phi_F
&=
1\otimes F\cdot\id\otimes\Delta(F)\cdot
\Phi\cdot
\Delta\otimes\id(F^{-1})\cdot F^{-1}\otimes 1
\end{align*}

A strict morphism $\Psi:A\rightarrow A'$ of \qbas is an
algebra homomorphism satisfying
$$\veps=\veps'\circ\Psi,\qquad
\Psi^{\otimes 2}\circ\Delta=\Delta'\circ\Psi
\qquad\text{and}\qquad
\Psi^{\otimes 3}(\Phi)=\Phi'$$
A morphism $A\rightarrow A'$ of \qbas is a pair $(\Psi,F')$ where
$F'$ is a twist of $A'$ and $\Psi$ is a strict morphism of $A$ to the
twisting of $A'$ by $F'$.

\subsubsection{}\label{sss:qtqba recap}

A \qba $(A,\Delta,\veps,\Phi)$ is \qt if it is endowed with an invertible
element $R\in A^{\otimes 2}$ called the $R$--matrix satisfying, for
any $a\in A$,
\begin{align}
\Delta\op(a)&=R\cdot\Delta(a)\cdot R^{-1}
\label{eq:Delta R}\\
\Delta\otimes\id(R)&=
\Phi_{312}\cdot R_{13}\cdot\Phi_{132}^{-1}\cdot R_{23}\cdot\Phi_{123}
\label{eq:R12}\\
\id\otimes\Delta(R)&=
\Phi_{231}^{-1}\cdot R_{13}\cdot\Phi_{213}\cdot R_{12}\cdot\Phi_{123}^{-1}
\label{eq:R23}
\end{align}

A twist $F$ of a \qtqba $A$ is a twist of the underyling quasibialgebra.
The twisting of $A$ by $F$ is the \qtqba $(A,\Delta_F,\veps,\Phi_F,R_F)$
where 
$$R_F=F_{21}\cdot R\cdot F^{-1}$$
A morphism $(\Psi,F'):A\rightarrow A'$ of \qt \qbas is a morphism of the
underlying \qbas such that $\Psi^{\otimes 2}(R)=R'_{F'}$.



\subsubsection{}\label{sss:KZ}

Let $\g$ be a complex, semisimple Lie algebra, with symmetric, invariant
tensor $\Omega\in(S^2\g)^\g$. Then, $U\g\fml$ is a \qtqba  with the standard
(cocommutative) coproduct $\Delta_0$ and counit $\veps$, $R$--matrix
given by $R^{\KKZ}=\exp(\hbar\Omega)$ and associator $\Phi^{\KKZ}$
constructed as follows.

Consider the differential equation with values in $U\g^{\otimes 3}\fml^o$
given by
\[\frac{dG}{dz}=\sfh\left(\frac{\Omega_{12}}{z}+\frac{\Omega_{23}}{z-1}\right)G\]
where $\Omega_{12}=\Omega\otimes 1$, $\Omega_{23}=1\otimes\Omega$,
and $\sfh=\hbar/\pi\iota$. This equation has a unique solution $G_0$ of the
form $H_0(z)\cdot z^{\sfh\Omega_{12}}$, where $H_0$ is holomorphic on
the disk $\{z\in\IC|\,|z|<1\}$ and such that $H_0(0)=1$, and $z^{\hbar\Omega
_{12}}=\exp(\sfh\Omega_{12}\log z)$, where $\log$ is the standard determination
of the logarithm. Similarly, there is a unique solution $G_1$ of the form
$H_1(z)(1-z)^{\sfh\Omega_{23}}$, where $H_1$ is holomorphic on $\{z
\in\IC|\,|z-1|<1\}$, and such that $H_1(1)=1$. The KZ associator $\Phi^
{\KKZ}$ is defined by
\[\Phi^{\KKZ}=G_1^{-1}(x)\cdot G_0(x)\]
where $x\in(0,1)$.

\subsection{Quasi--Coxeter algebras \cite{TL2}}

\subsubsection{Diagrams} 

A {\it diagram} is an 
undirected graph $D$ with no multiple
edges or loops. We denote the set of vertices of $D$ by $V(D)$, and
set $|D|=|V(D)|$. A {\it subdiagram} $B\subset D$ is a full subgraph
of $D$, that is a graph consisting of a subset $V(B)$ of vertices of
$D$, together with all edges of $D$ joining two elements of $V(B)$.
We shall often identify a subdiagram $B$ and its set of vertices $V(B)$.

Two subdiagrams $B_1,B_2\subseteq D$ are {\it orthogonal} if no two
vertices $\alpha_1\in B_1,\alpha_2\in B_2$ are joined by an edge in
$D$. $B_1$ and $B_2$ are {\it compatible} if either one contains the
other or they are orthogonal.

\subsubsection{Nested sets} 

Assume henceforth that $D$ is connected. A {\it nested set} on $D$
is a collection $\H$ of pairwise compatible, connected subdiagrams
of $D$ containing $D$.

We denote by $\ND$ the partially ordered set of nested sets on $D$,
ordered by reverse inclusion. $\ND$ has a unique maximal element
$\{D\}$. Its minimal elements are the maximal nested sets.

It is easy to see that a \mns $\F$ has the following properties
\begin{enumerate}
\item The cardinality of $\F$ is $|D|$.
\item If $B\in\F$, the maximal elements $\{B_i\}$ in $\F$ properly
contained in $B$ contain all the vertices of $B$ with the exception
of one, which will be denoted $\alpha^B_\F$. 
\end{enumerate}

The $B_i$ are in fact the connected components of the diagram
$B\setminus\alpha^B_\F$, and $\F$ may be obtained by taking
the connected  components $D_i$ of $D\setminus\alpha^D_\F$,
then those of each $D_i\setminus\alpha^{D_i}_\F$ and so on.

\subsubsection{Type $\sfA$}

If $D$ is the Dynkin diagram of type $\sfA_{n-1}$, with vertices labelled $1,
\ldots,n-1$, nested sets on $D$ are in bijection with bracketings of the
non associative monomial $x_1\cdots x_n$. The bijection is obtained by
mapping a connected subdiagram with vertices $i,\ldots,j-1$ to the pair
of parentheses $x_1\cdots x_{i-1}\left(x_i\cdots x_j\right)x_{j+1}\cdots x_n$
and noting that $B,B'\subseteq D$ are compatible if, and only if the
corresponding pairs of parentheses are. It follows that, in this case,
the poset $\ND$ is the face poset of Stasheff's associahedron.

\subsubsection{$D$--algebras}

Fix henceforth a commutative, unital ring $\sfk$.
A {\it $D$--algebra} is a $\sfk$--algebra $A$ endowed with subalgebras $A_B$
labelled by the non--empty connected subdiagrams $B\subseteq D$
such that the following holds
\begin{itemize}
\item $A_{B'}\subseteq A_{B\dprime}$ whenever $B'\subseteq B\dprime$.
\item $A_{B'}$ and $A_{B\dprime}$ commute whenever $B'$ and $B\dprime$
are orthogonal.
\end{itemize}

If $B_1,B_2\subseteq D$ are subdiagrams with $B_1$ connected,
we denote by $A_{B_1}^{B_2}$ the centraliser in $A_{B_1}$ of the subalgebras
$A_{B_2'}$ where $B_2'$ runs over the connected components of $B_2$.

A morphism of $D$--algebras $A,A'$ is a collection of algebra
homomorphisms $\Psi_\F:A\rightarrow A'$ labelled by the \mnss
on $D$ such that for any $\F$ and $B\in\F$, $\Psi_\F(A_{B})
\subseteq A'_{B}$.

\subsubsection{Pairs of maximal nested sets.}

An {\it elementary pair} $(\F,\G)$ of \mnss on $D$ is one for which
$\F$ and $\G$ differ by an element. In this case, there is a unique
minimal element $B$ of $\F\cap\G$, called the {\it support} of $(\F,
\G)$, which contains the single element in $\F\setminus\G$ (or
$\G\setminus\F$). The {\it central support} $\zsupp(\F,\G)$ of
$(\F,\G)$ is the union of the elements of $\F\cap\G$ properly
contained in $\supp(\F,\G)$.

In type $\sfA_{n-1}$, an elementary pair $(\F,\G)$ corresponds
to two complete parenthetisations of $x_1\cdots x_n$ which are
obtained by replacing a pair of parentheses $\cdots (x_i\cdots
x_j)\cdots$ with $\cdots (x_k\cdots x_\ell)\cdots$, where $i<k<j<\ell$,
and the support of $(\F,\G)$ to the smallest pair of parentheses
$\cdots (x_i\cdots x_\ell)\cdots$ which is consistent with both
pairs. 

\subsubsection{Weak \qc algebras}\label{sss:weak qc}

A {\it weak \qc algebra} $A$ of type $D$ is a $D$--algebra $A$ endowed
with invertible elements $\Phi_{\G\F}$ called {\it \DCP associators}
labelled by pairs of \mnss satisfying the following axioms\footnote
{The term weak \qc algebra is borrowed from \cite{ATL2}.}
\begin{itemize}
\Omit{
\item {\bf Orientation.} For any pair $(\G,\F)$ of \mnss on $D$
\[\Phi_{\G\F}=\Phi_{\F\G}^{-1}\]
}
\item {\bf Transitivity.} For any \mnss $\F,\G,\H$ on $D$
\[\Phi_{\H\F}=\Phi_{\H\G}\cdot\Phi_{\G\F}\]
\item {\bf Support.} For any elementary pair $(\G,\F)$ of
\mnss on $D$,
\[\Phi_{\G\F}\in A_{\supp(\G,\F)}^{\zsupp(\G,\F)}\]
\item {\bf Forgetfulness.} For any elementary pairs $(\G,\F)$, $(\G',\F')$
of \mnss on $D$ such that $\F\setminus\G=\F'\setminus\G'$ and
$\G\setminus\F=\G'\setminus\F'$,
\[\Phi_{\G\F}=\Phi_{\G'\F'}\]
\end{itemize}

A morphism of quasi--Coxeter algebras $A,A'$ of type $D$ is
a morphism $\{\Psi_\F\}$ of the underlying $D$--algebras such
that for any elementary pair $(\G,\F)$ of \mnss on $D$,
\[\Psi_\G\circ\Ad(\Phi^A_{\G\F})=\Ad(\Phi^{A'}_{\G\F})\circ\Psi_\F\]

\subsubsection{Completion}

Let $\Veck$ be the category of finitely--generated, free $\sfk$--modules
and $ \Mod{A}$ that of $A$--modules whose underlying $\sfk$--module
lies in $\Veck$. By definition, the completion of $A$ with respect to its \fd
representations is the algebra $\wh{A}$ of endomorphisms of the
forgetful functor
\[\bfF:\Mod{A}\rightarrow\Veck\]
An element of $\wh{A}$ is a collection $\Theta=\{\Theta_V\}$, with $\Theta
_V\in\End_{k}(V)$ for any $V\in\Mod{A}$, such that for any $U,V\in\Mod{A}$
and $f\in\Hom_A(U,V)$
$$\Theta_V\circ f=f\circ\Theta_U$$

\subsubsection{Labelled diagrams}

A {\it labelling} of $D$ is the assignement of an integer $m_{ij}\in\{2,3,
\ldots,\infty\}$ to any pair $\alpha_i,\alpha_j$ of distinct vertices of $D$
such that $m_{ij}=m_{ji}$ for any $\alpha_i\neq\alpha_j$ and $m_{ij}=2
$ if, and only if $\alpha_i$ and $\alpha_j$ are orthogonal.

If $D$ is labelled, the Artin or braid group $B_D$ is the group generated by
elements $S_i$ corresponding to the vertices $\alpha_i$ of $D$ with
relations
\[\underbrace{S_iS_j\cdots}_{m_{ij}}=
\underbrace{S_jS_i\cdots}_{m_{ij}}\]
for any $\alpha_i\neq\alpha_j$ such that $m_{ij}<\infty$.


\subsubsection{Quasi--Coxeter algebras}\label{sss:qc}

\comment{Not clear that breaking up into weak and full qc structure
is useful. Monodromy examples do not come in this form.}
Let $D$ be a labelled diagram. A \qc algebra of type $D$ is a weak
\qc algebra $A$ endowed with an invertible element
\[S^A_i\in\wh{A}_i\]
for any vertex $\alpha_i\in D$, where $\wh{A}_i$ is the completion
of $A_{\alpha_i}$ \wrt its \fd representations. The associators
$\Phi_{\G\F}$ and local monodromies are required to satisfy the
following 
\begin{itemize}
\item {\bf Braid relations.} For any pair $\alpha_i,\alpha_j$ of distinct
vertices of $D$ such that
$m_{ij}<\infty$, and pair 
$(\G,\F)$ of \mnss on $D$ such that $\alpha_i\in\F$ and $\alpha_j\in\G$,
\[\Ad(\Phi_{\G\F})(S^A_i)\cdot S^A_j\cdots=
S^A_j\cdot\Ad(\Phi_{\G\F})(S^A_i)\cdots\]
where the number of factors on each side is equal to $m_{ij}$.
\end{itemize}

A {\it morphism} $\{\Psi_\F\}$ of \qc algebras $A,A'$ is one of the
underlying weak \qc algebras such that $\Psi_\F(S_i^A)=S_i^{A'}$
for any \mns $\F$ on $D$, and ertex $\alpha_i\in D$ such that
$\alpha_i\in\F$. 


\subsubsection{Braid group representations}

A \qc algebra $A$ of type $D$ defines a family of homomorphisms
\[\pi_\F:\BD\to\wh{A}^\times\]
of the braid group $\BD$ to the set of invertible elements of the
completion of $A$ \wrt \fd representations. The homomorphisms
are labelled by the \mnss on $D$, and are defined as follows.

Let $\F$ be a \mns on $D$. For any $\alpha_i\in D$, choose
a \mns $\G_i$ such that $\alpha_i\in\G_i$ and set
\[\pi_\F(S_i)=\Phi_{\F\G_i}\cdot S^A_i\cdot\Phi_{\G_i\F}\]
The assignment $S_i\to\pi_F(S_i)$ is independent of the choice
of $\G_i$, and extends to a homomorphism $\BD\to\wh{A}^\times$
with the following properties
\begin{enumerate}
\item If $\alpha_i\in\F$, then $\pi_\F(S_i)=S^A_i$.
\item If $\G$ is another \mns on $D$ then, for any $b\in\BD$,
\[\pi_\G(b)=\Phi_{\G\F}\cdot \pi_\F(b)\cdot \Phi_{\F\G}\]
so that $\pi_\F$ and $\pi_\G$ are canonically equivalent.
\item If $\{\Psi_\F\}:A\rightarrow A'$ is a morphism of
quasi--Coxeter algebras, then for any \mns $\F$ and $b
\in\BD$,
$\Psi_\F(\pi^A_\F(b))=\pi^{A'}_\F(b)$. 
In particular, isomorphic quasi--Coxeter algebras yield
equivalent representations of $\BD$.
\end{enumerate}

\subsection{The Casimir connection} 

\subsubsection{}

Let $\g$ be a complex, simple Lie algebra, $\h\subset\g$
a Cartan subalgebra, and $\sfPhi=\{\alpha\}\subset\h^*$
the corresponding root system. Set
\[\hreg=\h\setminus\bigcup_{\alpha\in\sfPhi}\Ker(\alpha)\]
Choose root vectors $x_\alpha\in\g_\alpha$, $\alpha\in\sfPhi$,
such that $(x_\alpha,x_{-\alpha})=1$, where $(\cdot,\cdot)$
is a fixed non--degenerate, invariant bilinear form on $\g$,
and let
\[\Kalpha=x_\alpha x_{-\alpha}+x_{-\alpha}x_\alpha\in U\g\]
be the (truncated) Casimir operator of the three--dimensional
subalgebra $\sl{2}^\alpha\subseteq\g$ corresponding to $\alpha$.

Let $V$ be a \fd representation of $\g$, and $\VV$ the
holomorphically trivial vector bundle on $\hreg$ with fibre
$V$. The Casimir connection of $\g$ is the holomorphic
connection on $\VV$ given by
\begin{equation}\label{eq:nablak}
\nablak=
d-\frac{\nablah}{2}
\sum_{\alpha\in\sfPhi_+}\frac{d\alpha}{\alpha}\cdot\Kalpha
\end{equation}
where $\nablah\in\IC$ is a deformation parameter, and
$\alpha$ ranges over a chosen system of positive roots
$\sfPhi_+$.\footnote{the connection is independent of
the choice of $\sfPhi_+$} The connection $\nablak$ has
logarithmic singularities on the root hyperplanes and is
flat for any $\nablah$ \cite{MTL,DC,FMTV}.

The flat vector bundle $(\VV,\nablak)$ defines a one--parameter
family of monodromy representations of the pure braid group
$P_W=\pi_1(\hreg)$ on $V$ deforming the trivial action. We
explain in the next two paragraphs how to twist $(\VV,\nablak)$
so that it defines representations of the full braid group $B_W=
\pi_1(\hreg/W)$.

\subsubsection{Tits extension \cite{Ti}}

Let $W\subset GL(\h)$ be the Weyl group of $\g$. It is well--known
that $W$ does not act on $V$ in general, but that the triple exponentials
\[
\wt{s}_i =\exp(e_{\alpha_i})\exp(-f_{\alpha_i})\exp(e_{\alpha_i})
\]
corresponding to the simple roots $\alpha_1,\ldots,\alpha_r$
in $\sfPhi_+$ and a choice of root vectors $e_{\alpha_i},f_{\alpha_i}
$ such that $[e_{\alpha_i},f_{\alpha_i}]=h_{\alpha_i}$ give rise
to an action of an extension $\wt{W}$ of $W$ by the sign group
$\IZ_2^r$. 

As an abstract group, $\wt{W}$ is presented on generators
$\wt{s}_i$ subject to the relations
\begin{align*}
\underbrace{\wt{s}_i\wt{s}_j\cdots}_{m_{ij}}&=\underbrace{\wt{s}_j\wt{s}_i\cdots}_{m_{ij}}\\[-1.1ex]
\wt{s}_i^4&=1\\
\wt{s}_i^2\wt{s}_j^2&=\wt{s}_j^2\wt{s}_i^2\\
\wt{s}_i\wt{s}_j^2\wt{s}_i^{-1}&=\wt{s}_j^2(\wt{s}_i)^{-2\alpha_i(\cor{j})}
\end{align*}
for any $i\neq j$, where $m_{ij}$ is the order of product of simple
reflections $s_is_j$ in $W$. It is therefore independent of the
choices made, and a quotient of the braid group $B_W$ of
$W$. Its action on $V$ factors though that of $N_H\subset G$,
where $G$ is the simply--connected complex Lie group with
Lie algebra $\g$, $H\subset G$ its maximal torus with Lie algebra
$\h$, and $N_H$ the normaliser of $H$ in $G$. The image of
$\wt{W}$ in $N_H$, and therefore in $GL(V)$, depends upon the
choices of the root vectors, but different choices lead to subgroups
of which are canonically conjugate under an element in $H$.

\subsubsection{Twisting of $(\VV,\nablak)$ \cite{MTL}}

The flat vector bundle $(\VV,\nablak)$ is equivariant under the
Tits extension $\wt{W}$, and may be twisted into a $W$--equivariant,
flat vector bundle $(\wt{\VV},\wtnablak)$ on $\hreg$ as follows.
Let $\wt{\h}\reg\stackrel{p}{\rightarrow}\hreg$ be the universal
cover of $\hreg$ and $\hreg/W$. Since $\wt{W}$ is a quotient
of the braid group $B_W$, 
the latter acts on the flat vector bundle
$p^*(\VV,\nablak)$ on $\wt{\h}\reg$.
By definition, $(\wt{\VV},\wtnablak)$ is the quotient 
$p^*(\VV,\nablak)/P_W$, where $P_W=\pi_1(\hreg)$, 
and carries a residual action of $W=B_W/P_W$.
\Omit{First explanation could be perhaps removed?}

The flat vector bundle $(\wt{\VV},\wtnablak)$ may alternatively
be described as follows. Let $\IZ_2^r\cong Z\subset\wt{W}$ be
the subgroup generated by the elements $\{\wt{s}_i^2\}_{i=1}^r$. Since
$Z$ is a quotient of the pure braid group $P_W$, it is the deck
group of a Galois cover $\pi:\wt{\h}^Z\reg\to\hreg$ which, via
the projection $\hreg\to\hreg/W$, is also Galois cover of $\hreg/W$
with deck group $\wt{W}$. Let $\Z\to\hreg$ be the direct image
of the trivial line bundle over $\wt{\h}\reg^Z$.
Then, $\Z$ is a flat, $\wt{W}$--equivariant vector bundle of rank
$2^r$ over $\hreg$, and $(\wt{\VV},\wtnablak)$ is the bundle of
$Z$--coinvariants
\[\wt{\VV}=[\VV\otimes\Z]_Z\]

\subsection{}

The monodromy of $(\wt{\VV},\wtnablak)$, which we shall
abusively refer to as the monodromy of the Casimir connection,
yields a one--parameter family of representations $\mu_V^\nablah$
of $B_W$ on $V$ which is obtained as follows.
Fix a base point $\wt{x}_0\in\wt{\h}\reg$, and let $x_0,[x_0]$
be its images in $\hreg$ and $\hreg/W$ respectively.
The braid group $\pi_1(\hreg/W;[x_0])$ acts on
fundamental solutions $\Psi:\wt{\h}\reg\to GL(V)$ of
$p^*\nablak$ by $b\bullet\Psi(\wt{x})=b\cdot\Psi(b^{-1}\cdot\wt{x})$.
If $\Psi$ is a given fundamental solution, then $\mu_\Psi^\nablah
(b)=\Psi^{-1}\cdot b\bullet\Psi$ is locally constant function
with values in $GL(V)$ and the required monodromy.

It will be shown below that the formal Taylor series of $\mu_V^\nablah$
at $\nablah=0$ arises from an appropriate quasi--Coxeter structure on
$\Ug\fml$.

\subsection{The Casimir connection and \qc structure on $\Ug\fml$}

\subsubsection{Fundamental solutions of $\nablak$ \cite{DCP2,Ch1,Ch2}}
\label{sss:DCP soln}

Let $D$ be the Dynkin diagram of $\g$ relative to simple roots $\alpha_
1,\ldots,\alpha_r$. An {\it adapted family} of $\h^*$ is a collection $\beta
=\{x_B\}\subset\h^*$ labelled by the connected subdiagrams $B\subseteq
D$ such that, for any maximal nested set $\F$ on $D$, and $B\in\F$, the
elements $\{x_C\}_{\F\ni C\subseteq B}$ form a basis of the subspace
$\h^*_B$ of $\h^*$ spanned by the simple roots labelled by the vertices
of $B$. An example of such a family may be obtained by taking $x_B$
to be the sum of the positive roots in the root subsystem $\sfPhi_B\subset
\sfPhi$ corresponding to $B$.

For any \mns $\F$, let $\U$ denote the affine space $\IC^\F$
with coordinates $\{u_B\}_{B\in\F}$, and consider the map
$\rho_\F:\U\to\h$ given in the coordinates $\{x_B\}_{B\in\F}$
on $\h$ by
\[x_B=\prod_{\F\ni C\supseteq B}u_C\]
$\rho_\F$ is a birational map, with inverse
\[u_B=\left\{\begin{array}{ll}
x_D			&\text{if $B=D$}\\
x_B/x_{C(B)}	&\text{otherwise}
\end{array}\right.\]
where $C(B)\in\F$ is the unique minimal element properly
containing $B$. It restricts to an isomorphism between
the complement in $\U$ of the coordinate hyperplanes
$\{u_B=0\}$, and that in $\h$ of the hyperplanes
$\{x_B=0\}$. Moreover, $\rho_\F$ maps $\{u_B=0\}$
into the subspace $\sfPhi_B^\perp$.

The pull--back to $\U$ of a root $\alpha\in\sfPhi$ has the following
expression in terms of the coordinates $\{u_B\}$. Let $B\in\F$
be the unique minimal element such that $\alpha\in\sfPhi_B$.
Then, there are complex numbers $\{a_{B'}\}_{\F\ni B'\subseteq
B}$, with $a_B\neq 0$, such that
\[\alpha=
\sum a_{B'} x_{B'} =
a_B x_B\left(1+\sum_{\F\ni B'\subsetneq B}\frac{a_{B'}}{a_B} \frac{x_{B'}}{x_B}\right)=
a_B\cdot \prod_{\F\ni C\supseteq B}u_C\cdot P_\alpha\]
where $P_\alpha$ is a polynomial in the variables $\{u_{B'}\}_{\F\ni
B'\subsetneq B}$ such that $P_\alpha(0)=1$.

Set $\U_\F=\U\setminus\bigcup_{\alpha\in\sfPhi_+}\{P_\alpha=0\}$.
The pull--back of the connection $\nablak$ to $\U_\F$ has logarithmic
singularities on the divisor $\prod_{B\in\F}\{u_B=0\}$, with residue
on the hyperplane $u_B=0$ given by $\half{\nablah}\K_B$, where
\[\K_B=\sum_{\alpha\in\sfPhi_B\cap\sfPhi_+}\Kalpha\]
Let $p_\F\in\U_\F$ be the point with coordinates $u_B=0$, $B\in\F$. 
Then, for every simply--connected open set $\V\subset\U_\F$
containing $p_\F$, there is a unique holomorphic function $H_\F:
\V\to U\g\fml^o$ such that $H_\F(p_\F)=1$ and, for any determination
of the logarithm, the function
\[\Psi_\F=H_\F\cdot\prod_{B\in\F}u_B^{\frac{\nablah}{2} \K_B}\]
is a solution of $\nablak\Psi_F=0$.
The fundamental solution $\Psi_\F$ has good asymptotics near
$0\in\h$, when the latter is approached so that, for any $B\subset
C\in\F$, the roots in $\sfPhi_B$ go to zero faster than those in
$\sfPhi_C$.

\subsubsection{\DCP associators \cite{DCP2}}

Assume that the elements $\{x_B\}_{B\subseteq D}$ are real
and positive on the fundamental chamber
\[\C=\{t\in\h|\,\alpha(t)>0,\,\alpha\in\sfPhi_+\}\]
For any \mns $\F$, let $\V_\F\subset\U_\F$ be the complement
of the real, codimension one semialgebraic subvarieties $\{x_B\leq 0\}$,
$B\in\F$.
The preimage of the chamber $\C$ lies in $\V_\F$ since $x_B>0$ on $\C$. We shall 
henceforth only consider the standard determination of the logarithm,
so that $\log(x_B),\log(u_B)$, $B\in\F$ are well--defined and single--valued
on $\V_\F$. The fundamental solution $\Psi_\F$
is single--valued on the intersection of a neighborhood of $p_\F$ in $\U
_\F$ with $\V_\F$. Since $p_\F$ lies in the closure of $\C$, $\Psi_\F$
may be continued to a single--valued solution on $\C$.

Let now $\F,\G$ be two maximal nested sets. The \DCP associator
$\Phi_{\G\F}$ is the element of $U\g\fml$ defined by
\[\Phi_{\G\F}=
(\Psi_{\G}(y))^{-1}\cdot\Psi_\F(y)\]
for any $y\in\C$. 

\subsubsection{Quasi--Coxeter structure}\label{sss:DCP qC}

\begin{prop}\label{pr:QKnabla}
Set $\hbar=\pi\iota\nablah$. Then,
\begin{enumerate}
\item The associators $\Phi_{\G\F}$ and local monodromies
\begin{equation}\label{eq:Snablak}
S_{i,\kappa}^\nabla=
\wt{s}_i\cdot
\exp\left(\hbar/2\cdot\K_{\alpha_i}\right)
\end{equation}
endow $U\g\fml$ with a \qca structure $\Qnablak$ of type $D$.
\item For any \fd $\g$--module $V$ and \mns $\F$, the
representation
\[\pi_\F:B_W\to GL(V\fml)\]
obtained from the \qc structure $\Qnablak$ coincides with the
monodromy of $(\wt{\VV},\wt{\nabla}_\kappa)$ expressed in
the fundamental solution $\Psi_\F$.
\end{enumerate}
\end{prop}

\subsubsection{Modification of $\Q^\nabla_\kappa$}\label{sss:modification}

It will be convenient to alter the local monodromies of the \qc
structure $\Q^\nabla_\kappa$ as follows.

For any root $\alpha\in\sfPhi$, let $a_\alpha\in\IC[h_\alpha]\fml$
be an invertible element such that, for any $w\in W$, $wa_\alpha
=a_{w\alpha}$. It is easy to see that the local monodromies $(S_
{i,\kappa}^\nabla)^a=S_{i,\kappa}^\nabla\cdot a_{\alpha_i}$ satisfy
the braid relations of \ref{sss:qc} \wrt the associators $\Phi_{\G\F}$.
Indeed, since the $\Phi_{\G\F}$ are of weight zero, these relations
amount to checking that the following holds for any $i\neq j\in\bfI$
\[s_i(a_{\alpha_i})\cdot s_is_j(a_{\alpha_j})\cdot s_is_js_i(a_{\alpha_i})\cdots=
s_j(a_{\alpha_j})\cdot s_js_i(a_{\alpha_i})\cdot s_js_is_j(a_{\alpha_j})\cdots\]
where each side has $m_{ij}$ terms. This identity holds since both
sides are equal to the product $\prod_{\alpha\in\sfPhi_{ij}\cap\sfPhi_-}
a_\alpha$, where $\sfPhi_{ij}\subset\sfPhi$ is the rank two root
system generated by $\alpha_i$ and $\alpha_j$.

Set now $a_\alpha=\exp(\half{\hbar}\cdot\half{(\alpha,\alpha)}h_
\alpha^2)$, so that $(S_{i,\kappa}^\nabla)^a$ is given by
\begin{equation}
S^\nabla_{i,C}=
\wt{s}_i\cdot
\exp\left(\hbar/2\cdot C_{\alpha_i}\right)
\end{equation}
where $C_\alpha=\K_\alpha+\half{(\alpha,\alpha)}\cdot h_\alpha
^2$ is the Casimir operator of $\sl{2}^\alpha$. The following result
is a direct consequence of Proposition \ref{pr:QKnabla} and the
previous discussion.

\begin{corollary} The associators $\Phi_{\G\F}$ and local
monodromies $S^\nabla_{i,C}$ endow $U\g\fml$ with a
\qca structure $\Qnablac$  of type $D$.
\end{corollary}

\subsubsection{Modification of $\nablak$}

The \qc structure $\Qnablac$ encodes the monodromy of the
following connection
\begin{equation}\label{eq:nablac}
\nablac=
d-\frac{\nablah}{2}
\sum_{\alpha\in\sfPhi_+}\frac{d\alpha}{\alpha}\cdot C_\alpha
\end{equation}

More precisely, the connection $\nablac$ differs from the
Casimir connection $\nablak$ by the addition of the closed,
$S^2\h$--valued, $\wt{W}$--equivariant 1--form
\[a=
\half{\nablah}\sum_{\alpha\in
\sfPhi_+}\frac{d\alpha}{\alpha}\cdot \half{(\alpha,\alpha)}\cdot
h_\alpha^2\]
It follows that $\nablac$ is flat and may be twisted to a $W
$--equivariant flat connection $\wt{\nabla}_C$ on $\wt{\VV}$.
Since $\Psi$ is a horizontal section of $\nablak$ if, and only if
$\Psi\cdot\Theta$ is one of $\nablac$, where
\[\Theta=
\prod_{\alpha\in\sfPhi_+}\alpha^{\half{\nablah}\cdot\half{(\alpha,\alpha)}h_\alpha^2}\]
the following is a direct consequence of Proposition \ref{pr:QKnabla}.

\begin{corollary}
For any \fd $\g$--module $V$ and \mns $\F$, the representation
$\pi_\F:B_W\to GL(V\fml)$ obtained from the \qc structure $\Qnablac$
coincides with the monodromy of $(\wt{\VV},\wt{\nabla}_C)$
expressed in the fundamental solution $\Psi_\F\cdot\Theta$.
\end{corollary}

\begin{rem}
One can construct fundamental solutions $\Psi_{\F,C}$ of the connection
$\nablac$ exactly as in \ref{sss:DCP soln} and describe the monodromy
of $\nablac$ directly in terms of the associators $\Phi_{\G\F,C}=\Psi_{\G,C}
^{-1}\cdot\Psi_{\F,C}$. These, however, do not satisfy the central support
axiom of \ref{sss:weak qc} and therefore do not define a \qc structure on
$\Ug\fml$. The \qc structure $\Qnablac$ is therefore best seen as arising
from a modification of $\Qnablak$ rather than from the monodromy of
$\nablac$.
\end{rem}

\subsection{Quasi--Coxeter \qtqbas}

\subsubsection{$D$--quasibialgebras.}\label{ss:Dqba}

If $(A,\Delta,\veps)$ is a bialgebra and $n\in\IN$, we denote
by $\Delta^{(n)}:A\to A^{\otimes n}$ be the iterated coproduct
defined by  $\Delta^{(0)}=\veps$, $\Delta^{(1)}=\id$, and
\[\Delta^{(n+1)}=\Delta\otimes\id^{\otimes n-1}\circ\Delta^{(n)}\]
if $n\geq 1$. 

Let $D$ be a connected diagram. A {\it $D$--\ba} is a $D$--algebra
$A$ endowed with a bialgebra structure such that, for any connected
$B\subseteq D$, $A_B$ is a sub\ba of $A$. If $B'\subseteq B\subseteq
D$ are subdiagrams, with $B$ connected, and $n$ is an integer, we
set
\[(A^{\otimes n}_B)^{B'}=
\{\bfa\in A_B^{\otimes n}|\,[\Delta^{(n)}(a'),\bfa]=0,\,a'\in A_{B_i}\}\]
where $B_i$ ranges over the connected components of $B'$.

A {\it $D$--\qba} $(A,\{A_B\},\Delta,\veps,\{\Phi_B\},\{F_{(B;\alpha)}\})$
is a $D$--bialgebra endowed with the following additional data
\begin{itemize}
\item {\bf Associators.} For each connected subdiagram $B\subseteq
D$, an invertible element $$\Phi_B\in(A_B^{\otimes 3})^B$$
\item {\bf Structural twists.} For each connected subdiagram $B\subseteq
D$ and vertex $\alpha\in B$, a twist
$$F_{(B;\alpha)}\in(A_B^{\otimes 2})^{B\setminus\alpha}$$
\end{itemize}
satisfying the following axioms
\begin{itemize}
\item For any connected $B\subseteq D$, $(A_B,\Delta,\veps,\Phi_B)$
is a quasibialgebra.
\item For any connected $B\subseteq D$ and $\alpha\in B$,
\[(\Phi_B)_{F_{(B;\alpha)}}=\Phi_{B\setminus\alpha}\]
where $\Phi_{B\setminus\alpha}=\prod_{B'}\Phi_{B'}$, with the product
ranging over the connected components of $B\setminus\alpha$ if $B
\neq\alpha$, and $\Phi_\emptyset=1^{\otimes 3}$ otherwise.
\end{itemize}

The gist of the above axioms is the following. For any subdiagram
$B\subseteq D$, let $A_B$ be the algebra generated by the $A_
{B_i}$, where $B_i$ runs over the connected components of $B$
and set $\Phi_B=\prod_i\Phi_{B_i}$. A $D$--\qba gives rise to a
family of tensor categories
\[\C_B=\Rep(A_B,\Delta,\veps,\Phi_B)\]
labelled by the subdiagrams $B\subseteq D$. Moreover, the structural
twists give rise to restriction functors $\C_B\to\C_{B'}$, $B'\subseteq B$
in the following way. For any \mns $\F$ on $D$ containing the connected
components of $B,B'$, set
\[\F_{B'B}=\{C\in\F|B'_j\subseteq C\subseteq B_i\,\text{for some $i,j$}\}\]
where $B_i,B_j'$ are the connected components of $B,B'$. Define
the twist $F_{\F_{B'B}}\in A_B^{\otimes 2}$ by
\[F_{\F_{B'B}}=
\stackrel{\longrightarrow}{\prod_{C\in\F_{B'B}}}
F_{(C;\alpha^C_\F)}\]
where the product is taken with $F_{(C_1;\alpha^{C_1}_\F)}$ written
to the left of $F_{(C_2;\alpha^{C_2}_\F)}$ whenever $C_1\subset
C_2$. This does not specify the order of the factors uniquely, but two
orders satisfying this requirement are readily seen to yield the same
product. Then,
\[\left(\Phi_B\right)_{F_{\F_{B'B}}}=\Phi_{B'}\]
so that the restriction functor $\C_B\to\C_{B'}$ is endowed with a
collection of tensor structures labelled by any such $\F$.
\comment{One option would be to use right from the start the
notation $F_\F$ but then ask for the factorisation property of
twists.}

\subsubsection{} 

A {\it weak \qc \qba} of type $D$ is a set
\[\left((A,\{A_B\},\{\Phi_{\G\F}\},\Delta,\veps,\{\Phi_B\},\{F_{(B;\alpha)}\}\right)\]
where 
\begin{itemize}
\item $(A,\{A_B\},\{\Phi_{\G\F}\})$ is a weak \qc algebra of type $D$
\item $(A,\{A_B\},\Delta,\veps,\{\Phi_B\},\{F_{(B;\alpha)}\})$ is a $D$--\qba
\end{itemize}
and, for any pair $(\G,\F)$ of \mnss on $D$, the following holds
\[F_\G\cdot\Delta(\Phi_{\G\F})=\Phi_{\G\F}^{\otimes 2}\cdot F_\F\]
where, in the notation of \ref{ss:Dqba}, $F_\F=F_{\F_{\emptyset D}}$.

Thus, in a weak \qc \qba the tensor structures on the restriction functors
$\C_B\to\C_{B'}$ are naturally isomorphic via the associators $\Phi_{\G\F}$.

\subsubsection{} 

A {\it \qcqtqba of type $D$} is a weak \qc \qba $A$ of type $D$ endowed
with the following additional data
\begin{itemize}
\item {\bf R--matrices}. For any connected $B\subseteq D$, an invertible
element $R_B\in A_B^{\otimes 2}$ such that $(A_B,\Delta,\veps,\Phi_B,
R_B)$ is a quasitriangular quasibialgebra.
\item {\bf Local monodromies}. For any vertex $\alpha_i\in D$, an invertible
element $S_i\in\wh{A}_i$ such that $(A,\{\Phi_{\G\F}\},\{S_i\})$ is a \qc
algebra of type $D$.
\end{itemize}
The $R$--matrices and local monodromies are subject to the following
compatibility relation:
\begin{itemize}
\item {\bf Coproduct identity.} For any vertex $\alpha_i\in D$, the following
holds
\[\Delta_{F_{(\alpha_i;\alpha_i)}}(S_i)
=
(R_{\alpha_i})_{F_{(\alpha_i;\alpha_i)}}^{21}\cdot
S_i\otimes S_i\]
\end{itemize}
\comment{Say what this means.}

\section{Differential twists and quasi--Coxeter structures}
\label{se:diffl twist}

We define in this section the notion of {\it differential twist} of a semisimple
Lie algebra $\g$, and show that it gives rise to a \qcqtqba structure on $U
\g\fml$, which interpolates between the \qtqba structure underlying the
monodromy of the KZ connection and the \qc structure underlying that of
the Casimir one. The corresponding structural twists arise by comparing
the asymptotics of the differential twist for $\g$ when a given coordinate
$\alpha_i$ tends to $\infty$, to the differential twist for the subalgebra of
$\g$ generated by the root vectors corresponding to the simple roots
$\alpha_j\neq\alpha_i$.

\subsection{Notation}

For any subdiagram $B\subseteq D$, let $\g_B\subseteq\g$ be the
subalgebra generated by the root subspaces $\g_{\pm\alpha_i}$,
$\alpha_i\in B$, $\sfPhi_B\subset\sfPhi$ its root system, and $\ll_B
=\g_B+\h$ the corresponding Levi subalgebra of $\g$. Denote by
\[\Omega_B=x_a\otimes x^a,
\quad 
C_B=x_a\cdot x^a
\quad\text{and}\quad
r_B=
\sum_{\alpha\in\sfPhi_B\cap\sfPhi_+}x_\alpha\otimes x_{-\alpha}\]
where $\{x_a\}_a,\{x^a\}_a$ are dual basis of $\gD$
with respect to $(\cdot,\cdot)$, and $x_\alpha\in\g_\alpha$
are root vectors such that $(x_\alpha,x_{-\alpha})=1$,
the corresponding invariant
tensor, Casimir operator and standard solution of the
modified classical Yang--Baxter equation for 
$\g_B$ respectively. Let also $\Phi_B^{\KKZ}$
be the KZ associator for $\g_B$ defined in \ref{sss:KZ}.

\Omit{Abbreviate $\sl{2}^{\alpha_i}$,
$\Omega_{\alpha_i}$ and $C_{\alpha_{i}}$ to $\sl{2}^{i}$,
$\Omega_{i}$ and $C_{i}$ respectively and let $\wt{s}_i$
be the triple exponentials \eqref{eq:triple exp}.}

\subsection{Differential twist}

Let $\C_\IR=\{t\in\h|\,\alpha_i(t)>0,\,\,\forall i\in\bfI\}\subset\h_\IR$ be the
fundamental chamber of $\h$ and $\C=\C_\IR+i\h_\IR$ its complexification.
\comment{Recheck conventions: the differential twist I construct is more
likely defined on $\h_\IR+\iota\C_\IR$. If so, say that $F_\F$ are
defined on the quadrant $\C_\IR+i\C_\IR$.}

\begin{defn}
A {\it differential twist} for $\g$ is a holomorphic map $F:\C\to \Ug^{\otimes 2}\fml^o$
such that
\begin{enumerate}
\item\label{it:veps} $\veps\otimes\id(F)=1=\id\otimes\veps(F)$.
\item\label{it:norm}  $F\equiv 1^{\otimes 2}\mod\hbar$.
\item\label{it:Phi} $(\Phi\KKZ)_F=1^{\otimes 3}$.
\item\label{it:Alt} $\Alt_{2}F=\hbar\rD\mod\hbar^{2}$.
\item\label{it:PDE} F satisfies
\begin{equation}\label{eq:2-Casimir}
d F=\half{\hbar}\sum_{\alpha\in\Phi_+}\frac{d\alpha}{\alpha}
\Bigl((\Kalpha\otimes 1+1\otimes \Kalpha)F-F\Delta(\Kalpha)\Bigr)
\end{equation}
\Omit{\item $\Theta(F)=F^{21}$ where $\Theta$ is ....}
\end{enumerate}
\end{defn}

\begin{rem}
Condition \eqref{it:PDE} is compatible with \eqref{it:veps} in that
it implies that both $\veps\otimes\id(F)$ and $\id\otimes\veps(F)$ 
satisfy
\[d e=\half{\hbar}\sum_{\alpha\in\Phi_+}\frac{d\alpha}{\alpha}
\left[\Kalpha,e\right]\]
which is consistent with $e\equiv 1$. It is also compatible with
\eqref{it:Phi}--\eqref{it:Alt} in the
following sense. Define $f:\C\to\Ug^{\otimes 2}$ by $F=1^{\otimes 2}+\hbar
f\mod\hbar^2$. It follows from \eqref{it:Phi} and the fact that $\Phi\KKZ=1^
{\otimes 3}\mod\hbar^2$, that $\dH f=0$, where
\Omit{
\[\dH f=1\otimes f-\Delta\otimes\id(f)+\id\otimes\Delta(f)-f\otimes 1\]
is the Hochschild differential.}
$\dH:\Ug^{\otimes m}\to\Ug ^{\otimes(m+1)}$ is the Hochschild differential
\[\dH a=
1\otimes a+\sum_{i=1}^m(-1)^i \id^{\otimes (i-1)}\otimes\Delta\otimes\id^{\otimes (m-i)}(a)+(-1)^{m+1}a\otimes 1\]
Moreover, \eqref{it:PDE} and \eqref{it:norm} imply that
$$df=
\half{1}\sum_{\alpha\in\Phi_+}\Bigl(\Kalpha\otimes 1+1\otimes \Kalpha-\Delta(\Kalpha)\Bigr)=
\half{1}\sum_{\alpha\in\Phi_+}\frac{d\alpha}{\alpha}\dH(\Kalpha)$$
It follows that the \Ho cohomology class of $f$ is constant on $\C$
and, by \eqref{it:Alt}, equal to $\rD$.
\end{rem}

\subsection{Compatibility with DCP associators}\label{ss:diffl DCP}

Fix henceforth a positive, adapted family $\beta=\{x_B\}_{B\subseteq D}
\subset\h^*$. For any \mns $\F$, let $\Psi_\F:\C\to\Ug\fml^o$ be the
fundamental solution of $\nablak$ corresponding to $\F$ and $\beta$,
and $\Phi_{\G\F}=\Psi_\G^{-1}\cdot\Psi_\F$ the corresponding associators.

Let $F:\C\to\Ug^{\otimes 2}\fml^o$ be a differential twist for $\g$, and set
\[F_\F=(\Psi_\F^{\otimes 2})^{-1}\cdot F\cdot\Delta(\Psi_\F^{\otimes 2})\]

\begin{lemma}\label{le:diffl DCP}
The following holds
\begin{enumerate}
\item $F_\F$ is constant on $\C$.
\item $\veps\otimes\id(F_\F)=1=\id\otimes\veps(F_\F)$.
\item $(\Phi\KKZ)_{F_\F}=1^{\otimes 3}$.
\item $F_\F=\Phi_{\F\,\G}^{\otimes 2}\cdot F_\G\cdot\Delta(\Phi_{\G\,\F})^{-1}$.
\end{enumerate}
\end{lemma}
\begin{pf}
(1) Is clear.

(2) The differential equation $\nablak\Psi_\F=0$ implies that $d\veps
(\Psi_\F)=0$. Thus, $\veps(\Psi_\F)=\veps(H_\F\cdot\prod_{B\in\F}u_B
^{\hbar/2\K_B})=\veps(H_\F)$ is constant, and therefore equal to $\veps
(H_\F(p_\F))=1$. The claim now follows from the fact that $\veps\otimes
\id(u^{\otimes 2}\cdot F\cdot \Delta(u)^{-1})=\veps(u)\Ad(u)\left(\veps\otimes
\id(F)\right)$ for any invertible elements $u\in\Ug\fml$ and $F\in\Ug\fml
^{\otimes 2}$.

(3) Follows because $\Phi_{u^{\otimes 2}\cdot F\cdot \Delta(u)^{-1}}=
\Ad(u^{\otimes 3})(\Phi_F)$ for any invertible elements $u\in\Ug\fml$ and
$F\in\Ug\fml^{\otimes 2}$, and any associator $\Phi\in\Ug\fml^{\otimes 3}$.

(4) follows from the definiton of $F_\F,F_\G$.
\end{pf}

\subsection{Notation}\label{ss:preamble}

Fix $i\in\bfI$, let $\olPhi\subset\Phi$ be the root system
generated by the simple roots $\{\alpha_j\}_{j\neq i}$,
$\olg\subset\g$ the subalgebra spanned by the root
vectors and coroots $\{x_\alpha,\alpha^\vee\}_{\alpha
\in\olPhi}$, $\olg\supset\olh\subset\h$ its Cartan
subalgebra, and $\oll=\g\oplus\h$ the corresponding
Levi subalgebra of $\g$.

The inclusion of root systems $\olPhi\subset\Phi$ gives rise to
a projection $\pi:\h\to\olh$ determined by the requirement that
$\alpha(\pi(t))=\alpha(t)$ for any $\alpha\in\olPhi$. The kernel
of $\pi$ is the line $\IC\cow{i}$ spanned by the $i$th fundamental
coweight of $\g$. 

We shall coordinatise the fibres of $\pi$ by restricting
the simple root $\alpha_i$ to them. This amounts to
trivialising the fibration $\pi:\h\to\olh$ as $\h\isom\IC
\times\olh$ via $(\alpha_i,\pi)$. The inverse of this
isomorphism is given by $(w,\olmu)\to w\cow{i}+\imath
(\olmu)$, where $\imath:\ol\h\to\h$ is the embedding
with image $\Ker(\alpha_i)$ given by
\begin{equation}\label{eq:emb i}
\imath(\ol{t})=\ol{t}-\alpha_i(\ol{t})\cow{i}
\end{equation}

\Omit{Note that $\imath$ maps the fundamental coweights
$\{\lambda_j^{\olg,\vee}\}$ of $\olg$ to their counterparts
$\{\lambda_j^{\g,\vee}\}$ for $\g$, and that it differs
from the embedding $\olh\subset\h$ given by the
inclusion $\olg\subset\g$, which maps the coroots
$\alpha^\vee$ of $\olg$ to the corresponding ones
of $\g$.}

Denote by
\begin{equation}\label{eq:K olK}
\K=\sum_{\alpha\in\Phi_+}\Kalpha\aand
\ol{\K}=\sum_{\alpha\in\ol{\Phi}_+}\Kalpha
\end{equation}
the (truncated) Casimir operators of $\g$ and $\ol{\g}$.

\subsection{Asymptotics of the Casimir connection for $\alpha_i\to\infty$}\label{se:Fuchs infty}

Retain the notation of \ref{ss:preamble}. Fix $\olmu\in\olh$,
and consider the fiber of $\pi:\h\to\olh$ at $\olmu$. Since
the restriction of $\alpha\in\Phi$ to $\pi^{-1}(\olmu)$ is
equal to $\alpha(\cow{i})\alpha_i+\alpha(\imath(\olmu))$,
the restriction of the Casimir connection $\nablac$ to
$\pi^{-1}(\olmu)$ is equal to
\[\nabla_{i,\olmu}=
d-\half{\hbar}\sum_{\alpha\in\Phi_+\setminus\olPhi}
\frac{d\alpha_i}{\alpha_i-w_\alpha}\K_\alpha\]
where $w_\alpha=-\alpha(\imath(\olmu))/\alpha(\cow{i})$.
Set
\begin{equation}\label{eq:Rmu}
R_\olmu=\max\{|w_\alpha|\}_{\alpha\in\Phi\setminus\ol{\Phi}}
\end{equation}

\begin{prop}\label{pr:Fuchs infty}\hfill
\begin{enumerate}
\item For any $\olmu\in\olh$, there is a unique holomorphic
function
\[H_\infty:\{w\in \IP^1|\,|w|>R_\olmu\}\to\Ug\fml^o\]
such that $H_\infty(\infty)=1$ and, for any determination of
$\log(\alpha_i)$, the function $\Upsilon_\infty=H_\infty(\alpha_i)
\cdot \alpha_i^{\half{\hbar}(\K-\ol{\K})}$ satisfies
\[\left(d-\half{\hbar}\sum_{\alpha\in\Phi_+\setminus\olPhi}
\frac{d\alpha_i}{\alpha_i-w_\alpha}\K_\alpha\right)\Upsilon_\infty
=\Upsilon_\infty\,d\]
\item The function $H_\infty(\alpha_i,\olmu)$ is holomorphic on
the simply--connected domain $\D_\infty\subset\IP^1\times\ol{\h}$
given by
\begin{equation}\label{eq:D infty}
\D_\infty=\{(w,\ol\mu)|\,|w|>R_\olmu\}
\end{equation}
and, as a function on $\D_\infty$, $\Upsilon_\infty$ satisfies
\[\left(d-\half{\hbar}\sum_{\alpha\in\Phi_+}\frac{d\alpha}{\alpha}
\Kalpha\right)\Upsilon_\infty=
\Upsilon_\infty
\left(d-\half{\hbar}\sum_{\alpha\in\ol{\Phi}_+}\frac{d\alpha}{\alpha}
\Kalpha\right)\]
\item $H_\infty$ satisfies $H_\infty(t\alpha_i,t\olmu)=H_\infty(\alpha_i,
\olmu)$ for any $t\in\IC^*$.
\item If $\ol{\mu}=0$, $H_\infty(\alpha_i,0)=1$, $\Upsilon_\infty
(\alpha_i,0)=\alpha_i^{\half{\hbar}(\K-\ol{\K})}$, and both commute
with $\oll$.
\end{enumerate}
\end{prop}
\begin{pf}
(1) Denote the restriction of $\alpha_i$ to $\pi^{-1}(\olmu)$ by $w$.
$H=H_\infty$ is required to satisfy the ODE
\begin{equation}\label{eq:Fuchs infty}
\frac{dH}{dw}=
\half{\hbar}\left(
\sum_{\alpha\in\Phi_+\setminus\ol{\Phi}}\frac{\Kalpha}{w-w_\alpha}H
-H\frac{\K-\ol{\K}}{w}\right)
\end{equation}
all of whose singularities are contained in the disk $\{|w|\leq R_\olmu\}$.
Writing $H=\sum_{n\geq 0}\hbar^n H_n$ yields the recursive system of
ODEs
\[\frac{dH_n}{dw}=
\half{1}\left(
\sum_{\alpha\in\Phi_+\setminus\ol{\Phi}}\frac{\Kalpha}{w-w_\alpha}H_{n-1}
-H_{n-1}\frac{\K-\ol{\K}}{w}\right)\]
where $H_{-1}=0$, together with the boundary condition $H_n(\infty)
=\delta_{n0}$, which clearly possess at most one solution. For $n=0$,
this possesses the solution $H_0\equiv 1$. For $n\geq 1$, given that
\begin{equation}\label{eq:res at infty}
\frac{1}{w-w_\alpha}=\frac{1}{w}\left(1+\frac{w_\alpha}{w-w_\alpha}\right)
\end{equation}
the equation reads
\[\frac{dH_n}{dw}=
\frac{1}{2w}\left(
[\K-\ol{\K},H_{n-1}]+
\sum_{\alpha\in\Phi
_+\setminus\ol{\Phi}}\frac{w_\alpha\Kalpha}{w-w_\alpha}H_{n-1}\right)\]
This possesses the holomorphic solution
\[H_n(w)=\int_{\Gamma^w}
\frac{1}{2t}\left(
[\K-\ol{\K},H_{n-1}(t)]+
\sum_{\alpha\in\Phi
_+\setminus\ol{\Phi}}\frac{w_\alpha\Kalpha}{t-w_\alpha}H_{n-1}(t)\right)dt\]
where $\Gamma^w$ is the ray from $\infty$ to $w$, and $w$ is assumed
to be such that $|w|>R_\olmu$, provided the integral converges at $t=\infty$. For
$n=1$, this is the case since $H_0=1$ commutes with $\K-\ol{\K}$,
so the integrand is an $O(t^{-2})$ and, for $n \geq 2$, this follows since
$H_{n-1}(t)=O(t^{-1})$.

(2) The recursive construction of $H$ clearly shows that it is a holomorphic
function of $\olmu$. The fact that $\Upsilon_\infty$ satisfies the claimed
PDE follows from the integrability of $\nablac$ (a similar argument is
given in the proof of part (2) of Theorem \ref{th:Hn}).

(3) Follows by uniqueness.

(4) If $\olmu=0$, $w_\alpha=0$ for any $\alpha$, so that \eqref{eq:Fuchs infty}
simply reads
\[\frac{dH}{dw}=
\half{\hbar}
\frac{[\K-\ol{\K},H]}{w}\]
and, by uniqueness, $H(w,0)\equiv 1$. $\Upsilon_\infty(w,0)=w^{\half{\hbar}(\K
-\ol{\K})}$ commutes with $\oll$ because $\K-\ol{\K}$ does. Indeed, since $\K=
C-t_at^a$, where $C$ is the Casimir operator of $\g$ and $\{t_a\},\{t^a\}$ are
dual bases of $\h$,
\begin{equation}\label{eq:k-olk}
\K-\ol{\K}=C-\ol{C}-\frac{(\cow{i})^2}{\|\cow{i}\|^2}
\end{equation}
which commutes with $\oll$.
\end{pf}

\begin{rem}\label{rm:strong uniqueness}
We later use, we shall need the (obvious) fact that the uniqueness of
Proposition \ref{pr:Fuchs infty} holds under the weaker assumption
that the function $H_\infty(\alpha_i,\ol{\mu})$ is holomorphic on one
of the domains
\begin{equation}\label{eq:Dinftypm}
\D_\infty^\pm=
\{(w,\olmu)\in\IC\times\olh|\,\Im w\gtrless 0,\,|w|>R_\olmu\}
\end{equation}
and is such that $H_\infty(w,\ol{\mu})\to 1$ as $w\to\infty$
with $0<\delta<|\arg w|<\pi-\delta<\pi$.
\end{rem}

\subsection{Asymptotics at $\alpha_i=\infty$ and DCP solutions}

Let $\F$ be a \mns on $D$, set $\olcalF=\F\setminus\{D\}$
and $\alpha_i=\alpha^D_\F$. Let
\[\Psi_\F:\C\to U\g\fml^o
\aand
\Psi_{\olcalF}:\olC\to U\olg\fml^o\]
be the fundamental solutions of the Casimir connection for
$\g$ and $\olg=\g_{D\setminus\alpha_i}$ corresponding to
$\F$, $\olcalF$ respectively, and a positive, adapted family
$\{x_B\}_{B\subseteq D}$. Regard $\Psi_{\olcalF}$ as being
defined on $\C$ via the projection $\pi:\h\to\olh$. The result
below expresses $\Psi_\F$ in terms of $\Psi_{\olcalF}$ and
the solution $\Upsilon_\infty$ given by Proposition \ref
{pr:Fuchs infty}.

\begin{prop}\label{pr:infty factorisation}
The following holds
\[\Psi_\F=
\Upsilon_\infty\cdot\Psi_{\olcalF}\cdot
x_D(\cow{i})^{\half{\hbar}(\K-\olK)}\]
\end{prop}
\begin{pf}
Note first that $\alpha_i$ can be expressed as a linear combination of the form
$\sum_{B\in\F}a_B x_B$. Evaluating on $\cow{i}$ shows that $a_D=x_D(\cow{i})
^{-1}$, so that
\begin{equation}\label{eq:p alphai}
\alpha_i=
\frac{x_D}{x_D(\cow{i})}\left(1+\sum_{B\in\olcalF}\frac{x_B}{x_D}\right)=
\frac{x_D}{x_D(\cow{i})}\cdot p_{\alpha_i}(u)
\end{equation}
where $p_{\alpha_i}$ is a polynomial in the variables $\{u_B\}_{B\in\olcalF}$, such
that $p_{\alpha_i}(0)=1$. By construction,
\[\begin{split}
\Upsilon_\infty\cdot\Psi_{\olcalF}
&=
H_\infty(\alpha_i,\olmu)\cdot\alpha_i^{\half{\hbar}(\K-\olK)}\cdot
H_{\olcalF}(\ul{u})\cdot\prod_{B\in\olcalF} x_B^{\half{\hbar}(\K_B-\K_{B\setminus\alpha^B_\F})}\\
&=
H_\infty(\alpha_i,\olmu)\cdot
H_{\olcalF}(\ul{u})\cdot
p_i(u)^{\half{\hbar}(\K-\olK)}\cdot
\prod_{B\in\F} x_B^{\half{\hbar}(\K_B-\K_{B\setminus\alpha^B_\F})}\cdot
x_D(\cow{i})^{-\half{\hbar}(\K-\olK)}
\end{split}\]
where the second equality follows from \eqref{eq:p alphai} and the fact that
$\K,\olK$ commute with $H_{\olcalF}$ since the latter takes values in $U\olg
\fml$ and is of weight $0$. The claimed result follows from the uniqueness
of $\Psi_\F$, provided the function $H_\infty(\alpha_i,\olmu)\cdot H_{\olcalF}
(\ul{u})\cdot p_i(u)^{\half{\hbar}(\K-\olK)}$ is holomorphic
in the neighborhood of $u=\{u_B\}_{B\in\F}=0$, and equal to 1 at $u=0$.

This is clearly true for the factor $p_i(u)^{\half{\hbar}(\K-\olK)}$. For $H_{\olcalF}$
note that, for $B\in\olcalF$, $\ol{u}_B$ is equal to $u_B$, if $B$ is not a maximal
element of $\olcalF$, and to $x_B=u_B u_D$ otherwise. Thus, $H_{\olcalF}$
is holomorphic in the neighborhood of $u=0$, and equal to $1$ at $u=0$. 
Finally, note that the coordinates of $\olmu$ are $x_B$, $B\in\olcalF$, with
$x_B=x_D\prod_{\olcalF\ni C\supset B}u_C$. Thus, by the homogeneity 
of $H_\infty$,
\[\begin{split}
H_\infty(\alpha_i,\olmu)
&=
H_\infty\left(\frac{x_D}{x_D(\cow{i})}p_{\alpha_i}(u),
\{x_D\prod _{\olcalF\ni C\supset B}u_C\}_{B\in\olcalF}\right)\\
&=
H_\infty\left(\frac{p_{\alpha_i}(u)}{x_D(\cow{i})},
\{\prod _{\olcalF\ni C\supset B}u_C\}_{B\in\olcalF}\right)
\end{split}\]
which tends to $1$ as $u\to 0$ since $H_\infty(\alpha_i,0)=1$.
\end{pf}

\begin{corollary}
For any $\alpha\in B\subseteq D$, let $\Upsilon_\infty^{B,\alpha}$ be
the solution of the Casimir connection for $\g_B$ corresponding to the
simple root $\alpha$ given by Proposition \ref{pr:Fuchs infty}. Then,
\[\Psi_\F=
\stackrel{\longleftarrow}{\prod_{B\in\F}}\Upsilon_\infty^{B,\alpha^B_\F}\cdot
\prod_{B\in\F}x_B(\cow{\alpha^B_\F})^{\half{\hbar}(\K_B-\K_{B\setminus\alpha^B_\F})}\]
where the first product is ordered with $\Upsilon_\infty^{B,\alpha^B_\F}$
to the left of $\Upsilon_\infty^{C,\alpha^C_\F}$ if $B\supset C$.
\end{corollary}

\subsection{Relative twists}\label{ss:relative twist}

Let $F$ be a differential twist for $\g$, $\alpha_i\in D$ a simple
root, and $\Upsilon_\infty$ the solution of the Casimir equations
given by Proposition \ref{pr:Fuchs infty}, where we are using the
standard determination of $\log$. Define $F_\infty:\C\to\Ug^{\otimes 2}\fml^o$ by
\[F_\infty=(\Upsilon_\infty^{\otimes 2})^{-1}\cdot F\cdot\Delta(\Upsilon_\infty)\]
Then, $F_\infty$ satisfies
\begin{enumerate}
\item $\veps\otimes\id(F_\infty)=1=\id\otimes\veps(F_\infty)$.
\item $F_\infty\equiv 1^{\otimes 2}\mod\hbar$.
\item $(\Phi\KKZ)_{F_\infty}=1^{\otimes 3}$.
\item $\Alt_{2}F_\infty=\hbar\rD\mod\hbar^{2}$.
\item
\[d F_\infty=\hbar\sum_{\alpha\in\olPhi_+}\frac{d\alpha}{\alpha}
\Bigl((\Kalpha\otimes 1+1\otimes \Kalpha)\cdot F_\infty-F_\infty\cdot\Delta(\Kalpha)\Bigr)\]
\end{enumerate}

Let $\olC$ be the complexified chamber of $\olg$, and $\olF:
\olC\to U\olg^{\otimes 2}\fml^o$ a differential twist for $\olg$.
Since the projection $\pi:\h\to\olh$ maps $\C$ to $\olC$, we
may regard $\olF$ as a function on $\C$, and define
$F'_{(D;\alpha_i)}:\C\to U\g^{\otimes 2}\fml^o$ by
\[F'_{(D;\alpha_i)}=\olF^{-1}\cdot F_\infty\]

\begin{prop}\label{pr:relative twist}
$F'_{(D;\alpha_i)}$ satisfies the following properties
\begin{enumerate}
\item $\veps\otimes\id(F'_{(D;\alpha_i)})=1=\id\otimes\veps(F'_{(D;\alpha_i)})$.
\item $F'_{(D;\alpha_i)}\equiv 1^{\otimes 2}\mod\hbar$.
\item $(\Phi\KKZ)_{F'_{(D;\alpha_i)}}=\Phi_{D\setminus\alpha_i}$.
\item $\Alt_{2}F'_{(D;\alpha_i)}=\hbar\left(\rD-r_{D\setminus\alpha_i}\right)\mod\hbar^{2}$.
\item
\[d F'_{(D;\alpha_i)}=\hbar\sum_{\alpha\in\olPhi_+}\frac{d\alpha}{\alpha}
[\Delta(\Kalpha),F'_{(D;\alpha_i)}]\]
\item If $F'_{(D;\alpha_i)}$ is invariant under $\olg$, then it is constant
on $\C$.
\end{enumerate}
\end{prop}

\begin{pf}
(1), (2), and (4) are obvious. For (3), note that since $\Phi_{u^{\otimes 2}\cdot
F\cdot \Delta(u)^{-1}}=\Ad(u^{\otimes 3})(\Phi_F)$ for any invertible elements
$u\in\Ug\fml$ and $F\in\Ug\fml^{\otimes 2}$, and any associator $\Phi\in\Ug
\fml^{\otimes 3}$, we get
$$(\Phi_D)_{\olF\cdot F'_{(D;\alpha_i)}}=
(\Phi_D)_{F'_\infty}=
1^{\otimes 3}$$
whence
\[(\Phi_D)_{F'_{(D;\alpha_i)}}=
((\Phi_D)_{\olF\cdot F'_{(D;\alpha_i)}})_{\olF^{-1}}=
(1^{\otimes 3})_{\olF^{-1}}=
\Phi_{D\setminus\alpha_i}\]

(5) follows from the PDEs satisfied by $F_\infty$ and $\olF$.

(6) follows from (5). 
\end{pf}

\subsection{Centraliser property}

Let $\{F_B\}$ be a collection of differential twists for the subalgebras $\g_B
\subset\g$, where $B$ is a subdiagram of $D$, such that if $B$ has connected
components $\{B_i\}$, then $F_B=\prod_i F_{B_i}$.

\begin{defn}
The collection $\{F_B\}$ has the {\it centraliser property} if, for any $\alpha
\in B\subseteq D$, the relative twist $F_{(B,\alpha)}$ defined in \ref{ss:relative twist}
is invariant under $\g_{B\setminus\alpha}$ (and in particular constant).
\end{defn}

\subsection{Factorisation}\label{eq:factorisation}

Let $\{F_B\}_{B\subseteq D}$ be a collection of differential twists with the
centraliser property. For any $\alpha_i\in B\subseteq D$, set
\[F_{(B;\alpha_i)}=
\left(x_B(\cow{i})^{-\half{\hbar}(\K_B-\K_{B\setminus\alpha_i})}\right)^{\otimes 2}
\cdot F'_{(B;\alpha_i)}\cdot
\Delta\left(x_B(\cow{\alpha_i})^{\half{\hbar}(\K_B-\K_{B\setminus\alpha_i})}\right)
\]
where $F'_{(B;\alpha_i)}\in\Ug_B^{\otimes 2}\fml^o$ is the relative twist
defined in \ref{ss:relative twist}, and $\{x_B\}_{B\subseteq D}$ is a positive,
adapted family. The (constant) twist $F_{(B;\alpha_i)}$ is invariant under
$\g_{B\setminus\alpha_i}$, and has the properties (1)--(4) given in
Proposition \ref{pr:relative twist}.

\begin{lemma}\label{le:factorisation}
Let $\F$ be a \mns on $D$, and $F_\F$ the twist defined in \ref{ss:diffl DCP}.
Then, the following holds
\[F_\F=\stackrel{\longrightarrow}{\prod_{B\in\F}}F_{(B;\alpha^B_\F)}\]
where the product is taken with $F_{(B;\alpha^B_\F)}$ to the right of
$F_{(C;\alpha^C_\F)}$ if $B\supset C$.
\end{lemma}
\begin{pf}
By definition of $F_\F$,
\[\begin{split}
F_\F
&=
(\Psi_\F^{\otimes 2})^{-1}\cdot F\cdot \Delta(\Psi_\F)\\
&=
(x_D(\cow{i})^{-\half{\hbar}(\K_D-\K_{D\setminus\alpha^D_\F})})^{\otimes 2}\cdot
(\Psi_{\olcalF}^{\otimes 2})^{-1}\cdot
(\Upsilon_\infty^{\otimes 2})^{-1}\cdot
F\cdot
\Delta(\Upsilon_\infty)\cdot
\Delta(\Psi_{\olcalF})\cdot
\Delta\left(x_D(\cow{i})^{\half{\hbar}(\K_D-\K_{D\setminus\alpha^D_\F})}\right)\\
&=
(x_D(\cow{i})^{-\half{\hbar}(\K_D-\K_{D\setminus\alpha^D_\F})})^{\otimes 2}\cdot
(\Psi_{\olcalF}^{\otimes 2})^{-1}\cdot
\olF\cdot
F'_{(D;\alpha^D_\F)}\cdot
\Delta(\Psi_{\olcalF})\cdot
\Delta\left(x_D(\cow{i})^{\half{\hbar}(\K_D-\K_{D\setminus\alpha^D_\F})}\right)\\
&=
(\Psi_{\olcalF}^{\otimes 2})^{-1}\cdot
\olF\cdot
\Delta(\Psi_{\olcalF})\cdot
(x_D(\cow{i})^{-\half{\hbar}(\K_D-\K_{D\setminus\alpha^D_\F})})^{\otimes 2}\cdot
F'_{(D;\alpha^D_\F)}\cdot
\Delta\left(x_D(\cow{i})^{\half{\hbar}(\K_D-\K_{D\setminus\alpha^D_\F})}\right)\\
&=
\olF_{\olcalF}\cdot
F_{(D;\alpha^D_\F)}
\end{split}\]
where the second equality follows by Proposition \ref{pr:infty factorisation},
the third one by definition of $F_{(D;\alpha^D_\F)}$, and the fourth one
from the fact that $F'_{(D;\alpha^D_\F)}$ is invariant under $\g_{D\setminus
\alpha^D_\F}$ by the centraliser property, and therefore commutes with
$\Delta(\Psi_{\olcalF})$, and the fact that, by \eqref{eq:k-olk}, $\K_D-\K
_{D-\alpha^D_\F}$ commutes with $\g_{D\setminus\alpha^D_\F}$. The
result now follows by induction.
\end{pf}

\subsection{Quasi--Coxter \qtqba structure}

For any connected subdiagram $B\subseteq D$, let
$\left(\Ug_B\fml,\Delta_{0},\Phi^{\KKZ}_B,R^{\KKZ}_B\right)$
be the \qtqba structure underlying the monodromy of the KZ connection
for $\g_B$ (see \ref{sss:KZ}). Let also
\[\left(\Ug\fml,\{\Ug_B\fml\},\{\Sicnabla{i}\},\{\Phi_{\G\F}\}\right)\]
be the \qc structure underlying the monodromy of the (untruncated)
Casimir connection $\nablac$ (see \ref{sss:DCP qC}--\ref{sss:modification}),
relative to a choice of positive, adapted
family $\{x_B\}_{B\subseteq D}$. The following is the main result of this section

\begin{thm}\label{th:welding}
Let $\{F_B\}_{B\subseteq D}$ be differential twists for $\{\g_B\}
_{B\subseteq D}$ possessing the centraliser property, and such
that
\begin{itemize}
\item $F_B$ is of weight zero for any $B\subseteq D$.
\item The following holds for any simple root $\alpha_i$, 
\[\Theta^{\otimes 2}(F_{\alpha_i})=F_{\alpha_i}^{21}\]
where $\Theta$ is an automorphism of $\g$ acting by $-1$
on $\h$.
\end{itemize}
Define the relative twists $\{F_{(B;\alpha)}\}_{\alpha\in B\subseteq D}$
as in \ref{eq:factorisation}. Then,
\[\left(
\Ug\fml,\{\Ug_B\fml\},\{S_{i,C}\},\{\Phi_{\G\F}\},
\Delta_{0},\{\Phi^{\KKZ}_B\},\{R^{\KKZ}_B\},\{F_{(B;\alpha_i)}\}
\right)\]
is a \qcqtqba such that
\begin{align*}
\Sicnabla{i}	&= \wt{s}_{i}\cdot\exp(\hbar/2\cdot C_{\alpha_i})\\
\Phi_B^{\KKZ}	&= 1^{\otimes 3}\mod\hbar^2\\
R_B^{\KKZ}	&= \exp(\hbar\cdot\Omega_B)\\
\Alt_{2}F_{(B;\alpha_i)}&=
\hbar\cdot(r_B-r_{B\setminus\alpha_i})
\medspace\mod\hbar^{2}
\end{align*}
and $\Phi_{\G\F}$, $F_{(B;\alpha_i)}$ are
of weight 0.
\end{thm}
\begin{pf}
The identity $(\Phi^{\KKZ}_B)_{F_{(B;\alpha)}}=\Phi^{\KKZ}_
{B\setminus\alpha}$ was checked in Proposition \ref{pr:relative twist}.
The identity $F_\G=\Phi_{\G\F}^{\otimes 2}\cdot F_\F\cdot\Delta_0
(\Phi_{\F\G})$ follows from Lemma \ref{le:diffl DCP} and the
factorisation given by Lemma \ref{le:factorisation}. Given
that
\[\Delta_0(S_{i,C})=
e^{\hbar\Omega_{\alpha_i}}\cdot S_{i,C}\otimes S_{i,C}=
R^{\KKZ}_{\alpha_i}\cdot S_{i,C}\otimes S_{i,C}\]
the coproduct identity
\[(\Delta_0)_{F_{(\alpha_i;\alpha_i)}}(S_{i,C})
=
(R^{\KKZ}_{\alpha_i})_{F_{(\alpha_i;\alpha_i)}}^{21}\cdot
S_{i,C}\otimes S_{i,C}\]
is readily seen to be equivalent to $\Ad(S_{i,C}^{\otimes 2})
(F_{(\alpha_i;\alpha_i)})=F_{(\alpha_i;\alpha_i)}^{21}$ and
therefore to $\Ad(\wt{s}_{i})(F_{(\alpha_i;\alpha_i)})=
F_{(\alpha_i;\alpha_i)}^{21}$ since $C_{\alpha_i}$ is
central in $U\sl{2}^{\alpha_i}$. Since the restriction to
$\sl{2}^{\alpha_i}$ of $\Theta$ and $\Ad(\wt{s}_i)$ differ
by $\Ad(\exp(t))$ for some $t\in\h$, and $F_{(\alpha_i;
\alpha_i)}$ is of weight $0$, the coproduct identity is therefore
equivalent to $\Theta^{\otimes 2}(F_{(\alpha_i;\alpha_i)})=
F_{(\alpha_i;\alpha_i)}^{21}$. This in turn follows
from the assumption on $F_{\alpha_i}$ and the
fact that
\[F_{(\alpha_i;\alpha_i)}=
\alpha_i^{-\half{\hbar}(\K_{\alpha_i}\otimes 1+1\otimes\K_{\alpha_i})}\cdot
F_{\alpha_i}\cdot
\alpha_i^{\half{\hbar}\Delta_0(\K_{\alpha_i})}\]
\end{pf}

\section{The fusion operator}\label{se:fusion}  %

In this section, we construct a {\it fusion operator}, that is a joint solution
$J(z,\mu)$ of the coupled KZ--Casimir equations on $n=2$ points, with
prescribed asymptotics when $z=z_1-z_2\to\infty$. The coupling gives
rise to an irregular singularity in $z$, and due care needs to be taken to
deal with the corresponding Stokes phenomena. We shall in fact construct
a multicomponent version of the fusion operator, which solves the
coupled Casimir and KZ equations in an arbitrary number of points.
The differential twist for $\g$, which gives rise to a \qcqtqba structure 
on $U\g\fml$, as explained in Section \ref{se:diffl twist}, will be obtained
by taking the asymptotics of $J(z,\mu)$ when $z\to 0$.\footnote{The
name fusion operator originates from the work of Etingof
and Varchenko \cite{EV}, where a representation theoretic construction
is given for the the fusion operator of the loop algebra $\wh{\g}$. The
latter is a joint solution of the trigonometric KZ equations and dynamical
difference equations (a difference analogue of the Casimir equations),
and should degenerate to our $J(z,\mu)$.}

\subsection{The joint system}

For any $n\geq 2$, let
\[\gX_n=\IC^n\setminus\bigcup_{1\leq i<j\leq n}\{z_i=z_j\}\]
be the configuration space of $n$ ordered points in $\IC$. Consider the following
connection on the trivial vector bundle over $\gX_n\times\hreg$ with fibre $\Ug^
{\otimes n}\fml$
\begin{equation}\label{eq:DCKZ}
\nabla=d
-\hbar\sum_{1\leq i<j\leq n}\frac{d(z_i-z_j)}{z_i-z_j}\Omega_{ij}
-\frac{\hbar}{2}\sum_{\alpha\in\Phi_+}\frac{d\alpha}{\alpha}
\Delta^{(n)}(\Kalpha)-d\sum_{i=1}^n z_i\ad(\mu\ii)
\end{equation}

\noindent where $\hbar$ is a formal parameter, $\Omega_{ij}=\sum_a X_a\ii
{X^a}^{(j)}$, where $X\ii=1^{\otimes(i-1)}\otimes X\otimes 1^{\otimes (n-i)}$,
$\Delta^{(n)}:\Ug\to\Ug^{\otimes n}$ is the iterated coproduct and $\mu$ is
the embedding $\h\to\Ug$. The above connection is flat \cite{FMTV}. We give
an alternative proof of this in Appendix \ref{app:joint}.

\subsection{The change of variables}

We wish to construct horizontal sections of $\nabla$ having prescribed asymptotics
as $z_i-z_j\to\infty$ for any $i\neq j$. Consider to this end the change of variables
given by the map
\[\rho:\IC^\times\times\gX_n\to\gX_n,\qquad
\rho(\zeta;\zeta_1,\ldots,\zeta_n)=(\zeta\zeta_1,\ldots,\zeta\zeta_n)\]
\noindent
Since $d\log(z_i-z_j)=d\log\zeta+d\log(\zeta_i-\zeta_j)$, the pull--back of $\nabla$
under $\rho$ is equal to $\nabla_\zeta+\ul{\nabla}$, where
\begin{equation}\label{eq:nabla z}
\nabla_\zeta=
d_\zeta-\left(\sum_{i=1}^n\zeta_i\ad(\mu\ii)+\frac{\hbar\Omega}{\zeta}\right)d\zeta
\end{equation}
where $d_\zeta$ is the de Rham differential \wrt $\zeta$, $\Omega=\sum_{i<j}
\Omega_{ij}$, and
\[\ul{\nabla}=
\ul{d}-\hbar\negthickspace\negthickspace\sum_{1\leq i<j\leq n}d\log(\zeta_i-\zeta_j)\,\Omega_{ij}
-\frac{\hbar}{2}\sum_{\alpha\in\Phi_+}\frac{d\alpha}{\alpha}\,\Delta^{(n)}(\Kalpha)
-\zeta\ul{d}\sum_{i=1}^n\zeta_i\ad(\mu\ii)\]
where $\ul{d}$ is the de Rham differential \wrt $(\zeta_1,\ldots,\zeta_n)\in\gX_n$,
and $\mu\in\hreg$.

\subsection{}

\begin{lemma}\label{le:proj cent}
If $(\ul{\zeta},\mu)\in\gX_n\times\hreg$, the projection $[\Omega_{ij}]$
of $\Omega_{ij}$ onto the kernel of $\sum_i\zeta_i\ad(\mu\ii)$ is
\[\Omega_{ij}^\h=\sum_{a=1}^r(t_a)\ii(t^a)\jj\]
where $\{t_a\},\{t^a\}$ are dual bases of $\h$.
\end{lemma}
\begin{pf}
By definition, 
$\Omega_{ij}=\Omega_{ij}^\h+\sum_{\alpha\in\Phi} x_\alpha\ii x_{-\alpha}\jj$,
where $x_{\pm\alpha}\in\g_{\pm\alpha}$ are root vectors such that $(x_\alpha,x_
{-\alpha})=1$. Since $\ad(\mu^{(k)})\,x_\alpha\ii x_{-\alpha}\jj=\alpha(\mu)(\delta_{ki}
-\delta_{kj})x_\alpha\ii x_{-\alpha}\jj$, we have
\[
\sum_k\zeta_k\ad(\mu\kk)\,x_\alpha\ii x_{-\alpha}\jj=
\alpha(\mu)(\zeta_i-\zeta_j)x_\alpha\ii x_{-\alpha}\jj
\]
from which the conclusion follows.
\end{pf}

\subsection{Holomorphic fundamental solutions of $\nabla_\zeta$}

The connection \eqref{eq:nabla z} has an irregular singularity of Poincar\'e
rank one at $\zeta=\infty$. Let $\IH_\pm=\{\zeta\in\IC|\,\Im\zeta\gtrless 0\}$.
Let $\Omega^\h=\sum_{1\leq i<j\leq n}\Omega_{ij}^\h$.

\Omit{Let $H_0$ centralise $\sum_k\zeta_k\mu\kk$ and $\Omega^\h$.
\comment{If want $H_\pm$ to take values in $\A$, need the coefficient
of $\hbar^0$ to be a multiple of $1$, so asking that it centralises $\sum
\zeta_i\ad(\mu\ii)$ and $\Omega^\h$ is superfluous.}}
\Omit{Let $\gX_n^\IR=\gX_n\cap\IR^N$
and $\hreg^\IR=\{\mu\in\hreg|\alpha(\mu)\in\IR\text{ for any }\alpha\in\Phi\}$.
\comment{may want to have both the
$\zeta_i$ and $\mu$ purely imaginary, then conditions on $\zeta$ do not
change, and $z_i-z_j=\zeta(\zeta_i-\zeta_j)$ has either positive or negative
real part}}

Let $\C=\{t\in\h_\IR|\,\alpha(t)>0\text{ for any }\alpha\in\Phi_+\}$ be the
fundamental Weyl chamber of $\g$, and set $\C_n=\{\uzeta\in\gX_n\cap
\IR^n|\,\zeta_1>\zeta_2>\cdots>\zeta_n\}$. Set $\imath=\sqrt{-1}$. Let
$H^0=c1^{\otimes n}$ be a multiple of the identity.\comment{might be
easier to assume $H^0=1$ and to have a slightly separate uniqueness
result for $H^0=0$ perhaps.}\comment{{\bf Conventions.} $\mu\in\imath
\C_\h$ ('consistenly' with stability conditions), $\ul{\zeta}\in\imath\C_n$
('consistently' with choice for $\mu$), $\Im\zeta>0$. Thus, $\ul{z}=\zeta
\ul{\zeta}$ is such that $\Re(z_i-z_j)=\Re(\zeta(\zeta_i-\zeta_j))=-|\zeta
_i-\zeta_j|\Im(\zeta)<0$ if $i<j$, a domain which includes $\ul{z}\in -\C_n$,
which is consistent with labelling representations as $V_1\otimes\cdots
\otimes V_n$, and thinking of $V_i$ as sitting over the (real) point $z_i$,
which is Drinfeld's convention.}

\begin{thm}\label{th:Hn}\hfill
\begin{enumerate}
\item For any $(\ul{\zeta},\mu)\in\iCn\times\iC$, there is a unique holomorphic
function \[H_\pm:\IH_\pm\to\A\] such that $H_\pm(\zeta)$ tends to $H^0$
uniformly as $\zeta\to\infty$ in any sector of the form $|\arg(\zeta)|\in(\delta,
\pi-\delta)$,
$\delta>0$, and the $\End_{\IC\fml}(\A)$--valued function
\[\Psi_\pm(\zeta)=
H_\pm(\zeta)\,e^{\zeta\sum_i\zeta_i\ad(\mu\ii)}\,\zeta^{\hbar\Omega^\h}\]
is a fundamental solution of $\nabla_\zeta$.\comment{$H_\pm\to\ell(H_\pm)$?
YES, and $\Omega_\h\to\ell(\Omega_\h)$}
\end{enumerate}
Assume henceforth that $H^0=1$.
\begin{enumerate}
\item[(2)] As a function on $\IH_\pm\times\iCn\times\iC$, $H_\pm$
is real analytic, and satisfies
\begin{equation}\label{eq:induced}
\ul{\nabla}H_\pm=H_\pm\,\left(
\hbar\sum_{i<j}d\log(\zeta_i-\zeta_j)\,\Omega_{ij}^\h
+\frac{\hbar}{2}\sum_{\alpha\in\Phi_+}\frac{d\alpha}{\alpha}(\sum_{i=1}^n\,\Kalpha\ii)
+\zeta\ul{d}\sum_{i=1}^n\zeta_i\ad(\mu\ii)\right)
\end{equation}
\item[(3)]\label{it:rec Hn} The following holds\comment{Have switched notation by adding
the superscript $^{(n)}$}.
\[\lim_{\zeta_1\to+\imath\infty}H^{(n)}_\pm=1\otimes H^{(n-1)}
\aand
\lim_{\zeta_n\to-\imath\infty}H^{(n)}_\pm=H^{(n-1)}\otimes 1\]

\item[(4)] \comment{Something like }$H_\pm(t\zeta;\uzeta,\mu)=H_\pm(\zeta;t\uzeta,\mu)$
for any $t\in\IR_+^*$.
\item[(5)] Let $\sigma\in\Aut(\Ug^{\otimes n})$ be the algebra automorphism
given by $x_1\otimes\cdots\otimes x_n\to x_n\otimes\cdots\otimes x_1$. Then $\sigma(H_\pm)=
\Theta^{\otimes n}(H_\pm)$, where $\Theta\in\Aut(\g)$ is any \comment{(any?)} involution
acting as $-1$ on $\h$.
\end{enumerate}
\end{thm}

\noindent
The proof of Theorem \ref{th:Hn} is given in \S \ref{ss:fusion begin} -- \S \ref{ss:fusion end},
which occupy the rest of this section.

\subsection{The multicomponent fusion operator}\label{ss:multicomponent}

Let 
\[\J\nn_\pm:\imath\C_n\times\imath\C\to\End_{\IC\fml}(U\g^{\otimes n}\fml_{o})\]
be the real analytic function given by\comment{remove the dependence of
$\J\nn$ from $\zeta$ completely by using the 'scale invariance' of $H$, and
say exactly what kind of asymptotics $\J\nn$ has \wrt $\ul{z}$. Also, the two
fusion operators $\J\nn_\pm$ should be related by $\J\nn_+=(\J\nn_-)^{-1}$
or something of the sort.}
\begin{equation}\label{eq:Jnn}\begin{split}
\J\nn_\pm(\ul{z},\mu)
&=H_\pm(\zeta;\ul{\zeta},\mu)\,
e^{\sum_i\zeta\zeta_i\ad(\mu\ii)}\zeta^{\hbar\sum_{i<j}\Omega_{ij}^\h}
\prod_{i<j}(\zeta_i-\zeta_j)^{\hbar\Omega_{ij}^\h}\\
&=H_\pm(\zeta;\ul{\zeta},\mu)\,
e^{\sum_i z_i\ad(\mu\ii)}\prod_{i<j}(z_i-z_j)^{\hbar\Omega_{ij}^\h}
\end{split}
\end{equation}
where $H_\pm$ is given by Theorem \ref{th:Hn}.

\begin{defn}
The {\it multicomponent fusion operator} is the real analytic map 
\[J\nn_\pm:\imath\C_n\times \imath\C\to (U\g^{\otimes n}\fml)_{o}\]
given by\comment{remove dependence in $\zeta$. Also define
it on a complex domain, not the product of chambers.}
\[J\nn_\pm(\ulz,\mu)=\J_\pm\nn(\ulz,\mu)\,1^{\otimes n}=H_\pm(\zeta;\ul{\zeta},\mu)\cdot\prod_{i<j}(z_i-z_j)^{\hbar\Omega_{ij}^\h}\]
\end{defn}

It follows from Theorem \ref{th:Hn} that the fusion operator satisfies
\begin{align*}
d_{\ul{z}}J\nn_\pm &=
\left(\half{\hbar}\sum_{1\leq i<j\leq n}d\log(z_i-z_j)\Omega_{ij}+\sum_{i=1}^n dz_i\ad(\mu^{(i)})\right)J\nn_\pm \\
d_\h J\nn_\pm&=
\half{\hbar}\sum_{\alpha\in\sfPhi_+}\frac{d\alpha}{\alpha}
\left(\Delta^{(n)}(\Kalpha)\cdot J\nn_\pm-J\nn_\pm\cdot\left(\sum_{i=1}^n\Kalpha^{(i)}\right)\right)
+\sum_{i=1}^n z_i\ad(d_\h(\mu)^{(i)})J\nn_\pm
\end{align*}

\comment{
{\bf Notation} $\J\nn=\J\nn(\ul{z};\mu)$ is $\mathcal E$--valued,
$J\nn=J\nn(\ul{z};\mu)$ is $\mathcal A$--valued, $\ol{J}=\ol
{J}(\mu):=\rlim{z_1-z_2\to 0}J\twtw(\ul{z};\mu)$ is $\mathcal
A$--valued. Could also use the notation $\ol{\J}:=\rlim{z_1
-z_2\to 0}\J\twtw(\ul{z};\mu)$, but it isn't needed much.}

\subsection{Proof of Theorem \ref{th:Hn}}\label{ss:fusion begin}

(1) We carry out the proofs for $H_+$ only, those for $H_-$ being
identical, and therefore drop the subscript $+$. $H$ is required to
satisfy
\begin{equation}\label{eq:zeta ODE}
\frac{dH}{d\zeta}=
[\sum_i\zeta_i\mu\ii,H]+\hbar\frac{\Omega H-H\Omega^\h}{\zeta}
\end{equation}
Expanding $H(\zeta)=\sum_{n\geq 0}H_n(\zeta)\hbar^n$, yields the
recursive inhomogeneous ODEs
\begin{equation}\label{eq:recursive}
\frac{dH_n}{d\zeta}=[\sum_i\zeta_i\mu\ii,H_n]+\frac{\Omega H_{n-1}-H_{n-1}\Omega^\h}{\zeta}
\end{equation}
where $H_{-1}=0$ by convention, together with the requirement that
$H_n(\zeta)\to \delta_{n0}H^0$ as $\zeta\to\infty$ in $\IH$. We shall
treat the cases $n=0$, $n=1$ and $n\geq 2$ separately.

Let $Q=\IZ\Phi\subset(\h^\IR)^*$ be the root lattice and, for any
$\gamma\in Q$, let $\Ug_\gamma\subset\Ug$ be the weight space
corresponding to $\gamma$ for the adjoint action of $\h$ on $\Ug$.
Thus, $\Ug=\bigoplus_{\gamma\in Q}\Ug_\gamma$ and
\begin{equation}\label{eq:weight dec}
\Ug^{\otimes n}=
\bigoplus_{\ugamma\in Q^n}\Ug^{\otimes n}_\ugamma
\end{equation}
where, for $\ugamma=(\gamma_1,\ldots,\gamma_n)\in Q^n$, $\Ug
^{\otimes n}_\ugamma=\Ug_{\gamma_1}\otimes\cdots\otimes\Ug_
{\gamma_n}$. Notation: $\ugamma(\uzeta,\mu)=\sum_k\zeta_k
\gamma_k(\mu)$. 

\subsection{$\mathbf{n=0}$} 

The equation \eqref{eq:recursive} reduces to
\begin{equation}\label{eq:hgs}
\frac{dH_0}{d\zeta}=\sum_i\zeta_i\ad(\mu\ii)\,H_0
\end{equation}
and is therefore equivalent to $H_0(\zeta)=\exp(\zeta\sum\zeta_i\ad(\mu\ii))
\,C$, where $C\in\Ug^{\otimes n}$ is some constant element. Write $C=\sum_{\ugamma
\in Q^n}C_\ugamma$ so that $H_0=\sum_\ugamma e^{\zeta\ugamma(\uzeta,
\mu)}C_\ugamma$. Since $\ugamma(\uzeta,\mu)=\sum_k\zeta_k\gamma_k
(\mu)\in\IR$, the function $\exp
(\zeta\ugamma(\uzeta,\mu))$ has a limit as $\zeta\to\infty$ in $\IH$ if, and
only if $\ugamma(\uzeta,\mu)=0$, so that $H_0(\zeta)=\sum_{\ugamma:
\ugamma(\uzeta,\mu)=0}C_\ugamma$. Thus, $H_0$ is a constant function
of $\zeta$, and is therefore equal to its limit as $\zeta\to\infty$ in $\IH$. This
proves the existence and uniqueness of $H_0$, as well as the uniqueness
of $H_n$ for $n\geq 1$ since the homogeneous equation underlying \eqref
{eq:recursive} is \eqref{eq:hgs}.

\subsection{$\mathbf{n=1}$} 

Equation \eqref{eq:recursive} now yields
\begin{equation}\label{eq:n=1}
\frac{dH_1}{d\zeta}=
\sum_i\zeta_i\ad(\mu\ii)\,H_1+\frac{\ol{\Omega}H_0}{\zeta}
\end{equation}
where $\ol{\Omega}=\Omega-\Omega^\h$. Decomposing this equation
along \eqref{eq:weight dec}, and noting that $\ol{\Omega}=\sum_{i<j,\,
\alpha\in\Phi}\Omega_{ij}^\alpha$, where $\Omega_{ij}^\alpha=x_\alpha
\ii x_{-\alpha}\jj$, yields
\[\frac{dH_1^\ugamma}{d\zeta}=
\ugamma(\uzeta,\mu)\,H_1^\ugamma+
\sum_{\substack{\alpha\in\Phi,\\ i<j}}\delta_{\ugamma,\alpha\ii-\alpha\jj}\frac{\Omega_{ij}^\alpha H_0}{\zeta}
\]
where $H_1^\ugamma$ is the component of $H_1$ along $U\g^{\otimes n}
_\ugamma$, and $\alpha\ii=(0,\ldots,0,\alpha,0,\ldots,0)\in Q^n$, with the
$\alpha$ in the $i$th slot. This equation, together with the requirement
that $H_1^\ugamma\to 0$ as $\zeta\to\infty$, is clearly solved by setting
$H_1^\ugamma=0$ if $\gamma$ is not of the form $\alpha\ii-\alpha\jj$,
and by resorting to (2) of Proposition \ref{pr:basic ODE} otherwise since
in that case $\ugamma(\uzeta,\mu)=\alpha(\mu)(\zeta_i-\zeta_j)\neq 0$.

This yields a solution $H_1$ with values in $(U\g^{\otimes n})_{o_1}$,
which is a smooth function of $(\uzeta,\mu)\in\imath\C_n\times\imath\C$
since, again, $\ugamma(\uzeta,\mu)=\alpha(\mu)(\zeta_i-\zeta_j)\neq 0$
and moreover possesses an asymptotic expansion as $\zeta\to\infty$
with trivial constant term.
Explicitily
\[H_1(\zeta)=
\sum_{i<j,\,\alpha\in\Phi^+}\left(
\int_{-\infty}^0 e^{-t\alpha(\mu)(\zeta_i-\zeta_j)}\frac{\Omega_{ij}^{\alpha} H_0}{\zeta+t}dt
-\int_0^\infty e^{t\alpha(\mu)(\zeta_i-\zeta_j)}\frac{\Omega_{ij}^{-\alpha} H_0}{\zeta+t}dt
\right)\]
which has the asymptotic expansion
\begin{equation}\label{eq:H1 asy}
H_1\sim -\sum_{\substack{\alpha\in\Phi\\ i<j}}\frac{\Omega_{ij}^\alpha}{\alpha(\mu)(\zeta_i-\zeta_j)}\zeta^{-1}+O(\zeta^{-2})
\end{equation}

\subsection{$\mathbf{n\geq 2}$}

We now assume inductively that we have constructed $H_{n-1}$ with values in
$(U\g^{\otimes n})_{o_{n-1}}$\comment{not the same $n$}, which is a continuous
function of $(\uzeta,\mu)\in\imath\C_n\times\imath\C$, smooth away from a finite
collection $\{\Q_i\}_{i\in\I_{n-1}}$ of quadrics\footnote{such a collection is empty
if $n=2$, but will be seen to be a priori non--trivial thereafter}, and possesses an
asymptotic expansion in $\zeta\to\infty$ with leading term of the form $C_{n-1}
\zeta^{-(n-1)}$ away from $\bigcup_{i\in\I_{n-1}}\Q_i$ and $\zeta^{-1}$ otherwise, 
and set about constructing $H_n$ with the same properties.

Set $G_n=\Omega H^{n-1}-H^{n-1}\Omega^h\in(\Ug^{\otimes n})_{o_n}$.
Then decomposing equation \eqref{eq:recursive} along \eqref{eq:weight dec}
yields, for any $\ugamma\in Q^n$,
\[\frac{H_n^\ugamma}{d\zeta}=
\ugamma(\uzeta,\mu)H_n^\ugamma+\frac{G_n^\ugamma}{\zeta}\]

Since $G_n\sim\zeta^{-(n-1)}$, part (3) of Proposition \ref{pr:basic ODE}
shows the existence of a unique solution $H_n^\ugamma$ with values in
$V=(\Ug^{\otimes n})_{o_n}$. Explicitly
\[H^n(\zeta)=
-i\sum_{\ugamma\in Q^n}\int_0^\infty
e^{-it\sum_k\zeta_k\gamma_k(\mu)}\frac{G^n_\ugamma(\zeta+it)}{\zeta+it}dt\]
This solution is continuous on $\imath\C_n\times\imath\C$, and smooth on the
complement of the quadrics $\Q_i$, $i\in\I_{n-1}$ and
\[\mathcal Q_\ugamma=
\{(\uzeta,\mu)\in\imath\C_n\times\imath\C|\,\ugamma(\uzeta,\mu)=0\}\]
where $\ugamma$ ranges over the (finitely many) elements of $Q^n$ such
that $G_n^\gamma\neq 0$. Moreover, it has an asymptotic expansion in
$\zeta$ starting at $\zeta^{-n}$ on the complement of these quadrics, and
at $\zeta^{-1}$ elsewhere.

\subsection{}

Our next goal is to show that $H=\sum_{n\geq 0}\hbar^n H_n$ is smooth
on $\imath\C_n\times\imath\C$ and satisfies the PDE \eqref{eq:induced}.
We shall do so for the truncation $\ol{H}$ of $H$ mod $\hbar^p$\comment
{what exactly are we modding by?} and then let $p$ tend to $\infty$. By
construction, $\ol{H}$ is continuous on $\imath\C_n\times\imath\C$, and
smooth on the complement $\U$ of finitely many (real) quadrics. Thus, if
$\ul{D}$ is the flat connection
\begin{multline*}
\ul{D}=
\Omit{\ul{d}-\zeta\ul{d}\sum_{i=1}^n\zeta_i\ad(\mu\ii)
-\ell\left(\hbar\negthickspace\negthickspace\sum_{1\leq i<j\leq n}d\log(\zeta_i-\zeta_j)\,\Omega_{ij}
+\frac{\hbar}{2}\sum_{\alpha\in\Phi_+}\frac{d\alpha}{\alpha}\,\Delta^{(n)}(\Kalpha)\right)\\
+r\left(\hbar\sum_{i<j}d\log(\zeta_i-\zeta_j)\,\Omega_{ij}^\h
+\frac{\hbar}{2}\sum_{\alpha\in\Phi_+}\frac{d\alpha}{\alpha}(\sum_{i=1}^n\,\Kalpha\ii)\right)}
\ul{d}-\zeta\ul{d}\sum_{i=1}^n\zeta_i\ad(\mu\ii)
-\hbar\negthickspace\negthickspace\sum_{1\leq i<j\leq n}d\log(\zeta_i-\zeta_j)\,
\left(\ell(\Omega_{ij})-r(\Omega_{ij}^\h)\right)\\
-\frac{\hbar}{2}\sum_{\alpha\in\Phi_+}\frac{d\alpha}{\alpha}\,
\left(\ell(\Delta^{(n)}(\Kalpha))-r(\sum_{i=1}^n\,\Kalpha\ii)\right)
\end{multline*}
where $\ell,r$ denote left and right multiplication respectively, then $G=\ul{D}H$
is well--defined on $\IH\times\U$. To show that $D_\zeta G=0$, notice first that
\[D_\zeta G=D_\zeta\ul{D}H=\ul{D}D_\zeta H=0\]
By uniqueness, it therefore suffices to show that, as a function of $\zeta$,
$G$ is holomorphic on $\IH$ and tends to 0 as $\zeta\to\infty$ in $\IH$.
Where $D_\zeta=d_\zeta-\ad(\sum_i\zeta_i\mu\ii)-\ell(\hbar\Omega)+r(\hbar\Omega^\h)$
and, by integrability $[D_\zeta,\ul{D}]=0$.

Since $H=1+H^1\zeta^{-1}+O(\zeta^{-2})$, it is clear that $G$ has a well--defined
limit as $\zeta\to\infty$, given by
\[-[\ul{d}\sum_{i=1}^n\zeta_i\ad(\mu\ii),H^1]
-\hbar\negthickspace\negthickspace\sum_{1\leq i<j\leq n}d\log(\zeta_i-\zeta_j)\,
(\Omega_{ij}-\Omega_{ij}^\h)
-\frac{\hbar}{2}\sum_{\alpha\in\Phi_+}\frac{d\alpha}{\alpha}\,\left(
\Delta^{(n)}(\Kalpha)-\sum_{i=1}^n\,\Kalpha\ii)\right)\]

Since
\[(\Omega_{ij}-\Omega_{ij}^\h)=\sum_\alpha\Omega_{ij}^\alpha
\aand \Delta^{(n)}(\Kalpha)-\sum_{i=1}^n\,\Kalpha\ii=
2\sum_{i<j}(\Omega_{ij}^\alpha+\Omega_{ij}^{-\alpha})\]
the last two summands add up to
\[-\hbar\sum_{\substack{i<j\\\alpha\in\Phi}}
\left(d\log(\zeta_i-\zeta_j)+\frac{d\alpha}{\alpha}\right)
\Omega_{ij}^\alpha\]
On the other hand, by \eqref{eq:H1 asy}, 
\[\begin{split}
[\ul{d}\sum_{k=1}^n\zeta_k\ad(\mu\kk),H^1]
&=
-\hbar\sum_{\substack{\alpha\in\Phi\\ i<j}}
\sum_{k=1}^n[\ul{d}(\zeta_k\ad(\mu\kk)),\frac{\Omega_{ij}^\alpha}{\alpha(\mu)(\zeta_i-\zeta_j)}]\\
&=
-\hbar\sum_{\substack{\alpha\in\Phi\\ i<j}}
\left(\frac{d(\zeta_i-\zeta_j)}{\zeta_i-\zeta_j}+\frac{d\alpha}{\alpha}\right)\Omega_{ij}^\alpha
\end{split}\]
and we are done.\comment{explain why $H^1$ does not contain corrections
coming from $\hbar^q$, $q\geq 2$.}

Thus, $\ol{H}$ is a horizontal section of $\ul{D}$ on each of the
connected components of $\U$, and therefore extends to a smooth
horizontal section on the closure of each connected component in
$\imath\C_n\times\imath\C$. Since however $\ol{H}$ is continuous
on $\imath\C_n\times\imath\C$, this extension coincides with $\ol{H}$
thus showing that $\ol{H}$ is a smooth, horizontal section of $\ul{D}$
on the whole of $\imath\C_n\times\imath\C$.

\subsection{} 
\label{ss:fusion end}
\comment{(to be polished/corrected)}
Let us prove that $\lim_{\zeta_1\to+\imath\infty}H^{(n)}=1\otimes H^{(n-1)}$.
By (4) of Proposition \ref{pr:basic ODE}, $H\nn$ possesses an asymptotic
expansion $H\nn=G_0+G_1\zeta_1^{-1}+\cdots$. Plugging this into \eqref
{eq:zeta ODE}, and taking the coefficients of $\zeta_1$ and $\zeta_1^0$ respectively gives $
[\mu^{(1)},G_0]=0$ and
\[\frac{dG_0}{d\zeta}=\frac{[M_{\geq 2},G_0]+[\mu^{(1)},G_1]}{\zeta^2}+\frac{K G_0-G_0[K]}{\zeta}\]
Projecting onto $Z(\mu^{(1)})$ then yields
\begin{equation}\label{eq:intermediate}
\frac{dG_0}{d\zeta}=\frac{[M_{\geq 2},G_0]}{\zeta^2}+\frac{K_0 G_0-G_0[K]}{\zeta}
\end{equation}
where $K_0=P_1 K$. Write
\[K=\sum_{j\geq 2}\hbar\Omega_{1j}+\hbar\Omega_{\geq 2}=
\id\otimes\Delta^{(n-1)}(\hbar\Omega)+\hbar\Omega_{\geq 2}\]
Then $K_0=\id\otimes\Delta^{(n-1)}(\hbar\Omega^\h)+\hbar\Omega_{\geq 2}$. Note
that since the coefficients of \eqref{eq:intermediate} commute with the action of $\h
\otimes\h$ on $U\g\otimes U\g^{\otimes (n-1)}$, so does $G_0$ by uniqueness. It
follows that $G_0$ commutes with $\id\otimes\Delta^{(n-1)}(\hbar\Omega^\h)$ so
the above may be rewritten as
\[\frac{dG_0}{d\zeta}=\frac{[M_{\geq 2},G_0]}{\zeta^2}+\frac{\hbar\Omega_{\geq 2} G_0-G_0[\hbar\Omega_{\geq 2}]}{\zeta}\]
which is the differential equation for the sought for gauge transformation.

The proof that $\lim_{\zeta_n\to-\imath\infty}H^{(n)}=H^{(n-1)}\otimes 1$
is identical.

\section{The differential twist}\label{se:diffl from fusion} 

In this section, we construct a differential twist $J_\pm(\mu)$ for
$\g$ as the regularised limit of the fusion operator $J_\pm(z,\mu)$
when $z\to 0$, and prove that it kills the KZ associator.

\subsection{Fundamental solution of the KZ equations near $z=0$}

\begin{prop}\label{pr:Fuchs 0}\hfill
\begin{enumerate}
\item For any $\mu\in\h$, there is a unique holomorphic function
$H_0:\IC\to\A$ such that $H_0(0,\mu)\equiv 1$ and, for any determination
of $\log(z)$, the $\E$--valued function
$\Upsilon_0(z,\mu)=e^{z\ad\muone}\cdot H_0(z,\mu)\cdot z^{\hbar\Omega}$
satisfies
\[\left(d_z-\left(\hbar\frac{\Omega}{z}+\ad\muone\right) dz \right)\Upsilon_0=\Upsilon_0\,d_z\]
%
\item $H_0$ and $\Upsilon_0$ are holomorphic functions of $\mu$,
and $\Upsilon_0$ satisfies
\[\left(d_\h-\half{\hbar}\sum_{\alpha\in\Phi_+}\frac{d\alpha}{\alpha}
\Delta(\Kalpha)-z\ad(d\muone)\right)\Upsilon_0=
\Upsilon_0\left(d_\h-\half{\hbar}\sum_{\alpha\in\Phi_+}\frac{d\alpha}{\alpha}
\Delta(\Kalpha)\right)\]
\item $H_0,\Upsilon_0$ commute with the diagonal action of $\g$
if $\mu=0$.
\end{enumerate}
\end{prop}
\begin{pf}
(1) $H=H_0$ is required to satisfy
\[\begin{split}
\frac{dH}{dz}
&=\frac{\hbar}{z}\left(e^{-z\ad\muone}(\Omega)H-H\Omega\right)\\
&=\hbar\left(\frac{[\Omega,H]}{z}+\frac{e^{-z\ad\muone}-1}{z}(\Omega)H\right)
\end{split}\]
Writing $H=\sum_{n\geq 0}\hbar^n H_n$ yields the recursive system of
ODEs
\[\frac{dH_n}{dz}
=\left(\frac{[\Omega,H_{n-1}]}{z}+\frac{e^{-z\ad\muone}-1}{z}(\Omega)H_{n-1}\right)\]
where $H_{-1}=0$, together with the initial value condition $H_n(0)=\delta_{n0}$,
which clearly possesses at most one solution. For $n=0$, the solution is given
by $H_0\equiv 1$ and, for $n\geq 1$, by
\[H_n(z)=\int_0^z \frac{[\Omega,H_{n-1}(t)]}{t}+\frac{e^{-t\ad\muone}-1}{t}(\Omega)H_{n-1}(t)\,dt\]
provided the integral converges at $t=0$. For $n=1$, this is the case since
$H_0=1$ commutes with $\Omega$, so the integrand is an $O(1)$ and, for
$n\geq 2$, this follows since $H_{n-1}(t)=O(t)$.

(2) The recursive construction of $H$ shows that it, and therefore $\Upsilon_0$,
are holomorphic functions of $\mu\in\h$. The fact that $\Upsilon_0$ satisfies the
stated PDE follows by integrability.

(3) For $\mu=0$, $H_0=1$ and $\Upsilon_0=z^{\hbar\Omega}$ clearly commute
with $\g$.
\end{pf}

\subsection{Differential twist}\label{ss:diffl twist}

Let $J_\pm(z,\mu)$ be the fusion operator defined in \S \ref{ss:multicomponent}.

\begin{defn}
The differential twist of $\g$ is the smooth map $J_\pm:i\C\to\Ug^{\otimes 2}\fml_o$
given by
\[J_\pm(\mu)=\Upsilon_0(z,\mu)^{-1}\cdot J^{(2)}_\pm(z,\mu)\]
\end{defn}

It follows from \ref{ss:multicomponent} that $J_\pm(\mu)$ is independent of $z$,
and satisfies\comment{Give the 1--jet of the differential twist and 'prove' that it is
killed by $\veps\otimes 1$ and $1\otimes\veps$.}
\[d_\h J_\pm =
\half{\hbar}\sum_{\alpha\in\sfPhi_+}\frac{d\alpha}{\alpha}
\left(\Delta(\Kalpha)\cdot J_\pm-J_\pm\cdot(\Kalpha^{(1)}+\Kalpha^{(2)})\right)\]
The following is the main result of this section. Its proof will be given in \ref{ss:J Phi}
after some preparatory work.

\begin{thm}\label{th:J kills Phi}
The following holds
\[\Phi\KKZ\cdot\Delta\otimes\id(J_\pm)\cdot J_\pm\otimes 1=\id\otimes\Delta(J_\pm)\cdot 1\otimes J_\pm\]
\end{thm}

\comment{The only property of the fusion operator which is used to show that it kills $\Phi\KKZ$ is the inductive property
$\lim_{z_1\to\infty}\J\nn=1\otimes\J^{(n-1)}$ and $\lim_{z_n\to-\infty}
\J\nn=\J^{(n-1)}\otimes 1$. Highlight this?}

\subsection{}\label{ss:start bottom row}

The goal of \S \ref{ss:start bottom row}--\S \ref{ss:end bottom row}
is to show that the constants relating the fundamental solutions of
the dynamical KZ equations corresponding to Drinfeld's asymptotic
zones are the same as those relating their non--dynamical counterparts.
We begin by constructing the relevant fundamental solutions.

Fix $n\geq 2$, and let $\B_n$ be the set of complete bracketings
on the monomial $x_1\cdots x_n$. By convention, such bracketings
always contain the parentheses $\bigp$, and are easily seen to consist
of $n-1$ compatible pairs of parentheses. Let $\Cno=\{\ulz\in\IC^n|\,
\sum_i z_i=0\}$. For any $b\in\B_n$, the coordinates $z_P=z_i-z_j$,
where $P=\cdots(x_i\cdots x_j)\cdots\in b$, are a basis of $(\IC_0^n)^*$.
Let $\U_b\cong\Cno$ be the affine space with coordinates $\{u_P\}
_{P\in b}$. Consider the regular map $\rho_b:\U_b\to\Cno$ given by
\[z_P=\prod_{b\ni P'\supseteq P} u_{P'}\]
$\rho_b$ is a birational map, with inverse given by $u_P=z_P/z_Q$,
where $Q\in b$ is the smallest pair of parentheses strictly containing
$P$, and $z_Q=1$ if $P=\bigp$. It induces an isomorphism
\[\U_b\setminus\bigcup_{P\in b}\{u_P=0\}\isom
\Cno\setminus\bigcup_{P\in b}\{z_P=0\}\]
The pull--back $\rho_b^*\nabla\DKKZ$ is of the form
\[\rho_b^*\nabla\DKKZ=d-\hbar\sum_{P\in b}\dlog{u_P}\Omega_P-\R_b\]
where, for $P=\cdots(x_i\cdots x_j)\cdots$, $\Omega_P=\sum_{i\leq
k<l\leq j}\Omega_{kl}$, and $\R_b$ is an $\EE$--valued one--form
which is regular in the neighborhood of $0\in\U_b$. It follows from
this, or by direct inspection, that $[\Omega_P,\Omega_Q]=0$ for
any $P,Q\in b$.

Fix now $b\in\B_n$ and, for any $P\in b$, set $\Omega_{\ul{P}}
=\sum_{Q}\Omega_Q$, where $Q$ ranges over the maximal
elements of $b$ properly contained in $P$.

\begin{prop}\hfill
\begin{enumerate}
\item For any simply--connected neighborhood $\V_b$ of $0\in\U_b$,
there exists a unique holomorphic function $H_b:\V_b\to\EE$ such that
$H_b(0)=1$, and the function
\[\Psi_b=
H_b\prod_{P\in b}u_P^{\Omega_P}=
H_b\prod_{P\in b}z_P^{\Omega_P-\Omega_{\ul{P}}}\]
is a fundamental solution of $\rho^*_b\nabla\KKZ$.
\item $H_b$ is a holomorphic function of $\mu\in\h$, and satisfies
\[\left(d_\h-\half{\hbar}\sum_\alpha\dlog{\alpha}\Delta\nn(\Kalpha)-
\sum_i z_i\ad(d_\h\mu\ii)\right)H=
H\left(\half{\hbar}\sum_\alpha\frac{d\alpha}{\alpha}\Delta\nn(\Kalpha)\right)\]
\end{enumerate}
\end{prop}
\begin{pf}
(1) is standard, and proved as in Proposition \ref{pr:Fuchs 0}.\comment
{This is the general $n$ statement of the proposition proved in detail at
the start of the section for $n=2$. This is a little gauche.}

(2) Write the equation as $DH=0$, where $D$ is the flat connection
\comment{potential clash of notation, not clear whether the $\ad$ is in $\A$ or $\E$.}
\[D=d_\h-\half{\hbar}\sum_\alpha\dlog{\alpha}[\Left{\Delta\nn(\Kalpha)},-]-\sum_i z_i \ad(d\mu\ii)\]
Let $G=DH$. Then $\nabla\KKZ G=D\nabla\KKZ H=0$\comment{not quite, need to add
terms on the right for $\nabla\KKZ$ as well}. Moreover, $G$ is holomorphic near $u=0$
and $G(0)=0$\comment{check this}, whence by uniqueness $G\equiv 0$ as claimed.
\end{pf}

\subsection{}

Let now $b,b'\in\B_n$ be two bracketing on $n$ letters.

\begin{lemma}
The following holds
\[\Upsilon_b=\Upsilon_{b'}\,\Left{\Phi\KKZ}\]
where $X^\ell$ is the operator of left multiplication by $X$.
\end{lemma}
\begin{pf}
Set $\Phi=(\Upsilon_{b'})^{-1}\Upsilon_b$. $\Phi$ is a holomorphic function
of $\mu\in\h$ since both $\Upsilon_b$ and $\Upsilon_{b'}$ are, and satisfies
$\Phi(0)=\Left{\Phi_{b'b}}$ and
\[d\Phi=\frac{\hbar}{2}\sum_\alpha\dlog{\alpha}[\Left{\Delta(\Kalpha)},\Phi]\]
Fix $\mu\in\hreg$, and identify the line $\IC\mu\subset\hreg$ with $\IC$ via
$\IC\ni t\to t\mu$. The restriction of $\Phi$ to $\IC\mu$ satisfies
\[\frac{d\Phi}{dt}=\half{\hbar}[\Left{\Delta(\K)},\Phi]\]
where $\K=\sum_\alpha\Kalpha$. Thus $\ds{\Phi(\mu)=\Ad(Exp(\half{\hbar}
\Left{\Delta\nn(\K)}))\Phi(0)}=\Left{\Phi_{b'b}}$\comment{{\bf not correct} the
ODE satisfies by $\Psi$ is $\frac{d\Phi}{dt}=\frac{\hbar}{2t}[\Left{\Delta(\K)},\Phi]$}
since $\Phi_{b'b}$ is invariant.
as claimed. Since $\hreg\subset\h$ is dense, and $\Phi$ continuous on $\h$
it follows that $\Phi\equiv\Left{\Phi_{b'b}}$ as claimed.
\end{pf}

\subsection{}\label{ss:end bottom row}
For any \comment{$\A$ or $\E$--valued} solution of DKZ$_n$ $\Psi$,
and bracketing $b\in\B_n$, set $\rlim{b}\Psi=\Upsilon_b^{-1}\cdot\Psi$.

\begin{corollary}\label{co:bb'}
Let $\Psi$ be an $\E$--valued solution of DKZ$_n$, and $b,b'$ two
bracketings on $n$ letters. Then
\[\rlim{b'}\Psi=\Left{\Phi_{b'b}}\rlim{b}\Psi\]
\end{corollary}
\begin{pf}
Let $C_b=\rlim{b}\Psi$. Then, $\Psi=\Upsilon_bC_b=\Upsilon_{b'}\Left{\Phi_{b'b}}C_b$
whence $\lim_{b'}\Psi=\Left{\Phi_{b'b}}C_b$ as claimed.
\end{pf}

\subsection{The two normalised limits}

\begin{prop}\label{pr:Upsilon}\hfill
\begin{enumerate}
\item The function $\Upsilon\aBC=J^{(3)}\cdot 1\otimes(J^{(2)})^{-1}$
is regular at $z_2=z_3$, and
\[\Upsilon\aBC(z_1,z_2,z_2;\mu)=
\id\otimes\Delta\left(J\twtw(z_1,z_2;\mu)\right)\]
\item The function $\Upsilon\ABc=J^{(3)}\cdot (J^{(2)})^{-1}\otimes 1$
is regular at $z_1=z_2$, and
\[\Upsilon\ABc(z_1,z_1,z_3;\mu)=
\Delta\otimes\id\left(J\twtw(z_1,z_3;\mu)\right)\]
\end{enumerate}
\end{prop}
\begin{pf}
(1) As a function of $z_1$, $\Upsilon=\Upsilon\aBC$ satisfies
\[\frac{d\Upsilon}{dz_1}=
\left(\ad\mu\oo+
\hbar\frac{\Omega_{12}}{z_1-z_2}+\hbar\frac{\Omega_{23}}{z_1-z_3}\right)
\Upsilon\]
Moreover, by \eqref{eq:Jnn}
\[\begin{split}
\Upsilon
&=
H\thth(-)\cdot 
(z_2-z_1)^{\hbar\Omega_{12}^\h}(z_3-z_1)^{\hbar\Omega_{13}^\h}\cdot
1\otimes H\twtw(-)^{-1}\\
&=
H\thth(-)\cdot
(1-z_2/z_1)^{\hbar\Omega_{12}^\h}(1-z_3/z_1)^{\hbar\Omega_{13}^\h}\cdot
1\otimes H\twtw(-)^{-1}\cdot
(-z_1)^{\hbar\id\otimes\Delta(\Omega^\h)}\\
&=
\Xi(-)\cdot
(z_2-z_1)^{\hbar\Omega_{12}^\h}(z_3-z_1)^{\hbar\Omega_{13}^\h}
\end{split}\]
where the second equality follows from the fact that $H\twtw$ is of weight
zero, and the third defines the function $\Xi$. As a function of $z_1$, $\Xi$
is holomorphic\comment{where? Also, clarify how to get rid of the $\zeta$
dependence of $H\nn$} and, by \eqref{it:rec Hn} of Theorem \ref{th:Hn},
tends to 1 as $z_1\to\infty$\comment{which infinity?}.

$\Xi$ satisfies
\[\frac{d\Xi}{dz_1}=
\left(\ad\mu\oo+
\hbar\frac{\ell(\Omega_{12})-\err(\Omega_{12}^\h)}{z_1-z_2}+
\hbar\frac{\ell(\Omega_{13})-\err(\Omega_{13}^\h)}{z_1-z_3}\right)\Xi\]
\Omit{
\hbar\frac{\Omega_{12}}{z_1-z_2}+\hbar\frac{\Omega_{23}}{z_1-z_3}\right)\Xi
-\Xi\frac{\Omega_{12}^\h+\Omega_{23}^\h}{z_1}\]}

A proof similar to that of (1) of Theorem \ref{th:Hn}, and based on Proposition
\ref{pr:basic ODE} shows that there exists
a unique holomorphic function $\Xi$ \comment{$:?\to?$} satifsying the above equation
and tending to 1 as $z_1\to\infty$\comment{say which infinity}. That function is moreover holomorphic
in $z_2,z_3\in\C^2$\comment{not really, they have to be restrained a little
since they are poles of the equation} and regular at $z_2=z_3$, where it
coincides, by uniqueness, with $\id\otimes\Delta(H\twtw(z_1,z_2))$.
\end{pf}

\subsection{}

\renewcommand {\JJ}{J}

\begin{corollary}\label{co:rlimsJ}
The following holds
\[\rlim{((\cdot\cdot)\cdot)}J^{(3)}=\Delta\otimes\id(\JJ)\cdot\JJ\otimes 1
\aand
\rlim{(\cdot(\cdot\cdot))}J^{(3)}=\id\otimes\Delta(\JJ)\cdot 1\otimes \JJ\]
\end{corollary}
\begin{pf}
By definition,\comment{assuming it is understood that $J$ only depends
upon the differences $z_i-z_j$.}
\[\rlim{((\cdot\cdot)\cdot)}J^{(3)}=
\lim_{\substack{z_1-z_3\to 0\\[.5ex]\frac{z_1-z_2}{z_1-z_3}\to 0}}
(z_1-z_3)^{-\hbar\Delta\otimes\id(\Omega)}(z_1-z_2)^{-\hbar\Omega_{12}}J^{(3)}\]
Write $J^{(3)}=\Upsilon\ABc\cdot\J^{(2)}\otimes 1$, where $\Upsilon\ABc$
is defined in Proposition \ref{pr:Upsilon}.
Since
\[\Upsilon\ABc(z_1,z_2,z_3)=\Delta\otimes\id\left(J\twtw(z_1,z_2)\right)+(z_1-z_2)\R\]
where $\R$ is regular at $z_1=z_2$, and $\Omega$ commutes with
$\Delta(U\g)$, we have
\[(z_1-z_2)^{-\hbar\Omega_{12}}\Upsilon\ABc(z_1,z_2,z_3)=
\Delta\otimes\id\left(J\twtw(z_1,z_2)\right)(z_1-z_2)^{-\hbar\Omega_{12}}
+(z_1-z_2)^{-\hbar\Omega_{12}}\R\]
Since the second summand tends to zero as $z_1-z_2\to 0$, it follows
that
\[\begin{split}
\rlim{((\cdot\cdot)\cdot)}J^{(3)}
&=
\lim_{\substack{z_1-z_3\to 0\\[.5ex]\frac{z_1-z_2}{z_1-z_3}\to 0}}
\Delta\otimes\id\left((z_1-z_3)^{-\hbar\Omega}J\twtw(z_1,z_3)\right)
\cdot(z_1-z_2)^{-\hbar\Omega_{12}}J^{(2)}(z_2,z_3)\otimes 1\\
&=
\Delta\otimes\id(\JJ^{(2)})\cdot\JJ^{(2)}\otimes 1
\end{split}\]
The second one follows in a similar way.
\end{pf}

\subsection{Proof of Theorem \ref{th:J kills Phi}}\label{ss:J Phi}

By Corollary \ref{co:bb'},
$\rlim{(\cdot(\cdot\cdot))}J^{(3)}=\Phi\KKZ\cdot\rlim{((\cdot\cdot)\cdot)}J^{(3)}$.
The result now follows from Corollary \ref{co:rlimsJ}.

\section{The centraliser property}\label{se:centraliser}%

In this section, we prove that the differential twist for $\g$ obtained
from the fusion operator possesses the centraliser property. This
follows from a detailed analysis of the asymptotics of solutions of
the joint Casimir--KZ equations in $n=2$ points, the regime where
$z=z_1-z_2\to 0$, and a fixed root coordinate $\alpha_i$ tends to
infinity.


\subsection{}

Consider the joint KZ--Casimir connection when $n=2$. Since
we will only consider solutions with values in $(\Ugo{2})^\h$,
the connection reads
\begin{equation}\label{eq:DCKZ 2}
\nabla=d
-\hbar\Omega\frac{dz}{z}
-\frac{\hbar}{2}\sum_{\alpha\in\Phi_+}\frac{d\alpha}{\alpha}
\Delta(\Kalpha)-d(z\ad\muone)
\end{equation}
where $z=z_1-z_2\in\IC^\times$.
Fix $i\in\bfI$. We construct below a horizontal section of $\nabla$
with prescribed asymptotics as $\alpha_i\to\infty$, $\Im\alpha_i
\gtrless 0$.

Retain the notation of \ref{ss:preamble} and \ref{se:Fuchs infty},
and trivialise the fibration $\IC\cow{i}\to\h\to\olh$ by $\IC\times\olh
\ni(w,\olmu)\to w\cow{i}+\imath(\olmu)$, where $\imath:\olh\to\h$
is given by \eqref{eq:emb i}. For fixed $\olmu\in\olh$ and $z\in\IC$,
the restriction of the connection $\nabla$ to $\pi^{-1}(\olmu)\times
\{z\}$ is equal to
\begin{equation}\label{eq:nabla zi infty}
\nabla^i=
d_w-
\left(\half{\hbar}\sum_{\alpha\in\Phi_+\setminus\olPhi}\frac{\Delta(\Kalpha)}{w-w_\alpha}
+z\ad\cowo{i}\right)dw
\end{equation}
where $w_\alpha=-\alpha(\imath(\olmu))/\alpha(\cow{i})$.
Let $R_\olmu$ be given by \eqref{eq:Rmu}, and set
\begin{equation}\label{eq:Lambda}
\Lambda=\frac{\cow{i}\otimes\cow{i}}{\|\cow{i}\|^2}
\end{equation}

\subsection{Fuchs--Stokes solution}
\label{ss:Stokes}

Let $\D^\pm_\infty\subset\IC\times\olh$ be the domain given
by \eqref{eq:Dinftypm}.

\begin{prop}\label{pr:Stokes infty}\hfill
\begin{enumerate}
\item For any $\olmu\in\olh$ and $z\in\IR^\times$, there is
a unique holomorphic function
\[H_\infty^\pm:\{w\in\IC|\,\Im w\gtrless 0,\,|w|>R_{\olmu}\}\to\A\]
such that $H_\infty^\pm(w,\ol{\mu},z)\to 1$ as $\alpha_i\to
\infty$ with $0<<|\arg w|<<\pi$ and, for any determination
of $\log(\alpha_i)$, the $\E$--valued function
\[\Psi_\pm^\infty=
H_\infty^\pm(\alpha_i,\ol{\mu},z)\cdot
e^{z\alpha_i\ad(\cowo{i})}\cdot
\alpha_i^{\half{\hbar}(\onetwo{(\K-\ol{\K})})}\]
satisfies $\nabla^i\,\Psi^\pm_\infty=\Psi^\pm_\infty\,d_w$.
\Omit{
\[\left(d_w-\left(z\ad\cowo{i}+\half{\hbar}\sum_{\alpha\in\Phi_+\setminus\olPhi}
\frac{\K_\alpha}{\alpha_i-w_\alpha}\right)d\alpha_i\right)
\Psi^\pm_\infty=\Psi^\pm_\infty\,d\]}
\item The function $H_\infty^\pm(\alpha_i,\ol{\mu},z)$ is smooth
on $\D_\infty^\pm\times\IR^\times$, and $\Psi_\pm^\infty$ satisfies
\begin{multline*}
\left(d-\hbar\Omega\frac{dz}{z}-\half{\hbar}\sum_{\alpha\in\Phi_+}\dloga\Delta(\Kalpha)-d(z\ad\mu^{(1)})\right)
\Psi_\pm^\infty\\
=\Psi_\pm^\infty
\left(d-\hbar(\ol{\Omega}+\Lambda)\frac{dz}{z}-\half{\hbar}\sum_{\alpha\in\olPhi_+}\frac{d\alpha}{\alpha}\Delta(\Kalpha)-d(z\ad\imath(\olmu)^{(1)})\right)
\end{multline*}
\item $\Psi_\pm^\infty(\alpha_i,0,z)$ is invariant under the diagonal
action of $\oll$.
\item The function
\[\Psi_\pm^\infty(\alpha_i,\ol{\mu},z)\cdot e^{-z\alpha_i\ad(\cowo{i})}=
H^\infty_\pm(\alpha_i,\ol{\mu},z)\cdot\alpha_i^{\half{\hbar}(\onetwo{(\K-\ol{\K})})}\]
admits a limit for $z\to\infty$, which is equal to $\Upsilon_\infty^{\otimes 2}$,
where the latter is the solution given by Proposition \ref{pr:Fuchs infty}.
\end{enumerate}
\end{prop}
\subsection{}

\begin{pf}
(1) For fixed $z\in\IR^\times$, and $\ol{\mu}\in\ol{\h}$, $H=H^\infty_\pm$
is required to satisfy the ODE
\begin{equation}\label{eq:Stokes ODE}
\begin{split}
\frac{dH}{dw}
&=
z[\cowo{i},H]
+\half{\hbar}\left(
\sum_{\alpha\in\Phi_+\setminus\ol{\Phi}}\frac{\Delta(\Kalpha)H}{w-w_\alpha}
-\frac{H(\onetwo{(\K-\ol{\K})})}{w}\right)\\
&=
z[\cowo{i},H]
+\frac{\hbar}{2w}\left(
\Delta(\K-\ol{\K})H-H(\onetwo{(\K-\ol{\K})})
+
\sum_{\alpha\in\Phi_+\setminus\ol{\Phi}}\frac{w_\alpha\Delta(\Kalpha)}{w-w_\alpha}H\right)
\end{split}
\end{equation}
where we used \eqref{eq:res at infty}. Writing $H=\sum_{n\geq 0}\hbar^n H_n$,
this is equivalent to the recursive system of ODEs
\begin{multline*}
\frac{dH_n}{dw}
=
z[\cowo{i},H_n]\\
\phantom{\gamma(\cow{i})\gamma(\cow{i})}
+\frac{\hbar}{2w}\left(
\Delta(\K-\ol{\K})H_{n-1}-H_{n-1}(\onetwo{(\K-\ol{\K})})
+
\sum_{\alpha\in\Phi_+\setminus\ol{\Phi}}\frac{w_\alpha\Delta(\Kalpha)}{w-w_\alpha}H_{n-1}\right)
\end{multline*}
where $H_{-1}=0$, together with the condition that $H_n\to\delta_{n0}$
as $w\to\infty$ in $\IH_\pm$ with $0<<|\arg(w)|<<\pi$.

Let $\sfQ\subset\h^*$ be the root lattice, and $U\g^{\otimes 2}=\bigoplus
_{\gamma\in \sfQ}U\g^{\otimes 2}_\gamma$ the weight decomposition
\wrt the adjoint action of $\h$ acting on the first tensor copy. In terms
of the components $H_n^{\gamma}$ of $H_n$, $\gamma\in \sfQ$, the
above equation reads
\begin{multline}\label{eq:Stokes rec}
\frac{dH^\gamma_n}{dw}
=
z\gamma(\cow{i})H_n^\gamma\\
\phantom{\gamma(\cow{i})\gamma(\cow{i})}
+\frac{\hbar}{2w}\left[
\Delta(\K-\ol{\K})H_{n-1}-H_{n-1}(\onetwo{(\K-\ol{\K})})
+
\sum_{\alpha\in\Phi_+\setminus\ol{\Phi}}\frac{w_\alpha\Delta(\Kalpha)}{w-w_\alpha}H_{n-1}\right]^\gamma
\end{multline}
We shall treat the cases $n=0$, $n=1$ and $n\geq 2$ separately.

{$\mathbf{n=0}$}. In this case, $H_0\equiv 1$ is a solution of \eqref
{eq:Stokes rec}, which is unique by Proposition \ref{pr:basic ODE}.

{$\mathbf{n=1}$}. Given that $H_0=1$, the equation reads
\[\frac{dH^\gamma_1}{dw}
=
z\gamma(\cow{i})H_n^\gamma
+\frac{\hbar}{2w}\left[
\Delta(\K-\ol{\K})-(\onetwo{(\K-\ol{\K})})
+
\sum_{\alpha\in\Phi_+\setminus\ol{\Phi}}\frac{w_\alpha\Delta(\Kalpha)}{w-w_\alpha}\right]^\gamma\]
By Proposition \ref{pr:basic ODE}, this has a unique solution with the
required limiting behaviour unless $\left[\Delta(\K-\ol{\K})-(\onetwo
{(\K-\ol{\K})})\right]^\gamma\neq 0$ and $z\gamma(\cow{i})=0$.
Since $z\neq 0$, this is ruled out by the fact that
\[\Delta(\K-\ol{\K})=
\onetwo{(\K-\ol{\K})}+\sum_{\alpha\in\Phi\setminus\ol{\Phi}}x_\alpha\otimes x_{-\alpha}\]
and that $\alpha(\cow{i})\neq 0$ for any $\alpha\in\Phi\setminus\ol{\Phi}$.

{$\mathbf{n\geq 2}$}. The existence and uniqueness of $H_n^
\gamma$ follows from Proposition \ref{pr:basic ODE} since, by
induction, the inhomogeneous term of \eqref{eq:Stokes rec} is
an $O(w^{-2})$.

\Omit{
(2) It follows by Proposition \ref{pr:basic ODE} that $H$ is a smooth
function of $z\in\IR^\times$ and $\ol{\mu}\in\ol{\h}$. The fact that
$\Psi_\pm^\infty$ satisfies the claimed PDE follows by integrability.}

(2) It follows by Proposition \ref{pr:basic ODE} that $H$ is a smooth
function of $z\in\IR^\times$ and $\ol{\mu}\in\ol{\h}$. The fact that
$\Psi_\pm^\infty$ satisfies the claimed PDE follows by integrability.
\Omit{
Specifically, writing any root $\alpha\in\h^*\cong(\IC\oplus\olh)^*$
as $\alpha=\alpha(\cow{i})\alpha_i+\alpha\circ\imath\pi$, we see
that $\nabla$ decomposes as $\nabla^i+\ol{\nabla}$, where
$\nabla^i$ is given by \eqref{eq:nabla zi infty}, and 
\begin{multline*}
\ol{\nabla}=
\left(d_z-\left(\hbar\frac{\Omega}{z}+\ad(w\cow{i}+\imath(\olmu))^{(1)}\right)dz\right)\\
+
\left(d_\olh-\half{\hbar}\sum_{\alpha\in\Phi_+}
\frac{(\imath\pi)^*d\alpha}{\alpha(\cow{i})w+(\imath\pi)^*\alpha}\Delta(\Kalpha)
+zd(\ad\imath\olmu^{(1)})\right)
\end{multline*}
We wish to show that $\ol{\nabla}H_\infty^\pm=H_\infty^\pm\ol{\nabla}^\infty$,
where 
\[\ol{\nabla}^\infty=
d_z+d_{\olh}
-\alpha_i\ad{\cow{i}}^{(1)}
-\left(\hbar(\ol{\Omega}+\Lambda)\frac{dz}{z}+
\half{\hbar}\sum_{\alpha\in\olPhi_+}\dloga\Delta(\Kalpha)+
d(z\ad\imath\olmu^{(1)})\right)\]}

(3) When $\ol{\mu}=0$, $w_\alpha=-\alpha(\imath(\mu))/\alpha
(\cow{i})=0$ for any $\alpha\in\Phi\setminus\ol{\Phi}$, and the
connection \eqref{eq:nabla zi infty} is equal to 
\[d_w-
\left(\half{\hbar}\frac{\Delta(\K-\ol{\K})}{w}+
z\ad\cowo{i}\right)dw\]
and it follows from \eqref{eq:k-olk} that $\K-\ol{\K}$ is invariant
under $\oll$.

(4) By Proposition \ref{pr:basic ODE}, $H$ possesss an asymptotic
expansion \wrt $z=\infty$, locally uniformly in $w$. Plugging $H(w,
\ol{\mu},z)=H^0(w,\ol{\mu})+H^1(w,\ol{\mu})z^{-1}+O(z^{-2})$ into
\eqref{eq:Stokes ODE}, and taking the coefficients of $z$ and $z^0$
yields $[\cowo{i},H^0]=0$ and
\[\frac{dH^0}{dw}
=
[\cowo{i},H^1]
+\half{\hbar}\left(
\sum_{\alpha\in\Phi_+\setminus\ol{\Phi}}\frac{\Delta(\Kalpha)H^0}{w-w_\alpha}
-\frac{H^0(\onetwo{(\K-\ol{\K})})}{w}\right)\]
Projecting onto the $0$ eigenspace of $\ad\cowo{i}$, and noting
that $\Delta(\Kalpha)=\onetwo{\Kalpha}+x_\alpha\otimes x_{-\alpha}
+x_{-\alpha}\otimes x_\alpha$, where $\Kalpha^{(1)},\Kalpha^{(2)}$
commute with $\cowo{i}$, and $[\cowo{i},x_{\pm\alpha}\otimes x_
{\mp\alpha}]=\pm\alpha(\cow{i})x_{\pm\alpha}\otimes x_{\mp\alpha}$
which is non--zero if $\alpha\in\Phi\setminus\ol{\Phi}$, therefore yields
\[\frac{dH^0}{dw}
=
\half{\hbar}\left(
\sum_{\alpha\in\Phi_+\setminus\ol{\Phi}}
\frac{(\onetwo{\Kalpha})H^0}{w-w_\alpha}
-\frac{H^0(\onetwo{(\K-\ol{\K})})}{w}\right)\]
which is precisely the differential equation \eqref{eq:Fuchs infty} satisfied
by the holomorphic part $H_\infty(w,\ol{\mu})^{\otimes 2}$ of $\Upsilon_
\infty(w,\ol{\mu})^{\otimes 2}$. The fact that $H^0=H_\infty^{\otimes 2}$
now follows by Remark \ref{rm:strong uniqueness}.
\end{pf}

\begin{rem}\label{rm:recurrence}
Note that, by \eqref{eq:emb i} the gauge transform of the connection
$\nabla$ by $\Upsilon^\infty_\pm$ can be written as the sum of 
commuting terms
\[d-\left(\hbar\ol{\Omega}\frac{dz}{z}-\half{\hbar}\sum_{\alpha\in\olPhi_+}\frac{d\alpha}{\alpha}\Delta(\Kalpha)-d(z\ad\olmu^{(1)})\right)
-\left(\hbar\Lambda\frac{dz}{z}-\ad\cowo{i}d(z\jmath^*\alpha_i)\right)\]
where the first two summands are the connection $\nabla$ for $\olg$,
and $\jmath:\olh\to\h$ is the canonical embedding corresponding to the
inclusion $\olg\subset\g$.
\end{rem}


\subsection{Recurrence}%
\label{ss:recurrence}

Fix $i\in\bfI$, and let $\Upsilon_\infty,\Upsilon_\infty^+$ be the solutions
of the Casimir equations given by Propositions \ref{pr:Fuchs infty} and
\ref{pr:Stokes infty} respectively. The following result relates the fusion
operators of $\g$ and $\olg$.
\comment{Since the fusion operator is initially only defined for $\mu\in
\iota\C$, only the solution $\Upsilon^+(\alpha_i,\olmu,z)$ is used for
both $\J^\pm$.}

\begin{thm}\label{th:recurrence}
The following holds for any $z\in\IR_{\gtrless 0}$ and $\mu\in\iota\C$
\[\J^\pm_\g(z,\mu)=
\Upsilon_\infty^+(\alpha_i,\ol{\mu},z)\cdot
\J^\pm_{\ol{\g}}(z,\olmu)\cdot
e^{-z\alpha_i(\olmu)\ad\cowo{i}}\cdot
(\pm z)^{\hbar\Lambda}\cdot
(\Upsilon_\infty(\alpha_i,\olmu)^{\otimes 2})^{-1}\]
\end{thm}
\begin{pf}
By construction,\comment{, and Remark \ref{rm:strong uniqueness of J},
which should be saying that $J$ is unique under the weaker assumption
that it tends to 1 if $\IR\ni z\to\pm\infty$.}
$\J^\pm_\g(z,\mu)$ is the unique solution of
\[\left(d_z-\left(\hbar\frac{\Omega}{z}+\ad\mu^{(1)}\right)dz\right)\J^\pm_\g(z,\mu)
=\J^\pm_\g(z,\mu)\,d_z\]
which is of the form
$\J^\pm_\g(z,\mu)=
H_\g^\pm(z,\mu)\cdot e^{z\ad\mu^{(1)}}\cdot(\pm z)^{\hbar\Omega_\h}$,
where
\[H_\g^\pm:\{z|\,\Re z\gtrless 0\}\times\iota\C\to\A\]
is holomorphic in $z$, and such that $H_\g^\pm(z,\mu)\to 1$ as
$\IR\ni z\to\pm\infty$.

On the other hand, applying $d_z-\left(\hbar\Omega/z+\ad\mu^{(1)}\right)dz$
to the \rhs of the stated identity yields, by Proposition \ref{pr:Stokes infty}, and
Remark \ref{rm:recurrence}
\[\begin{split}
\phantom{=}&
\Upsilon^+_\infty(\alpha_i,\ol{\mu},z)\cdot
\left(d_z
-\left(\hbar\frac{\ol{\Omega}}{z}+\ad\olmu^{(1)}\right)dz
-\left(\hbar\frac{\Lambda}{z}-\alpha_i(\olmu)\ad\cowo{i}\right)dz\right)\cdot
\J_{\ol{\g}}^\pm(z,\olmu)\\
\phantom{=}&
\phantom{\Upsilon^+_\infty(\alpha_i,\ol{\mu},z)\cdot
(d_z-(\hbar\frac{\ol{\Omega}+\Lambda}{z}+\ad(\olmu-}
\cdot e^{-z\alpha_i(\olmu)\ad\cowo{i}}\cdot
(\pm z)^{\hbar\Lambda}\cdot
(\Upsilon_\infty(\alpha_i,\olmu)^{\otimes 2})^{-1}\\[1.1ex]
=&
\Upsilon^+_\infty(\alpha_i,\ol{\mu},z)\cdot
\J_{\ol{\g}}^\pm(z,\olmu)\cdot
\left(d_z-\left(\hbar\frac{\Lambda}{z}-\alpha_i(\olmu)\ad\cowo{i}\right)dz\right)\cdot
e^{-z\alpha_i(\olmu)\ad\cowo{i}}\cdot
(\pm z)^{\hbar\Lambda}\\
\phantom{=}&
\phantom{\Upsilon^+_\infty(\alpha_i,\ol{\mu},z)\cdot
(d_z-(\hbar\frac{\ol{\Omega}+\Lambda}{z}+\ad(\olmu-}
\cdot(\Upsilon_\infty(\alpha_i,\olmu)^{\otimes 2})^{-1}\\[1.1ex]
=&
\Upsilon^+_\infty(\alpha_i,\ol{\mu},z)\cdot
\J_{\ol{\g}}^\pm(z,\olmu)\cdot
e^{-z\alpha_i(\olmu)\ad\cowo{i}}\cdot
(\pm z)^{\hbar\Lambda}\cdot
(\Upsilon_\infty(\alpha_i,\olmu)^{\otimes 2})^{-1}\,d_z
\end{split}\]
where the first equality follows from the fact that $\Lambda$
and $\cow{i}$ commute with $\olg$, and the second from the
fact that $\Upsilon_\infty(\alpha_i,\olmu)$ is independent of $z$.

Moreover, if $z\in\IR_{\gtrless 0}$, Proposition \ref{pr:Stokes infty}
(4) implies that
\[\begin{split}
&\Upsilon^+_\infty(\alpha_i,\ol{\mu},z)\cdot
\J_{\ol{\g}}^\pm(z,\olmu)\cdot
e^{-z\alpha_i(\olmu)\ad\cowo{i}}\cdot
(\pm z)^{\hbar\Lambda}\cdot
(\Upsilon_\infty(\alpha_i,\olmu)^{\otimes 2})^{-1}\\[1.1ex]
=&
(\Upsilon_\infty(\alpha_i,\olmu)^{\otimes 2}+O(z^{-1}))\cdot e^{z\alpha_i\cowo{i}}\cdot
H_\olg^\pm(z,\olmu)\cdot e^{z\ad\olmu^{(1)}}\cdot (\pm z)^{\hbar\Omega_{\olh}}\\
\phantom{=}&
\phantom{(\Upsilon_\infty(\alpha_i,\olmu)^{\otimes 2}+O(z^{-1}))\cdot e^{z\alpha_i\cowo{i}}\cdot H_\olg^\pm}
\cdot e^{-z\alpha_i(\olmu)\ad\cowo{i}}\cdot
(\pm z)^{\hbar\Lambda}\cdot
(\Upsilon_\infty(\alpha_i,\olmu)^{\otimes 2})^{-1}\\[1.1ex]
=&
(\Upsilon_\infty(\alpha_i,\olmu)^{\otimes 2}+O(z^{-1}))\cdot 
H_\olg^\pm(z,\olmu)\cdot
(\Upsilon_\infty(\alpha_i,\olmu)^{\otimes 2})^{-1}
\cdot e^{z\ad\mu^{(1)}}\cdot (\pm z)^{\hbar\Omega_\h}
\end{split}\]
where we used the fact that $\cow{i}$ commutes with $\olg$, that
$\mu=(\alpha_i-\alpha_i(\olmu))\cowo{i}+\olmu$ and that $\Omega
_{\olh}+\Lambda=\Omega_\h$ commute with $\Upsilon_\infty(\alpha
_i,\olmu)^{\otimes 2}$ since $\Upsilon_\infty(\alpha_i,\olmu)$ is of
weight $0$. The result now follows by uniqueness.
\end{pf}

\begin{rem}
The operator $\J^\pm_{\ol{\g}}(z,\olmu)\cdot e^{-z\alpha_i(\olmu)
\ad\cowo{i}}\cdot (\pm z)^{\hbar\Lambda}$ may be thought of as
the fusion operator of the Levi subalgebra $\oll=\olg+\h\subset\g$.
\end{rem}

\subsection{Centraliser property}

Retain the notation of \ref{ss:recurrence}. The following relates
the differential twists of $\g$ and $\olg$.
\comment{Note on conventions: $\J$ is initially only defined for
$\mu\in\iota\C$ but then can be continued to a single--valued
function on $\h_\IR+\iota\C$. 'Dually' the fundamental DCP
solutions are single--valued on $\C$ and can be continued
to single valued functions on $\C+i\h_\IR$ (we are using the
standard determination of the $\log$ throughout). Thus,
to compare/multiply $\J$ and the DCP solutions one can
do so on the intersection $\C+\iota\C$ of those two domains.}

\begin{thm}\label{th:centraliser}
The following holds
\[\Delta(\Upsilon_\infty(\alpha_i,\olmu))^{-1}\cdot
\J^\pm_\g(\mu)\cdot
\Upsilon_\infty(\alpha_i,\olmu)^{\otimes 2}=
C^\pm\cdot\J^\pm_\olg(\olmu)\]
where $C^\pm\in\E$ commutes with the diagonal action of $\oll$.
\end{thm}
\begin{pf}
By definition, the \lhs is equal to
\[\Delta(\Upsilon_\infty(\alpha_i,\olmu))^{-1}\cdot\Upsilon_{0,\g}(z,\mu)^{-1}\cdot
\J^\pm_\g(z,\mu)\cdot
\Upsilon_\infty(\alpha_i,\olmu)^{\otimes 2}\]
where $\Upsilon_{0,\g}(z,\mu)$ is given by Proposition \ref{pr:Fuchs 0}.
On the other hand, by Theorem \ref{th:recurrence},
\[\begin{split}
\J^\pm_\olg(\olmu)
&=
\Upsilon_{0,\olg}(z,\olmu)^{-1}\cdot\J^\pm_\olg(z,\olmu)\\
&=
\Upsilon_{0,\olg}(z,\olmu)^{-1}\cdot e^{z\alpha_i(\olmu)\ad\cowo{i}}\cdot(\pm z)^{-\hbar\Lambda}\cdot\J^\pm_\olg(z,\olmu)\cdot e^{-z\alpha_i(\olmu)\ad\cowo{i}}\cdot(\pm z)^{\hbar\Lambda}\\
&=
\Upsilon_{0,\olg}(z,\olmu)^{-1}\cdot e^{z\alpha_i(\olmu)\ad\cowo{i}}\cdot(\pm z)^{-\hbar\Lambda}\cdot
\Upsilon_\infty^+(\alpha_i,\ol{\mu},z)^{-1}\cdot
\J^\pm_\g(z,\mu)\cdot
\Upsilon_\infty(\alpha_i,\olmu)^{\otimes 2}
\end{split}\]
where the second equality follows from the fact that $\Lambda,\cow{i}$
commute with $\olg$.

We wish to compare the functions
\begin{align*}
\Upsilon_{0\infty}&=\Upsilon_{0,\g}(z,\mu)\cdot\Delta(\Upsilon_\infty(\alpha_i,\olmu))\\
\Upsilon_{\infty 0}^\pm&=
\Upsilon_\infty^+(\alpha_i,\ol{\mu},z)\cdot
e^{-z\alpha_i(\olmu)\ad\cowo{i}}\cdot(\pm z)^{\hbar\Lambda}\cdot\Upsilon_{0,\olg}(z,\olmu)
\end{align*}
It follows from Propositions \ref{pr:Fuchs 0} and \ref{pr:Fuchs infty} for
$\Upsilon_{0 \infty}$, and Propositions \ref{pr:Stokes infty} and \ref{pr:Fuchs 0}
for $\Upsilon_{\infty 0}^\pm$, that both are holomorphic functions of $\olmu\in\olh$,
which satisfy
\[\left(d_\h-\half{\hbar}\sum_{\alpha\in\Phi_+}\frac{d\alpha}{\alpha}
\Delta(\Kalpha)-z\ad(d\muone)\right)\Upsilon
=
\Upsilon
\left(d-\half{\hbar}\sum_{\alpha\in\ol{\Phi}_+}\frac{d\alpha}{\alpha}
\Delta(\Kalpha)\right)\]
and are such that their value at $\olmu=0$ commutes with $\oll$.

Set $C^\pm=\Upsilon_{0\infty}^{-1}\cdot\Upsilon_{\infty 0}^\pm$,
so that
\[\Delta(\Upsilon_\infty(\alpha_i,\olmu))^{-1}\cdot
\J^\pm_\g(\mu)\cdot
\Upsilon_\infty(\alpha_i,\olmu)^{\otimes 2}=
C^\pm\cdot\J^\pm_\olg(\olmu)\]
Then, $C^\pm$ is independent of $z$ and the coordinate $\alpha_i$,
is a holomorphic function of $\olmu\in\olh$ which satisfies
\[dC^\pm=-\half{\hbar}\sum_{\alpha\in\olPhi}\dloga[\Delta(\Kalpha),C^\pm]\]
and is such that $C^\pm(\olmu=0)$ commutes with $\oll$. We claim that
that $C^\pm$ is a constant function of $\olmu$, and therefore that it
commutes with $\oll$ for any $\olmu\in\olh$. Fix $\olmu\in\olh\reg$,
then $c(t)=C^\pm(t\olmu)$, $t\in\IC$, satisfies
\[\frac{dc}{dt}=-\half{\hbar}\frac{[\Delta(\ol{\K}),c]}{t}\]
so that $\Ad(t^{\hbar\Delta(\olK)/2})c$ is a constant $c_0$.
Since $c(t)=c^0+O(t)$, where $c^0$ commutes with $\oll$,
$c_0=c^0+\Ad(t^{\hbar\Delta(\olK)/2})O(t)=c^0+O(t)$,
whence $c_0=c^0$ and $c(t)=\Ad(t^{-\hbar\Delta(\olK)/2})c^0
=C^\pm(0)$ as claimed.\comment{There must be a very mild
sign/phase difference between $C^\pm$. Work it out.}
\end{pf}

\section{Quasi--Coxeter quasitriangular quasibialgebra structure on $\Ug$} %
\label{se:diffl qcqtqba}

The following is the main result of this paper. It shows the existence of 
a \qcqtqba structure on $\Ug\fml$ interpolating between the \qtqba
structure underlying the KZ equations and the \qca one underlying
the Casimir connection.

\begin{thm}\label{th:existence}
There exists a \qcqtqba structure 
on $\Ug\fml$ of the form
$$\Bigl(
\Ug\fml,\{\Ug_B\fml\},\{\Sicnabla{i}\},\{\Phi_{\G\F}\},
\Delta_0,
\{R_B\},\{\Phi_B\},\{F_{(B;\alpha_i)}\}
\Bigr)$$
where $\Delta_0$ is the cocommutative coproduct on $\Ug$,
\begin{align*}
\Sicnabla{i}&=\wt{s}_{i}\cdot\exp(\hbar/2\cdot C_{i}),
\\
R_B&=\exp(\hbar\cdot\Omega_B),
\\
\Alt_{2}F_{(B;\alpha_i)}&=
\hbar\cdot(\rB-r_{\g_{B\setminus\{\alpha_i\}}})
\medspace\mod\hbar^{2}
\end{align*}
and $\Phi_{\G\F}$, $F_{(B;\alpha_i)}$ are
of weight 0. Moreover, $\Phi_B\in 1^{\otimes 3}+\hbar^2\Ug_D\fml^{\otimes 3}$
is the
associator for the KZ equations corresponding to $\gB$ 
and the $\Phi_{\G\F}$ are the \DCP associators of the (truncated)
Casimir connection $\nablak$.
\end{thm}
\begin{pf}
Let $J_\pm(\mu)$ be the differential twist obtained from the fusion
operator in \ref{ss:diffl twist}. $J_\pm(\mu)$ kills the KZ associator
by Theorem \ref{th:J kills Phi} and satisfies the centraliser property
by Theorem \ref{th:centraliser}. The result now follows from Theorem
\ref{th:welding}.
\end{pf}

\section{The quantum group $\Uhg$} %
\label{se:Uhg}

\subsection{} 

Let $\Uhg$ be the \DJ quantum group corresponding to $\g$ and the bilinear
form $(\cdot,\cdot)$. Thus, $\Uhg$ is the algebra over $\IC\fml$ topologically
generated by elements $E_{i},F_{i},H_{i}$, $i\in\bfI=\{1,\ldots,r\}$, subject to
the relations\footnote{we follow here the conventions of \cite{Lu}.}
\begin{gather*}
[H_{i},H_{j}]=0
\\[1.1 ex]
[H_{i},E_{j}]= a_{ij}E_{j}
\qquad
[H_{i},F_{j}]=-a_{ij}F_{j}
\\[1.1 ex]
[E_{i},F_{j}]=
\delta_{ij}\frac
{q_i^{H_{i}}-q_i^{-H_{i}}}
{q_i-q_i^{-1}}
\intertext{where $a_{ij}=\alpha_j(\cor{i})$, $q=e^{\hbar}$, $q_i=q^{(\root{i},\root{i})/2}$, and the $q$--Serre
relations}
\sum_{k=0}^{1-a_{ij}}(-1)^{k}
\bin{1-a_{ij}}{k}{i}
E_{i}^{k}E_{j}E_{i}^{1-a_{ij}-k}=0
\\[1.1 ex]
\sum_{k=0}^{1-a_{ij}}(-1)^{k}
\bin{1-a_{ij}}{k}{i}
F_{i}^{k}F_{j}F_{i}^{1-a_{ij}-k}=0
\intertext{where for any $k\leq n$,}
[n]_{i}=\frac{q_i^{n}-q_i^{-n}}{q_i-q_i^{-1}}\\
[n]_{i}!=[n]_{i}[n-1]_{i}\cdots[1]_{i}
\qquad\text{and}\qquad
\bin{n}{k}{i}=\frac{[n]_{i}!}{[k]_{i}![n-k]_{i}!}
\end{gather*}

\subsection{$D$--algebra structure on $\Uhg$}\label{ss:D on Uhg}

For any connected $B\subseteq D$ with vertex set $\bfJ\subseteq\bfI$,
let $\Uhg_B\subseteq\Uhg$ be the subalgebra topologically generated
by the elements $\{E_j,F_j,H_j\}_{j\in\bfJ}$. Then, $\Uhg_B$ is
the \DJ quantum group corresponding to Levi subalgebra $\g_B\subset
\g$ generated by the root subspaces $\g_{\pm\alpha_j}$, $j\in\bfJ$, and
the restriction of the bilinear form $(\cdot,\cdot)$ to it. If $\bfJ=\{j\}$, we
denote $\Uhg_B$ by $\Uhsl{2}^j$.

It is clear that the assignment $B\to\Uhg_B$ defines a $D$--algebra
structure on $\Uhg$.

\subsection{Quantum Weyl group operators}\label{ss:qW}

For any $i\in\bfI$, let $\Sikh{i}$ be the operator acting on a \fd 
$\Uhg$--module $\V$ as\footnote{the element $\Sikh{i}$ is, in
the notation of \cite[\S 5.2.1]{Lu}, the operator $T_{i,+1}''$.}
\[\Sikh{i}\medspace v=
\sum_{\substack{a,b,c\in\IZ : \\a-b+c=-\lambda(\cor{i})}}
(-1)^{b}q_i^{b-ac}
E_i^{(a)}F_i^{(b)}E_i^{(c)}
v\]
where
\[E_{i}^{(a)}=\frac{E_i^{a}}{[a]_i!}
\qquad\qquad
F_{i}^{(a)}=\frac{F_i^{a}}{[a]_i!}\]
and $v\in\V$ if of weight $\lambda\in\h^*$. By \cite[\S 39.4]{Lu},
the operators $\Sikh{i}$ satisfy the braid relations
\begin{equation}\label{eq:braid reln}
\underbrace{\Sikh{i} \Sikh{j}\cdots}_{m_{ij}}=
\underbrace{\Sikh{j} \Sikh{i}\cdots}_{m_{ij}}
\end{equation}
for any $i\neq j\in\bfI$ such that the order $m_{ij}$ of $s_is_j$
is finite, and therefore define an action of the braid group $B_
W$ on $\V$.

The following modification of the operators $\Sikh{i}$ will also
be needed. Set
\[\Sich{i}=
\Sikh{i}\cdot q_i^{H_i^2/4}\]
It follows as in \ref{sss:modification} that the operators $\Sich{i}$
also satisfy the braid relations \eqref{eq:braid reln}. We refer to
either $\{\Sikh{i}\}$ or $\{\Sich{i}\}$ as the \qW operators of $\Uhg$.

The subscripts $\kappa$ and $C$ are justified by the following.
Let $\exp(\pi\iota H_i)$ and $C_i$ be the sign and Casimir
operators of $\Uhsl{2}^{i}$, that is the central elements
of $\wh{\Uhsl{2}^{i}}$ acting on the indecomposable representation
$\V_{m}$ of dimension $m+1$ as multiplication by $(-1)^{m}$
and $\half{(\alpha_i,\alpha_i)}\cdot\half{m(m+2)}$
respectively. Let $\kappa_i=C_i-(\alpha_i,\alpha_i)/4 H_i^2$
be the truncated Casimir operator of $\Uhsl{2}^{i}$.
\Omit{, that is the
operator acting on the subspace of $\V_{m}$
of weight $\ol{\jmath}=-m+2j$, $j=0\ldots m$, as multiplication by
$$q_i^{2j(m-j)+m}=
q_i^{\half{1}(m-\ol{\jmath})(m+\ol{\jmath})+m}=
q_i^{\half{m(m+2)}+\half{\ol{\jmath}}^{2}}$$}

\begin{lemma}
The following holds
\[(\Sikh{i})^2=\exp(\pi\iota H_i)\cdot q^{\kappa_i}
\aand
(\Sich{i})^2=\exp(\pi\iota H_i)\cdot q^{C_i}\]
\end{lemma}
\begin{pf}
The first identity is proved in \cite[Prop. 5.2.2.(b)]{Lu},
the second is an immediate consequence.
\end{pf}

\subsection{Quasi--Coxeter structure on $\Uhg$}\label{ss:qc on Uhg}

It follows from \ref{ss:D on Uhg} and \ref{ss:qW} that the
assignments
\[\left(\Uhg\right)_B=\Uhg_B
\qquad
S_i^{\Uhg}=\text{$\Sikh{i}$ or $\Sich{i}$}
\aand
\Phi^{\Uhg}_{\G\F}=1\]
endow $\Uhg$ with two \qca structures $\Q^\hbar_\kappa$
and $\Q^\hbar_C$ of type $D$ respectively.

\subsection{Quasitriangular quasibialgebra structure}\label{ss:DJ}

$\Uhg$ is a topological Hopf algebra 
coproduct given by
\begin{align*}
\Delta(E_{i})&=E_{i}\otimes 1+q_i^{H_{i}}\otimes E_{i}\\
\Delta(F_{i})&=F_{i}\otimes q_i^{-H_{i}}+1\otimes F_{i}\\
\Delta(H_{i})&=H_{i}\otimes 1+1\otimes H_{i}
\end{align*}

For any $B\subset D$, let
\[R\Bh\in 1^{\otimes 2}+\hbar\Uhg_B^{\otimes 2}\]
be the universal $R$--matrix of $\Uhg_B$ \cite{Dr1,Dr2}.
For $B=\alpha_i$, we denote $R\Bh$ by $R\ih$. Let
$\alpha_i\in D$, and $S\ih\in\wh{\Uhsl{2}^i}$ the \qW
element defined in \ref{ss:qW}. It follows by \cite
[Prop. 5.3.4]{Lu} that\footnote{Specifically, from
the fact that $\Delta(\Sikh{i})=
(\ol{R}\ih)^{21}\cdot\Sikh{i}\otimes\Sikh{i}$, 
where $\ol{R}\ih=q_{i}^{-\half{H_{i}\otimes H_{i}}}\cdot
R\ih$} 
\[\Delta(\Sikh{i})=
(\ol{R}\ih)^{21}\cdot\Sikh{i}\otimes\Sikh{i}\]

\subsection{}\label{ss:DL}

Label the Dynkin diagram $\Dg$ by attaching to each pair of
distinct vertices $\alpha_i\neq\alpha_j$ the order $m_{ij}$ of
the product $s_is_j\in W$. The ollowing is a direct consequence
of \ref{ss:qc on Uhg} and \ref{ss:DJ}

\begin{prop}
For any \mnss $\F,\G$ on $D$ and $\alpha_i\in B\subseteq D$,
set
\[\Phi_{\G\F}=1,
\qquad
F_{(B;\alpha_i)}=1^{\otimes 2}
\qquad\text{and}\qquad
\Phi_B=1^{\otimes 3}\]
Then,
\[\left(
\Uhg,\{\Uhg_B\},\{\Sich{i}\},\{\Phi_{\G\F}\},
\Delta,\{R\Bh\},\{F_{(B;\alpha_i)}\},\{\Phi_B\}
\right)\]
is a \qcqtqba structure $\QQch$ of type $D$ on $\Uhg$.
\end{prop}

\section{The monodromy theorem}\label{se:monodromy}

\subsection{}\label{se:main qcqtqba}

Let $\QQch$ be the \qcqtqba structure on the quantum group $\Uhg$
obtained in Section \ref{se:Uhg}, and underlying its $R$--matrix and
\qW representations. Let $\Q\Q_C^\nabla$ be the \qcqtqba structure
on $\Ug\fml$ obtained in Section \ref{se:diffl qcqtqba} which underlies
the monodromy of the KZ and Casimir connections of $\g$.

\begin{thm}\label{th:main qcqtqba}
$(\Uhg,\QQch)$ and $(\Ug\fml,\QQdh)$ are isomorphic as
quasi--Coxeter, quasitriangular quasibialgebras.
\end{thm}
\begin{pf}
By \cite[Thm. 8.3]{TL2}, $\QQch$ is isomorphic to a \qcqtqba
structure of type $D$ on $\Ug\fml$ of the form
\[\Bigl(
\Ug\fml,\{\Ug_B\fml\},\{S\ic\},\{\Phi_{\G\F}\},
\Delta_{0},\{\Phi_B\},\{R\BKZ\},\{F_{(B;\alpha_i)}\}
\Bigr)\]
where $\Delta_{0}$ is the cocommutative coproduct
on $\Ug$,
\begin{align*}
S\ic		&= \wt{s}_{i}\cdot\exp(\hbar/2\cdot C_i)\\
\Phi_B	&= 1^{\otimes 3}\mod\hbar^2\\
R\BKZ	&= \exp(\hbar\cdot\Omega_B)\\
\Alt_{2}F_{(B;\alpha_i)}&=
\hbar\cdot(\rB-r_{\g_{B\setminus\{\alpha_i\}}})
\medspace\mod\hbar^{2}
\end{align*}
and $\Phi_{\G\F}$, $F_{(B;\alpha_i)}$ are of weight 0. Since
$\QQdh$ is also of this form by Theorem \ref{th:existence},
the result follows from the rigidity of such structures \cite
[Thm. 9.1]{TL2}.
\end{pf}

\subsection{Monodromy}

\begin{thm}\hfill
\begin{enumerate}
\item $(\Ug\fml,\Qnablac)$ and $(\Uhg,\Qhc)$ are isomorphic as quasi--Coxeter algebras.
In particular, if $V$ is a \fd $\g$--module, the monodromy of the Casimir connection $\nablac$
on $V\fml$ is equivalent to the action of the braid group $B_W$ on any quantum deformation
of $V$ given by Lusztig's \qW operators $\Sich{i}$.
\item $(\Ug\fml,\Qnablak)$ and $(\Uhg,\Qhk)$ are isomorphic as quasi--Coxeter algebras.
In particular, if $V$ is a \fd $\g$--module, the monodromy of the (truncated) Casimir connection
$\nablak$ on $V\fml$ is equivalent to the action of the braid group $B_W$ on any quantum
deformation of $V$ given by Lusztig's \qW operators $\Sikh{i}$.
\end{enumerate}
\end{thm}
\begin{pf}
The first statement is an immediate consequence of Theorem \ref{th:main qcqtqba}. The
isomorphism of \qc structures $\Qnablac\cong\Qhc$ gives rise to one between $\Qnablak$
and $\Qhk$ since, by construction, the \DCP associators of $\Qnablac$ are the same as
those of $\Qnablak$, and the isomorphism is equivariant for the $\h$--action, and therefore
compatible with the modifications
\[\Sich{i}=\Sikh{i}\cdot q_i^{H_i^2/4}
\aand
\Sicnabla{i}=\Siknabla{i}\cdot\exp\left(\half{\hbar}\half{(\alpha_i,\alpha_i)}h_i^2\right)\]
\end{pf}

\appendix

\section{The basic ODE}\label{se:basic ODE} %

Fix $R\geq 0$, and let
\[\IH^R_\pm=\{z\in\IC|\,\text{$\Im z\gtrless 0$ and $|z|>R$}\}\]
be the complements in the upper and lower half--planes $\IH_\pm$
of the closed disk of radius $R$.
Let $V$ be a \fd complex vector space, and $k:\IH^R_\pm\to V$
a holomorphic function possessing an asymptotic expansion of the form
$k\sim\sum_{n\geq 0}k_nz^{-n}$ on $\IH^R_\pm$. Thus, for any $\delta>0$ and
$n\in\IN$, there is a constant $C=C(\delta,n)$ such that, for any $z\in\IH^R
_\pm$ with $\delta\leq|\arg(z)|\leq\pi-\delta$, the following holds
\[\|k(z)-\sum_{m=0}^nk_mz^{-m}\|\leq C|z|^{-(n+1)}\]

Let $\lambda\in\IR$ and consider the inhomogeneous ODE
\Omit{\comment{this ODE has an irregular singularity at $z=\infty$.
Given the simple contour integral expression of its solutions, it might
be possible to explicitly compute its Stokes factors in terms of the
change of asymptotics in $k$ across the real axis. Not entirely
clear this makes sense since the ODE is not homogeneous.}}

\begin{equation}\label{eq:basic}
\frac{dh}{dz}=\lambda h+\frac{k}{z}
\end{equation}
We seek holomorphic solutions $h:\IH^R_\pm\to V$ satisfying
the boundary condition 
\begin{equation}\label{eq:boundary}
\text{$h(z)\to 0$ as $z\to\infty$ on any sector $\delta<|\arg(z)|<\pi-\delta$}
\end{equation}

For any $\theta\in[-\pi,\pi]$, let $\Gamma_\pm^\theta$ be the
contour given by the ray $\pm\iota R+e^{\iota \theta}\cdot\IR_
{\geq 0}$, oriented from $\infty$ to $\pm\iota R$, followed by
the interval from $\pm\iota R$ to $0$.

\begin{prop}\label{pr:basic ODE}\hfill
\begin{enumerate}
\item The equation \eqref{eq:basic} has at most one solution
such that \eqref{eq:boundary} holds.

\item If $k_0\neq 0$ and $\lambda=0$, no solution to \eqref
{eq:basic}--\eqref{eq:boundary} exists.

\item If $k_0=0$ or $\lambda\neq 0$, \eqref{eq:basic}--\eqref
{eq:boundary} have a unique solution given by the Laplace
integral
\begin{equation}\label{eq:Laplace}
h(z)=\int_\Gamma e^{-\lambda t}\frac{k(z+t)}{z+t}dt
\end{equation}
where $\Gamma$ is a path from $\infty$ to $0$ such that
$\Gamma+\IH_\pm^R\subset\IH_\pm^R$, and such that
the integral is convergent as $t$ tends to infinity along
$\Gamma$. Specifically,

\begin{enumerate}
\item If $k_0\neq 0$ and $\lambda\neq 0$, $\Gamma$
can be chosen as $\Gamma_\pm^{\pm\theta}$, with
\[\pi>\theta>\pi/2\text{  if  }\lambda<0
\quad\text{and}\quad 
\pi/2>\theta>0\text{  if  }\lambda>0\]
As a function of $\lambda$, $h$ is smooth on $\IR^*$.

\item If $k_0=0$, $\Gamma$ can be chosen as $\Gamma_\pm
^{\pm\theta}$, with
\[\pi>\theta\geq \pi/2\text{  if  }\lambda\leq 0
\quad\text{and}\quad 
\pi/2\geq\theta>0\text{  if  }\lambda\geq 0\]
As a function of $\lambda$, $h$ is continuous on $\IR$, and
smooth on $\IR^*$.
\end{enumerate}

\item If $\lambda\neq 0$, $h$ has an asymptotic expansion
\wrt $z$ which is given by
\begin{equation}\label{eq:asy h}
h(z)\sim\sum_{n\geq 1}\frac{(n-1)!}{z^n}
\left(\sum_{p=0}^{n-1}\frac{k_p}{p!}\frac{1}{(-\lambda)^{n-p}}\right)
\end{equation}
and is valid uniformly on compact subsets of $\IR^*\ni\lambda$.

\item If $\lambda=0$ and $k_0=0$, $h$ has an asymptotic expansion
given by\comment{The continuity of $h$ at $\lambda=0$
does not contradict the discontinuity of asymptotic expansions \wrt $z$
since a function $h_\lambda(z)$ can depend continuously upon a parameter
$\lambda$ and have an asymptotic expansion \wrt z which is discontinuous
in $\lambda$, \eg $h_\lambda(z)=e^{-\lambda/z}$.}
\[h(z)\sim -\sum_{n\geq 1}k_n\frac{z^{-n}}{n}\]

\item As a function of $\lambda$, $h$ has an asymptotic expansion
as $\lambda\to\infty$ given by
\[h\sim-\sum_{n\geq 1}\left(\frac{k(z)}{z}\right)^{(n-1)}\lambda^{-n}\]
which is valid uniformly on compact subsets of $\IH^R_\pm$.
\end{enumerate}
\end{prop}
\begin{pf}
(1) The solutions of the underlying homogeneous equation are given
by $g=e^{\lambda t}g_0$, where $g_0\in V$ is a constant. If $\lambda
=0$, then $g\equiv g_0$ is equal to its limit as $z\to\infty$, and $g=0$.
If $\lambda\neq 0$, the function $e^{\lambda t}$ does not have a limit
as $z\to\infty$ in $\IH^R_\pm$, whence $g_0=0$.

(2) Write $k=k_0+\ol{k}$, where $\ol{k}=k-k_0=O(z^{-1})$. Any solution
of \eqref{eq:basic} with $\lambda=0$ is of the form $C+k_0\ln z-\int_
\gamma\frac{\ol{k}(t)}{t}dt$ where $C$ is a constant, and $\gamma$
is a path in $\IH^R_\pm$ with a fixed starting point and ending in $z$.
Since $\int_\gamma\frac{\ol{k}(t)}{t}dt=O(z^{-1})$, no such solution
admits a limit as $z\to\infty$ if $k_0\neq 0$.

(3) Integration by parts readily shows that the integral \eqref{eq:Laplace}
is a solution of \eqref{eq:basic}. If $k_0\neq 0$ and $\lambda\neq 0$,
$\Gamma_\pm^{\pm\theta}$ satisfies the required conditions since
integrability at $\infty$ is guaranteed by the exponential factor $e^
{-\lambda\Re t}$, so long as $\lambda\cos\theta>0$. If, on the other
hand, $k_0=0$, integrability is guaranteed by $k(z)/z=O(z^{-2})$,
so long as $e^{-\lambda\Re t}$ remains bounded, and the given
$\Gamma_\pm^{\pm\theta}$ satisfy the required conditions.

It is clear that the function defined by \eqref{eq:Laplace} is a
smooth function of $\lambda\in\IR^*$, since the convergence
of the derivatives of integral is guaranteed by the factor $\exp
(-\lambda\Re t)$. If $k_0=0$, we may choose $\Gamma=
\Gamma_\pm^{\pm\theta/2}$ for any $\lambda\in\IR$, and the
continuity of $h(z)$ at $\lambda=0$ follows from the
Riemann--Lebesgue Lemma.

(4) Set $g(\zeta)=k(\zeta)/\zeta$. If $\lambda\neq 0$, integration
by parts shows that, for any $m\geq 0$,
\begin{equation}\label{eq:h exp}
h(z)
=
\int_\Gamma e^{-\lambda t}g(z+t)dt\\
=
-\sum_{p=0}^{m-1}\frac{1}{\lambda^{p+1}}g^{(p)}(z)+
\frac{1}{\lambda^m}\int_\Gamma e^{-\lambda t}g^{(m)}(z+t)dt
\end{equation}
Let $k(z)\sim \sum_{n\geq 0} k_nz^{-n}$ be the asymptotic
expansion of $k$. Then, for any $p\geq 0$
\[g^{(p)}(z)\sim
(-1)^p\sum_{n\geq 0}k_n\frac{(n+p)!}{n!}z^{-(n+p+1)}=
(-1)^p\sum_{n\geq p+1}k_{n-p-1}\frac{(n-1)!}{(n-p-1)!}z^{-n}\]
It follows that
\[\begin{split}
-\sum_{p=0}^{m-1}\frac{1}{\lambda^{p+1}}g^{(p)}(z)
&=
-\sum_{p=0}^{m-1}\frac{(-1)^p}{\lambda^{p+1}}
\sum_{n=p+1}^m k_{n-p-1}\frac{(n-1)!}{(n-p-1)!}z^{-n}+O(z^{-m-1})\\
&=\phantom{-}
\sum_{n=1}^m z^{-n}(n-1)!\sum_{p=1}^n
\frac{(-1)^p}{\lambda^p}\frac{k_{n-p}}{(n-p)!}+O(z^{-m-1})
\end{split}\]
locally uniformly in $\lambda\in\IR^*$.

To estimate the second summand in \eqref{eq:h exp}, note
that
\[\begin{split}
\left|\int_\Gamma e^{-\lambda t}g^{(m)}(z+t)dt\right|
&\leq \int_0^\infty e^{-\lambda\Re\varphi(s)}\left|g^{(m)}(z+\varphi(s))\right||\varphi'(s)|ds\\
&\leq C\, d(z,-\Gamma)^{m+1}\int_0^\infty e^{-\lambda\Re\varphi(s)}ds\\
&=C\, d(z,-\Gamma)^{m+1}\left(R+\frac{1}{\lambda\cos\theta}\right)
\end{split}\]
where $t=\varphi(s)$ is the parametrisation of $\Gamma=\Gamma
_\pm^\theta$ with the opposite orientation given by $\varphi(s)=\pm
\imath s$ for $s\in[0,R]$ and $\varphi(s)=\pm\imath R+(s-R)e^{\imath
\theta}$ for $s\geq R$, and the second inequality follows from the
fact that $g^{(m)}(z)=O(z^{-m-1})$.

Let $\Gamma_\pm^\theta\subset K^\theta_\pm\subset\IH_\pm$ be
the convex cone bounded by the rays $e^{\imath\theta}$ and $e^
{\pm\imath(\pi-\theta)}$, and $S_\pm^{|z|}\subset\IH_\pm$ the half--circle
with radius $|z|$. Then, $d(z,-\Gamma^\theta_\pm)\geq d(S_\pm
^{|z|},-K^\theta_\pm)\geq |z||\sin\theta|$, where the bound is
attained when $w\in S_\pm^{|z|}$ lies on the real axis. It follows that
\begin{equation}\label{eq:summand}
\left|\int_\Gamma e^{-\lambda t}g^{(m)}(z+t)dt\right|
\leq C\, |z|^{-m-1}|\sin\theta|\left(R+\frac{1}{\lambda\cos\theta}\right)
\end{equation}
where $C$ is independent of $\lambda$.

(5) If $k_0=0$ and $\lambda=0$, then
\[h=\int_\Gamma\frac{k(z+t)}{z+t}
\sim\sum_{n\geq 1}\int_\Gamma k_n(z+t)^{-n-1}dt=
-\sum_{n\geq 1}\frac{k_n}{n}z^{-n}\]

(6) The existence of the claimed asymptotic expansion of $h$ \wrt
$\lambda$ is guaranteed by \eqref{eq:h exp}, provided $\int_\Gamma
e^{-\lambda t}g^{(m)}(z+t)dt=O(\lambda^{-1})$. Integrating by parts,
we have
\[\int_\Gamma e^{-\lambda t}g^{(m)}(z+t)dt=
\frac{1}{\lambda}\left(-g^{(m)}(z)+\int_\Gamma e^{-\lambda t}g^{(m+1)}(z+t)dt\right)\]
and the required estimate now follows from \eqref{eq:summand}.
\end{pf}

\section{The constant $C^\pm$ revisited}\label{app:revisited}  %

The goal of this section is to show that the constant $C^\pm$ relating
the differential twists of $\g$ and $\olg$ given by Theorem \ref{th:centraliser}
can be computed as the monodromy from $0$ to $\infty$ of an ODE
on $\IP^1$ with regular singularities at $0,1$ and an irregular singularity
at $\infty$. This gives in particular a canonical, transcendental construction
of a relative twist for the pair of KZ associators $(\Phi,\olPhi)$, \ie an
element $J\in 1+\hbar(\Ug[[\hbar]]^{\otimes 2})^\olg$ such that $(\Phi)
_{J}=\olPhi$, in the spirit of Drinfeld's construction of the KZ associator.
\Omit{Highlight the killing of \Phi\KKZ in a corollary. Maybe pose the
problem of computing it explicitly.}

\subsection{} 

Consider the connection $\nabla$ on $\IC\times\h$ given by \eqref
{eq:DCKZ 2}. Fix $\olmu\in\olh$, and coordinatise $\pi^{-1}(\olmu)$
by $w=\alpha_i$ as in \ref{ss:Stokes}. The restriction of $\nabla$ to
$\IC\times\pi^{-1}(\olmu)$ then reads
\[\nabla=d
-\hbar\Omega\frac{dz}{z}
-\frac{\hbar}{2}\sum_{\alpha\in\Phi_+\setminus\ol{\Phi}_+}\frac{\Delta(\Kalpha)}{w-w_\alpha}dw
-\ad\cowone{i}d(zw)-\ad\imath(\olmu)^{(1)}dz\]
Our first goal is to construct two canonical solutions of $\nabla$ with
prescribed asymptotic behaviour as $z\to 0$, $w\to\infty$, and $zw
\to 0,\infty$ respectively.

\subsection{Blow--up coordinates}

To this end, let $\rho:(v_1,v_2)\in\IC\times\IP^1\to(z,w)\in\IC\times\IP^1$
be the rational map given by \footnote{note that, unlike previous changes
of coordinates, $\rho$ mixes the configuration coordinate $z$ with the
Cartan coordinate $w$.}
\begin{equation}\label{eq:blowup} 
z=\frac{v_1v_2}{v_2+1}\aand w=\frac{v_2+1}{v_1}
\end{equation}
$\rho$ is a birational isomorphism, with inverse given by
\begin{equation}\label{eq:inv blowup}
v_1=z+1/w\aand v_2=zw
\end{equation}
and restricts to an isomorphism of $\{(v_1,v_2)\in\IC\times\IP^1|\,v_1\neq 0,
v_2\neq -1\}$ onto $\{(z,w)\in\IC\times\IP^1|\,w\neq 0,(z,w)\neq(0,\infty),zw
\neq -1\}$. Since $v_2=zw$, the asymptotic zones $z\sim 0$, $w\sim\infty$
and $zw\sim 0,\infty$ correspond respectively to the neighborhoods of the
points $(v_1,v_2)=(0,0)$ and $(0,\infty)$.

Given that $w-w_\alpha=(v_2+1-w_\alpha v_1)/v_1$, we get
\[d\log (w-w_\alpha)=-\frac{dv_1}{v_1}-w_\alpha\frac{dv_1}{v_2+1-w_\alpha v_1}+\frac{dv_2}{v_2+1-w_\alpha v_1}\]
It follows that the pulled--back connection $\rho^*\nabla$ is given by
\[\begin{split}
\rho^*\nabla=d
&-\left(\frac{\hbar}{2}\frac{2\Omega-\Delta(\K-\ol{\K})}{v_1}
-\frac{\hbar}{2}\sum_{\alpha\in\Phi_+\setminus\ol{\Phi}_+}\frac{w_\alpha\Delta(\Kalpha)}{v_2+1-w_\alpha v_1}
+\frac{v_2}{v_2+1}\ad\imath(\olmu)^{(1)}\right)dv_1\\
&-\left(\frac{\hbar\Omega}{v_2}
-\frac{\hbar\Omega}{v_2+1}
+\frac{\hbar}{2}\sum_{\alpha\in\Phi_+\setminus\ol{\Phi}_+}\frac{\Delta(\Kalpha)}{v_2+1-w_\alpha v_1}
+\ad{\cow{i}}^{(1)}
+\frac{v_1}{(v_2+1)^2}\ad\imath(\olmu)^{(1)}
\right)dv_2
\end{split}\]
where $\K,\ol{\K}$ are given by \eqref{eq:K olK}.

\subsection{} 

The proof of the following result is similar to that of Proposition \ref{pr:Fuchs 0}
and therefore omitted.

\begin{prop}\label{pr:core}\hfill
\begin{enumerate}
\item
For any $\olmu\in\olh$ and $v_2\in\IP^1\setminus\{-1\}$, there is a unique
holomorphic function
\[I_0:\{v_1\in\IC|\,|v_1|<R_\olmu^{-1}\cdot|v_2+1|\}\to\A\]
such that $I_0(0,v_2,\ol{\mu})=1$ and, for any determination of $\log(v_1)$,
the $\E$--valued function
\[\Xi_0(v_1,v_2,\ol{\mu})=
e^{\frac{v_1v_2}{v_2+1}\ad\imath(\olmu)^{(1)}}\cdot
I_0(v_1,v_2,\ol{\mu})\cdot
v_1^{-\half{\hbar}({\Delta(\K-\ol{\K})-2\Omega})}\]
satisfies $\rho^*\nabla_{\partial_{v_1}}\Xi_0=\Xi_0\,d_{v_1}$, where
$\rho^*\nabla_{\partial_{v_1}}$ is given by
\begin{equation}\label{eq:core conn}
d_{v_1}-
\left(\frac{\hbar}{2}\frac{2\Omega-\Delta(\K-\ol{\K})}{v_1}
-\frac{\hbar}{2}\sum_{\alpha\in\Phi_+\setminus\ol{\Phi}_+}\frac{w_\alpha\Delta(\Kalpha)}{v_2+1-w_\alpha v_1}
+\frac{v_2}{v_2+1}\ad\imath(\olmu)^{(1)}\right)dv_1
\end{equation}
\item $I_0$ and $\Xi_0$ are holomorphic functions of $v_2$ and $\olmu$,
and $\Xi_0$ satisfies
\[\rho^*\nabla\,\Xi_0=\Xi_0\,\left(d_{v_1}+\nabla_0+\ol{\nabla}_C\right)\]
where
\begin{align}
\nabla_0
&=
d_{v_2}-
\left(\frac{\hbar\Omega}{v_2}
+\frac{\hbar}{2}\frac{\Delta(\K-\ol{\K})-2\Omega}{v_2+1}
+\ad\cowone{i}
\right)dv_2
\label{eq:nabla0}\\
\ol{\nabla}_C
&=
d_{\ol{\h}}
-\half{\hbar}\sum_{\alpha\in\ol{\Phi}_+}\frac{d\alpha}{\alpha}\Delta(\Kalpha)
\label{eq:nablacc}
\end{align}
\end{enumerate}
\end{prop}

\subsection{} 

Let $\nabla_0$ be the connection given by \eqref{eq:nabla0}.

\begin{prop}\label{pr:Fuchs and Stokes}\hfill
\begin{enumerate}
\item The coefficients of the connection $\nabla_0$ commute with
the diagonal adjoint action of the Levi subalgebra $\oll$ and with
$\Delta(\K-\ol{\K})-2\Omega$.
\item There is a unique holomorphic function \[G_0:\{v_2\in\IC|\,|v_2|<1\}\to\A\]
such that $G_0(0)=1$ and, for any determination of $\log v_2$, the
$\E$--valued function
\[\Psi_0=
e^{v_2\ad{\cow{i}}^{(1)}}\cdot G_0(v_2)\cdot v_2^{\hbar\Omega}\]
satisfies $\nabla_0\Psi_0=\Psi_0 d_{v_2}$.
\item There is a unique holomorphic function
\[G_\infty^\pm:\{v_2\in\IC|\,\Im v_2\gtrless 0, |v_2|>1\}\to\A\]
such that $G_\infty^\pm(v_2)\to 1$ as $v_2\to\infty$ with $0<<|\arg(v_2)|
<<\pi$, and, for any determination of $\log(v_2)$, the function 
\[\Psi_\infty^\pm=
G_\infty^\pm(v_2)
\cdot e^{v_2\ad\cowone{i}}
\cdot v_2^{\half{\hbar}\left((\K-\ol{\K})^{(1)}+(\K-\ol{\K})^{(2)}\right)}\]
satisfies $\nabla_0\Psi_\infty^\pm=\Psi_\infty^\pm d_{v_2}$.
\item The functions $\Psi_0,\Psi_\infty^\pm$ commute with the diagonal
adjoint action of $\oll$ and with $\Delta(\K-\ol{\K})-2\Omega$.
\end{enumerate}
\end{prop}
\begin{pf}
(1) Since $\Omega$ commutes with the diagonal action of $\g$,
and $\cow{i}$ commutes with $\ol{\g}$, it suffices to show that
$\Delta(\K-\ol{\K})-2\Omega$ commutes with $\oll$ and with
$\cowo{i}$. The commutation with $\cowo{i}$ follows from the
fact that
\begin{equation}
\begin{split}
\label{eq:Delta K olK}
\Delta(\K-\ol{\K})
&=
(\K-\ol{\K})^{(1)}+(\K-\ol{\K})^{(2)}+
2\sum_{\alpha\in\Phi\setminus\ol{\Phi}}x_\alpha\otimes x_{-\alpha}\\
&=
(\K-\ol{\K})^{(1)}+(\K-\ol{\K})^{(2)}+2\left(\Omega-\ol{\Omega}-\Lambda\right)
\end{split}
\end{equation}
The fact that $\K-\olK$ commutes with $\oll$ follows from \eqref{eq:k-olk}.

(2) is proved in a similar way to Proposition \ref{pr:Fuchs 0}.

(3) is proved in a similar way to Proposition \ref{pr:Stokes infty}, and
relies on the fact that the connection $\nabla_0$ is of the form
\[d_{v_2}-\left(\ad\cowone{i}+\half{\hbar}\frac{\Delta(\K-\olK)}{v_2}+O(v_2^{-2})\right)dv_2\]
and that, by (1), the projection of $\Delta(\K-\olK)$ onto the kernel
of $\ad\cowone{i}$ is $(\K-\ol{\K})^{(1)}+(\K-\ol{\K})^{(2)}$.

(4) follows from (1). 
\end{pf}

\subsection{}

Let $C^\pm\in\E$ be the constant relating the differential twists of
$\g$ and $\olg$ given by Theorem \ref{th:centraliser}, $\nabla_0$
the connection \eqref{eq:nabla0}, and $\Psi_0,\Psi_\infty^\pm$ its
fundamental solutions given by Proposition \ref{pr:Fuchs and Stokes}.

\begin{thm}
The following holds,
\[C^\pm=\Psi_0^{-1}\cdot \Psi_\infty^\pm\]
\end{thm}
\begin{pf}
%
The proof of Theorem \ref{th:centraliser} shows that $C^\pm=
\Upsilon_{0\infty}^{-1}\cdot\Upsilon_{\infty 0}^\pm$, where
\begin{align*}
\Upsilon_{0\infty}&=\Upsilon_{0,\g}(z,\mu)\cdot\Delta(\Upsilon_\infty(\alpha_i,\olmu))\\
\Upsilon_{\infty 0}^\pm&=
\Upsilon_\infty^+(\alpha_i,\ol{\mu},z)\cdot
e^{-z\alpha_i(\olmu)\ad\cowo{i}}\cdot(\pm z)^{\hbar\Lambda}\cdot\Upsilon_{0,\olg}(z,\olmu)
\end{align*}
and $\Upsilon_{0,\g}(z,\mu)$, $\Upsilon_\infty(\alpha_i,\olmu)$, $\Upsilon
_\infty^+(\alpha_i,\ol{\mu},z)$ are the functions given by Propositions \ref
{pr:Fuchs 0}, \ref{pr:Fuchs infty} and \ref{pr:Stokes infty} respectively.
%
Let $\Xi_0$ and $\Psi_0,\Psi_\infty^\pm$ be the functions given by
Propositions \ref{pr:core} and \ref{pr:Fuchs and Stokes} respectively.
We claim that
\[\Upsilon_{0\infty}=\Xi_0\cdot\Psi_0
\aand
\Upsilon_{\infty 0}^\pm=\Xi_0\cdot\Psi_\infty^\pm\]
so that $C^\pm=\Psi_0^{-1}\cdot \Psi_\infty^\pm$ holds.

To see that the first claimed identity holds, write
\[\Upsilon_{0\infty}=
\wt{H}_{0,\g}(z,w)\cdot z^{\hbar\Omega}\cdot \Delta(H_\infty(w))\cdot w^{\half{\hbar}\Delta(\K-\olK)}\]
where $\wt{H}_{0,\g}=e^{z\ad\muone}H_{0,\g}$, and we have suppressed the
dependence in $\olmu\in\olh$ which will be held fixed throughout the
argument. We have
\[\begin{split}
\Xi_0\cdot\Psi_0
&=
\wt{I}_0(v_1,v_2)\cdot v_1^{-\half{\hbar}(\Delta(\K-\olK)-2\Omega)}\cdot
G_0(v_2)\cdot v_2^{\hbar\Omega}\\
&=
J_0(z,w)\cdot \wt{I}_0(1/w,0)\cdot w^{\half{\hbar}\Delta(\K-\olK)}\cdot z^{\hbar\Omega}\\
&=
J_0(z,w)\cdot z^{\hbar\Omega}\cdot \wt{I}_0(1/w,0)\cdot w^{\half{\hbar}\Delta(\K-\olK)}
\end{split}\]
where $\wt{I}_0=e^{\frac{v_1v_2}{v_2}\ad\imath(\olmu)^{(1)}}\cdot I_0$, $v_1,v_2$
are expressed in terms of $z,w$ through \eqref{eq:inv blowup}, the function
$J_0(z,w)=
\wt{I}_0(v_1,v_2)\cdot v_1^{-\half{\hbar}(\Delta(\K-\olK)-2\Omega)}\cdot G_0(v_2)
\cdot w^{-\half{\hbar}(\Delta(\K-\olK)-2\Omega)}
\cdot \wt{I}_0(1/w,0)^{-1}$ is holomorphic near $z=0$ for fixed $w$ and such that $J_0
(0,w)=1$, and the third equality follows from the fact $\wt{I}_0(v_1,0)$ commutes with
$\Omega$ since the coefficients of the connection \eqref{eq:core conn} of which $\wt{I}
_0(v_1,v_2)$ is a horizontal section do for $v_1=0$. Since $\Xi_0\cdot\Psi_0$ is a horizontal
section of $\nabla_{\partial_z}$, so is $J_0(z,w)\cdot z^{\hbar\Omega}=\Xi_0\cdot\Psi_0
\cdot\left(\wt{I}_0(1/w,0)\cdot w^{\half{\hbar}\Delta(\K-\olK)}\right)^{-1}$ and it follows
by the uniqueness of Proposition \ref{pr:Fuchs 0} that $J_0(z,w)=\wt{H}_{0,\g}(z,w)$.
We are therefore reduced to proving that $\Delta(H_\infty(w))\cdot w^{\half{\hbar}
\Delta(\K-\olK)}=\wt{I}_0(1/w,0)\cdot w^{\half{\hbar}\Delta(\K-\olK)}$, which follows
from the uniqueness of Proposition \ref{pr:Fuchs infty}.\\

%
The second identity is proved in a similar way. We may assume that $z\gtrless 0$.
\comment{I think I actually need to assume that $z>0$.} We have
\begin{multline*}
\Upsilon_{\infty 0}^\pm=
H_\infty^\pm(w,z)\cdot e^{zw\ad\cowone{i}}\cdot w^{\half{\hbar}(\onetwo{(\K-\olK)})}\\
\cdot e^{-z\alpha_i(\olmu)\ad\cowo{i}}\cdot(\pm z)^{\hbar\Lambda}
\cdot H_{0,\olg}(z,\olmu)\cdot z^{\hbar\ol{\Omega}}
\end{multline*}
and
\[\begin{split}
\Xi_0\cdot\Psi_\infty^\pm
&=
\wt{I}_0(v_1,v_2)
\cdot v_1^{-\half{\hbar}(\Delta(\K-\olK)-2\Omega)}
\cdot G_\infty^\pm(v_2)
\cdot e^{v_2\ad\cowone{i}}\cdot v_2^{\half{\hbar}(\onetwo{(\K-\olK)})}\\
&=
\wt{I}_0(v_1,v_2)
\cdot G_\infty^\pm(v_2)
\cdot v_1^{-\half{\hbar}(\Delta(\K-\olK)-2\Omega)}
\cdot e^{zw\ad\cowone{i}}\cdot (zw)^{\half{\hbar}(\onetwo{(\K-\olK)})}\\
&=
\wt{J}_0(w,z)
\cdot \wt{I}_0(z,\infty)
\cdot z^{-\half{\hbar}(\Delta(\K-\olK)-2\Omega)}
\cdot e^{zw\ad\cowone{i}}
\cdot (zw)^{\half{\hbar}(\onetwo{(\K-\olK)})}\\
&=
\wt{J}^\pm_0(w,z)
\cdot e^{zw\ad\cowone{i}}
\cdot w^{\half{\hbar}(\onetwo{(\K-\olK)})}
\cdot \wt{I}_0(z,\infty)
\cdot z^{\hbar(\ol{\Omega}+\Lambda)}
\end{split}\]
where the second equality follows from the fact that $G_\infty^\pm$
commutes with $\Delta(\K-\olK)-2\Omega$ by Proposition \ref{pr:Stokes infty},
the function $\wt{J}^\pm_0(w,z)$ is defined by 
\begin{multline*}
\wt{J}^\pm_0(w,z)=
\wt{I}_0(v_1,v_2)
\cdot G_\infty^\pm(v_2)
\cdot v_1^{-\half{\hbar}(\Delta(\K-\olK)-2\Omega)}
\cdot z^{\half{\hbar}(\Delta(\K-\olK)-2\Omega)}
\cdot \wt{I}_0(z,\infty)^{-1}
\end{multline*}
and, for fixed $z$, tends to $1$ as $w\to\infty$ with $\Im w\gtrless 0$
and $0<<|\arg w|<<\pi$, and the last equality uses the fact that 
$\wt{I}_0(z,\infty)$ commutes with $\cowone{i}$ and $\onetwo{(\K-\olK)}$
since the coefficients of the connection \eqref{eq:core conn} do for
$v_2=\infty$, and \eqref{eq:Delta K olK}.
By the uniqueness of Proposition \ref{pr:Stokes infty}, $\wt{J}^\pm_0(w,z)
=H_\infty^\pm(w,z)$, which reduces the stated claim to proving that
\[e^{-z\alpha_i(\olmu)\ad\cowo{i}}\cdot(\pm z)^{\hbar\Lambda}
\cdot H_{0,\olg}(z,\olmu)\cdot z^{\hbar\ol{\Omega}}=
\wt{I}_0(z,\infty)
\cdot z^{\hbar(\ol{\Omega}+\Lambda)}\]
In turn, this follows from the uniqueness of Proposition \ref{pr:Fuchs 0}.
\end{pf}

\section{The dynamical KZ and Casimir equations}\label{app:joint}

For any $n\geq 2$, let
\[\gX_n=\IC^n\setminus\bigcup_{1\leq i<j\leq j}\{z_i=z_j\}\]
be the configuration space of $n$ ordered points in $\IC$. Consider the following
connection on the trivial vector bundle over $\gX_n\times\hreg$ with fibre $\Ug^
{\otimes n}$
\[\nabla=d
-\sfh\sum_{1\leq i<j\leq n}\frac{d(z_i-z_j)}{z_i-z_j}\Omega_{ij}
-\frac{\sfh}{2}\sum_{\alpha\in\Phi_+}\frac{d\alpha}{\alpha}
\Delta^{(n)}(\Kalpha)-td(\sum_{i=1}^n z_i\mu^{(i)})\]
Above, $\sfh,t$ are complex parameters, 
$\Omega_{ij}=\sum_a X_a^{(i)}{X^a}^{(j)}$, where $X^{(i)}=1^{\otimes(i-1)}\otimes
X\otimes 1^{\otimes (n-i)}$, $\Delta^{(n)}:\Ug\to\Ug^{\otimes n}$ is the iterated coproduct
and $\mu$ is the embedding $\h\to\Ug$.

The coupling term $d(\sum_{i=1}^n z_i\mu^{(i)})$ was shown to lead to a consistent 
system in \cite{FMTV} when $t=\sfh$. In the present paper, it is important to keep
$\sfh$ and $t$ as independent parameters, since we consider $\sfh$ formal while
$t=1$ throughout. \Omit{It does not seem to have been realised in \cite{FMTV} that
the coupling term is a total differential\comment{But how is this relevant to us?}.}

\begin{prop}\label{pr:joint flat}
The connection is integrable for any $\sfh,t\in\IC$.
\end{prop}
\begin{pf}
Consider more generally a connection on the product $X\times Y$ of two manifolds
with values in an algebra $U$, which is of the form $\nabla=\nabla^X+\nabla^Y+t\lambda$
where $\nabla^X$ is an $A$--valued connection on $X$ of the form $d+\sfh A_X$,
with $A_X$ an $A$--valued one--form on $X$, $\nabla^Y$  is an $A$--valued connection
on $Y$ of the form $d+\sfh A_Y$, with $A_Y\in\Omega^1(X,A)$, and $\lambda\in
\Omega^1(X\times Y,A)$. It is then easy to check that $\nabla$ is flat for any $t,\sfh\in
\IC$ if, and only if
\begin{gather*}
d_X A_X=0=A_X\wedge A_X\\
d_Y A_Y=0=A_Y\wedge A_Y\\
[A_X,A_Y]=0\\
d \lambda=0=\lambda\wedge\lambda\\
[A_X+A_Y,\lambda]=0
\end{gather*}

In the case at hand, $X=\gX_n$ and $A_X=\sum_{i<j}d\log(z_i-z_j)\Omega_{ij}$ are
well--known to satisfy $d_X A_X=0=A_X\wedge A_X$. Similarly, $Y=\hreg$ and
$A_Y=\half{1}\sum_{\alpha\in\Phi_+}d\log\alpha\Kalpha$ satisfy $d_Y A_Y=0=
A_Y\wedge A_Y$ \cite{MTL,TL2,FMTV}. Finally, $\lambda=d(\sum_i z_i\mu^{(i)})$
clearly satisfies $d\lambda=0=\lambda\wedge\lambda$ since it is exact and abelian.
It therefore remains to check that $[A_X+A_Y,\lambda]=0$. Write to this end 
$\lambda=\lambda_X+\lambda_Y$, where
\[\lambda_X=\sum_i dz_i\mu^{(i)}\aand
\lambda_Y=\sum_i z_id\mu^{(i)}\]
Then,\comment{the form computations below are modulo multiplicative constants}
\[\begin{split}
[A_X,\lambda_X]&=
\sum_{i<j,k}\frac{d(z_i-z_j)}{z_i-z_j}\wedge dz_k\,[\Omega_{ij},\mu^{(k)}]\\
&=
\sum_{i<j}\frac{d(z_i-z_j)}{z_i-z_j}\wedge
\Bigl(dz_i\,[\Omega_{ij},\mu^{(i)}]+dz_j\,[\Omega_{ij},\mu^{(j)}]\Bigr)\\
&=
\sum_{i<j}\frac{d(z_i-z_j)}{z_i-z_j}\wedge
(dz_i-dz_j)\,[\Omega_{ij},\mu^{(i)}]\\
&=
0
\end{split}\]
where the second equality follows from the fact that $\mu^{(k)}$ commutes with
$\Omega_{ij}$ unless $k\in\{i,j\}$ and the third from the fact that $\Omega_{ij}$
commutes with $\mu^{(i)}+\mu^{(j)}$. Next,
\[\begin{split}
[A_X,\lambda_Y]
&=
\sum_{i<j}\frac{d(z_i-z_j)}{z_i-z_j}\wedge\Bigl(z_i[\Omega_{ij},d\mu^{(i)}]+z_j[\Omega_{ij},d\mu^{(j)}]\Bigr)\\
&=
\sum_{i<j}\frac{d(z_i-z_j)}{z_i-z_j}\wedge(z_i-z_j)[\Omega_{ij},d\mu^{(i)}]\\
&=
\sum_{i<j}d(z_i-z_j)\wedge[\Omega_{ij},d\mu^{(i)}]
\end{split}\]
which completes the computation of the commutator $[A_X,\lambda]$.

To compute $[A_Y,\lambda]$, write
\[\Delta^{(n)}(\Kalpha)=\sum_{i=1}^n\Kalpha^{(i)}+2\sum_{i<j}(\Kalphap^{(ij)}+\Kalpham^{(ij)})\]
where $\Kalphapm^{(ij)}=x_{\pm\alpha}^{(i)}x_{\mp\alpha}^{(j)}$. Then,\comment{The factor of
$2$ coming from $\Delta^{(n)}(\Kalpha)=2\sum_{i<j}$ is cancelled by the factor of $1/2$ coming
from $A_Y$.}
\[\begin{split}
[A_Y,\lambda_X]
&=
\sum_{\alpha,i<j}\frac{d\alpha}{\alpha}[\Kalphap^{(ij)}+\Kalpham^{(ij)},dz_i\mu^{(i)}+dz_j\mu^{(j)}]\\
&=
\sum_{\alpha,i<j}\frac{d\alpha}{\alpha}\wedge d(z_i-z_j)[\Kalphap^{(ij)}+\Kalpham^{(ij)},\mu^{(i)}]\\
&=
\sum_{\alpha,i<j}d\alpha\wedge d(z_i-z_j)\Bigl(-\Kalphap^{(ij)}+\Kalpham^{(ij)}\Bigr)
\end{split}\]
where the first equality follows from the fact that $\Kalpha^{(i)}$ commutes with any $\mu^{(k)}$,
the second from the fact that $\Kalphap^{(ij)}+\Kalpham^{(ij)}$ is of weight zero and the third from
the fact that $[\mu^{(i)},\Kalphapm^{(ij)}]=\pm\alpha(\mu)\Kalphapm^{(ij)}$\comment{maybe change
notation for $\mu$? $\mu$ tends to denote an element in $\h^*$ rather than in $\h$.} Finally,
\[\begin{split}
[A_Y,\lambda_Y]
&=
\sum_{\alpha,i<j}\frac{d\alpha}{\alpha}[\Kalphap^{(ij)}+\Kalpham^{(ij)},z_id\mu^{(i)}+z_jd\mu^{(j)}]\\
&=
\sum_{\alpha,i<j}\frac{d\alpha}{\alpha}\wedge(z_i-z_j)[\Kalphap^{(ij)}+\Kalpham^{(ij)},d\mu^{(i)}]\\
&=
\sum_{\alpha,i<j}\frac{d\alpha}{\alpha}\wedge(z_i-z_j)d\alpha\Bigl(-\Kalphap^{(ij)}+\Kalpham^{(ij)}\Bigr)\\
&=
0
\end{split}\]
To conclude, we need to show that $[A_X,\lambda_Y]+[A_Y,\lambda_X]=0$. This follows by
writing $\Omega_{ij}=\sum_{\alpha\in\Phi_+}\Bigl(\Kalphap^{(ij)}+\Kalpham^{(ij)}\Bigr)+\Omega_{ij}^\h$,
where $\Omega_\h=t_a\otimes t^a$, with $\{t_a\},\{t^a\}$ dual bases of $\h$, so that
\[\begin{split}
[A_X,\lambda_Y]
&=
\sum_{i<j}d(z_i-z_j)\wedge[\Omega_{ij},d\mu^{(i)}]\\
&=
\sum_{i<j,\alpha}d(z_i-z_j)\wedge[\Kalphap^{(ij)}+\Kalpham^{(ij)},d\mu^{(i)}]\\
&=
\sum_{i<j,\alpha}d(z_i-z_j)\wedge d\alpha\Bigl(-\Kalphap^{(ij)}+\Kalpham^{(ij)}\Bigr)
\end{split}\]
\end{pf}

\end{document}